\documentclass[letter,11pt]{amsart}
\usepackage{amsfonts,textcomp,amssymb,amsmath,amsthm,appendix,hyperref}
\usepackage{color,ulem}
\usepackage{fullpage}
\usepackage{graphicx}
\usepackage[active]{srcltx}
\usepackage{mathabx}
\usepackage{dsfont,mathrsfs,stmaryrd,wasysym,amsbsy}
\usepackage{enumerate}

\newcommand{\be}{\begin{equation}}
\newcommand{\ee}{\end{equation}}
\newcommand{\ba}{\begin{array}}
\newcommand{\ea}{\end{array}}

\def\P{\mathbb{P}}

\def\R{{\mathbb R}}
\def\N{{\mathbb N}}

\def\P{{\mathbb P}}

\usepackage{graphicx}

\theoremstyle{plain}
\newtheorem{theorem}{Theorem}[section]

\newtheorem{lemma}[theorem]{Lemma}
\newtheorem{proposition}[theorem]{Proposition}

\theoremstyle{definition}
\newtheorem{remark}[theorem]{Remark}
\newtheorem{definition}[theorem]{Definition}

\numberwithin{equation}{section}

\begin{document}

\title[Multilevel Dyson Brownian motions]{Multilevel Dyson Brownian motions via the superposition principle}
\author{Benjamin Budway and Mykhaylo Shkolnikov}
\address{ORFE Department, Princeton University, Princeton, NJ 08544.}
\email{bb2584@princeton.edu}
\address{Department of Mathematical Sciences, Carnegie Mellon University, Pittsburgh, PA 15232.}
\email{mshkolni@gmail.com}
\footnotetext[1]{M.~Shkolnikov is partially supported by the NSF grant DMS-2108680.}

\begin{abstract}\normalsize
Multilevel Dyson Brownian motions (MDBMs) combine Dyson Brownian motions of different dimensions into a single process in a canonical way. This paper completes the theory of MDBMs for $\beta\ge2$. Specifically, we use the superposition principle of Figalli and Trevisan to construct the MDBMs for all $\beta>2$ in a unified manner. This also extends their stochastic differential equation representation, first discovered in \cite{gorin2014multilevel}, to all $\beta>2$ and proves the uniqueness of the MDBMs for all $\beta>2$. Finally, we show that their limit as $\beta\downarrow2$ is given by the $\beta=2$ MDBM of \cite{warren2005dysons}, commonly referred to as the Warren process. 
\end{abstract}

\maketitle

\section{Introduction}

Let $N\ge2$ and consider the following stochastic process in the space of $N\times N$ Hermitian matrices: As the diagonal entries, take independent real standard Brownian motions; for the entries above the diagonal, choose the real and the imaginary parts to be independent real standard Brownian motions multiplied by $\frac{1}{\sqrt{2}}$. \textsc{Dyson} has famously shown in \cite{10.1063/1.1703862} that the (real) ordered eigenvalues $X_1(t)\le X_2(t)\le\cdots\le X_N(t)$ of such a process obey the stochastic differential equation (SDE)
\begin{align}\label{DBM}
\mathrm{d}X_i(t) = \theta\sum_{j\neq i}\frac{1}{X_i(t)-X_j(t)}\,\mathrm{d}t + \mathrm{d}B_i(t),\quad i=1,\,2,\,\ldots,\,N
\end{align}
started from $X_1(0)=X_2(0)=\cdots=X_N(0)=0$, where $\theta=1$ and $B = (B_i)_{i=1}^N$ is a vector of independent real standard Brownian motions. It turns out (see \cite[Theorem 3.1]{article}) that \eqref{DBM} admits a unique strong solution for all $\theta>0$, called $\theta$-Dyson Brownian motion, or $\theta$-DBM for short. (When the value of $\theta$ is apparent from the context, we omit the dependence on $\theta$ and simply write DBM.)  

\medskip

Next, \textit{fix} a $t>0$ and consider not only the eigenvalues $X_1(t)\le X_2(t)\le\cdots\le X_N(t)$ of the full $N\times N$ Hermitian matrix, but also the eigenvalues $X^k_1(t)\le X^k_2(t)\le\cdots\le X^k_k(t)$ of its top left $k\times k$ submatrices for $k=1,\,2,\,\ldots,\,N-1$. By a well-known theorem of \textsc{Cauchy}, these eigenvalues interlace in the sense that $X_i^{k+1}(t) \le X_i^k(t) \le X_{i+1}^{k+1}(t)$ for all $1 \leq i \leq k \leq N-1$, and thus induce a probability distribution on the Gelfand-Tsetlin cone 
\begin{equation}
\mathcal{K} := \mathcal{K}_N := \big\{(x_i^k)_{1\leq i \leq k \leq N}:\, x_i^{k+1}\leq x_i^k \leq x_{i+1}^{k+1} \text{ for all } 1 \leq i \leq k \leq N-1 \big\}. 
\end{equation}
The resulting probability measure is referred to as the $\theta=1$ Hermite corners process of variance $t$. It can be characterized by two properties: The density of the random vector $(X^N_i(t))_{i=1}^N:=(X_i(t))_{i=1}^N$ is proportional to 
\begin{align}\label{Hermiteensemble}
\prod_{1 \leq i < j \leq N}(x_j^N - x_i^N)^{2\theta}\,\prod_{i=1}^N \,\exp\bigg(\!-\frac{(x_i^N)^2}{2t} \bigg), 
\end{align}
and the Markov transition density from $X^k(t):=(X_i^k(t))_{i=1}^k=x^k$ to $X^{k-1}(t)$ is given by
\begin{align}\label{DA}
\lambda^{k,\theta}(x^{k},x^{k-1}) := \frac{\Gamma(k\theta)}{\Gamma(\theta)^k}\prod_{1\leq i < j \leq k-1} (x_j^{k-1}-x_i^{k-1}) \,\prod_{i=1}^{{k-1}}\prod_{j=1}^k |x_j^k - x_i^{k-1}|^{\theta-1} \,\prod_{1\leq i < j \leq k} (x_j^k - x_i^k)^{1-2\theta},
\end{align}
the so-called Dixon-Anderson density, for $k=2,\,3,\,\ldots,\,N$. Beyond the case $\theta=1$, the described probability measure on $\mathcal{K}$ has a random matrix intepretation also for $\theta=\frac{1}{2}$ and $\theta=2$, yielding the joint distribution of the ordered eigenvalues of the full matrix and of its top left $k\times k$ submatrices for $k=1,\,2,\,\ldots,\,N-1$ in the (scaled) Gaussian Orthogonal Ensemble (GOE) and in the (scaled) Gaussian Symplectic Ensemble (GSE), respectively. For more details on these two random matrix ensembles, see, e.g., \cite[Chapters 6, 7]{Meh2004}. 

\medskip 

The above suggest the existence of a process $(X^{k,\theta}_i)_{1\le i\le k\le N}$ in $\mathcal{K}$ such that $(X^{k,\theta}_i)_{i=1}^k$ is a $\theta$-DBM for all $k=1,\,2,\,\ldots,\,N$ and that $(X^{k,\theta}_i(t))_{1\le i\le k\le N}$ is the $\theta$-Hermite corners process of variance $t$ for all $t>0$. For $\theta=1$, such a process $(X^{k,\theta}_i)_{1\le i\le k\le N}$ has been constructed by \textsc{Warren} in the seminal paper \cite{warren2005dysons}. His process is an appropriate (obliquely) reflected Brownian motion (RBM) in $\mathcal{K}$, which exhibits the two aforementioned properties when started from $0\in\mathcal{K}$. It is worth mentioning that Warren's process is different from the (non-Markovian) process given by the eigenvalues of the full matrix and its $k\times k$ top left submatrices for $k=1,\,2,\,\ldots,\,N-1$ in the Hermitian matrix-valued process described at the beginning of this introduction (see \cite{adler2012consecutive}). For $\theta>1$, the authors of~\cite{gorin2014multilevel} obtain, via a discrete-space Markov chain approximation based on Jack polynomials, a process ($Y^{mu}$ in their notation) featuring the two desired properties. For $\theta\ge2$, it is further shown in \cite{gorin2014multilevel} that their process satisfies the SDE
\begin{align}\label{MLSDE}
\mathrm{d}X_i^{k,\theta}(t) = (\theta-1)\bigg(\sum_{j=1}^{k-1} \frac{1}{X_i^{k,\theta}(t) - X_j^{k-1,\theta}(t)} - \sum_{j \neq i} \frac{1}{X_i^{k,\theta}(t)- X_j^{k,\theta}(t)}\bigg)\,\mathrm{d}t + \mathrm{d}B_i^k(t),\;1\le i\le k\le N,
\end{align}
where $(B^k_i)_{1\le i\le k\le N}$ is a vector of independent real standard Brownian motions. In \cite[Theorem~4.1]{gorin2014multilevel}, the existence of a unique $\mathcal{K}$-valued weak solution to \eqref{MLSDE} is proved for $\theta>1$ and for initial conditions in the interior of $\mathcal{K}$ and $0\in\mathcal{K}$. However, due to difficulties associated with hitting the boundary of~$\mathcal{K}$, the authors of \cite{gorin2014multilevel} identify $Y^{mu}$ with the solution to \eqref{MLSDE} only for $\theta \geq 2$.

\medskip

A natural family of initial distributions on $\mathcal{K}$ that contains the Dirac mass at $0\in\mathcal{K}$ consists of the so-called $\theta$-Gibbs distributions. These are the distributions of the form
\begin{equation}
\mu_0\times\Lambda^{N,\theta}
:=\mu_0(\mathrm{d}x^N)\,\Lambda^{N,\theta}(x^N,\mathrm{d}(x^k_i)_{1\le i\le k\le N-1})
:=\mu_0(\mathrm{d}x^N)\,\prod_{k=0}^{N-2} \lambda^{N-k,\theta}(x^{N-k},x^{N-k-1})\,\mathrm{d}x^{N-k-1}.
\end{equation}
Our first main result handles all $\theta>1$ and all $\theta$-Gibbs initial distributions (even those assigning a positive mass to the boundary of $\mathcal{K}$), in particular bridging the gap in understanding for $\theta\in(1,2)$.  

\begin{theorem}\label{MainTheorem}
Fix a $\theta>1$ and a probability measure $\mu_0$ on the Weyl chamber 
\begin{equation}
\mathcal{W}:=\mathcal{W}_N:=\{ x \in\mathbb{R}^N\!: x_1 \le x_2 \le \cdots \le x_N \}, 
\end{equation}
and for $t>0$ let $\mu_t$ be the time $t$ distribution of a $\theta$-DBM started from $\mu_0$. Then, there exists a $\mathcal{K}$-valued weak solution $X^\theta$ to the SDE \eqref{MLSDE} with the initial distribution $\mu_0 \times \Lambda^{N,\theta}$ such that 
\begin{enumerate}[(a)]
\item $X^{k,\theta}:= (X_i^{k,\theta})_{i=1}^k$ is a $\theta$-DBM for all $k=1,\,2,\,\ldots,\,N$, 
\item and for all $t\ge0$, $X^{\theta}(t)$ has the distribution $\mu_t\times\Lambda^{N,\theta}$.
\end{enumerate}
Moreover, the $\mathcal{K}$-valued weak solution satisfying (b) is unique in distribution and Markovian. 
\end{theorem}

\smallskip

The main tool we rely on to prove Theorem \ref{MainTheorem} is the superposition principle of \textsc{Trevisan} \cite{10.1214/16-EJP4453}. Given an SDE, a superposition principle guarantees that a measure-valued solution $(\nu_t)_{t\in[0,T]}$ of the associated Fokker-Planck equation can be lifted to a solution of the associated martingale problem with the fixed time distributions $(\nu_t)_{t\in[0,T]}$. The literature on superposition principles goes back to \textsc{DiPerna} and \textsc{Lions} who studied the deterministic case in \cite{DiPerna1989-nw} (see also \cite{Ambrosio2008}). \textsc{Figalli} \cite{FIGALLI2008109} was then the first to extend the result to the stochastic case, for bounded drift and diffusion coefficients. His result was, in turn, extended by \textsc{Trevisan} in \cite{10.1214/16-EJP4453} to cover drift and diffusion coefficients integrable with respect to the given fixed time distributions. A further extension can be found in \cite{bogachev2019ambrosiofigallitrevisan}, where drift and diffusion coefficients are allowed to grow faster at infinity. 

In the setting of Theorem \ref{MainTheorem}, checking that $\mu_t\times\Lambda^{N,\theta}$, $t\ge0$ solve the correct Fokker-Planck equation requires inverse moment estimates for the spacings in the $\theta$-DBM, as well as in the  Dixon-Anderson process. We hope that the new estimates for the Dixon-Anderson density can be combined with inputs about other processes, such as the $\beta$-Laguerre and the $\beta$-Jacobi processes (see, e.g., \cite{Assiotis}), to construct further multilevel processes akin to the ones in Theorem \ref{MainTheorem}.

\medskip

Given Theorem \ref{MainTheorem}, it is natural to wonder whether the processes $\{X^\theta\}_{\theta>1}$ therein converge to the Warren process as $\theta\downarrow1$. Our second contribution is to answer this question in the affirmative, at least when the limiting initial distribution $\mu_0$ does not assign mass to the boundary $\partial\mathcal{W}$ of $\mathcal{W}$. 

\begin{theorem}\label{convergence}
For each $\theta>1$, let $\mu_0^{\theta}$ be a probability measure on $\mathcal{W}$. Assume that as $\theta\downarrow1$, the measures $\mu_0^{\theta}$ converge weakly to $\mu_0$, a probability measure on $\mathcal{W}$ with $\mu_0(\partial \mathcal{W})=0$. If $X^\theta$ is the weak solution to \eqref{MLSDE} with the initial distribution $\mu_0^\theta \times \Lambda^{N,\theta}$ described in Theorem \ref{MainTheorem}, then as $\theta\downarrow1$, the processes $X^{\theta}$ converge in distribution to the Warren process $X^1$ with the initial distribution $\mu_0 \times \Lambda^{N,1}$.
\end{theorem}

\begin{remark}
Theorem \ref{convergence} along with our estimates for the Dixon-Anderson density (see Section~\ref{sec:DA}) and standard martingale problem techniques render the family $\{X^{\theta}\}_{\theta \geq 1}$ weakly continuous in $\theta$.
\end{remark}

\smallskip 

The Warren process $X^1$ is a weak solution to an SDE with reflection (SDER), whose definition we recall in Section \ref{sec:Thm2}. To identify a limit point of $\{X^\theta\}_{\theta>1}$ as $\theta\downarrow1$ with $X^1$, we use the submartingale problem formulation of SDERs due to \textsc{Kang} and \textsc{Ramanan} \cite{kang2014submartingale}. The concept of submartingale problems was originally introduced by \textsc{Stroock} and \textsc{Varadhan} for smooth domains in \cite{b87357027e204ddaa1b2245d90aab3f7}, which prompted a large body of work on the equivalence between SDERs and submartingale problems, see, e.g., \cite{51caa9f40b3b4ff69bb203303bd67a71}, \cite{Williams1987-aj}, \cite{Kwon1992TheSP}, \cite{6c0596bd-71ac-3c41-b92c-86b5281f6f20}, \cite{ca64bbce-df8d-3bcc-8136-a83b2fe7457d}. To the best of our knowledge, \cite{kang2014submartingale} is the first work to establish the equivalence between SDERs and submartingale problems in enough generality to include $X^1$. 

\medskip

Central to the construction of the Warren process in \cite{warren2005dysons} is the fact that the semigroups of $\theta=1$ DBMs of different dimensions are \textit{intertwined} through the Dixon-Anderson (conditional) density. Thus, the Warren process can be viewed as a coupling of $\theta=1$ DBMs of different dimensions that realizes their intertwining, an idea originally proposed by \textsc{Diaconis} and \textsc{Fill} in \cite{DiFi} for Markov chains. More generally, it has been proved in \cite{ramanan2016intertwinings} that, for any $\theta>0$, the semigroups of $\theta$-DBMs of different dimensions are intertwined (while analogues of this statement for the $\beta$-Laguerre and the $\beta$-Jacobi ensembles can be found in \cite{Assiotis}). In the same vein, then, the process $X^\theta$ of Theorem \ref{MainTheorem} can be regarded as a coupling that realizes the latter intertwining. 

\medskip 

The rest of the paper is structured as follows. In Section \ref{sec:DA}, we present useful results for the Dixon-Anderson density, including the aforementioned inverse moment estimates for spacings. In Section~\ref{sec:super}, we recall the superposition principle of \cite{10.1214/16-EJP4453} and show that it applies to $(\mu_t \times \Lambda^{N,\theta})_{t\ge0}$ of Theorem~\ref{MainTheorem}. In Section \ref{sec:conseq}, we rely on the results of Section \ref{sec:super} in order to prove Theorem \ref{MainTheorem}. In Section \ref{sec:Thm2}, we review the results on submartingale problems from \cite{kang2014submartingale} and combine them with Theorem \ref{MainTheorem} to prove Theorem \ref{convergence}. Throughout, we rely on various estimates for the $\theta$-DBM, the statements and proofs of which are collected in Appendix \ref{sec:app}.

\section{Preliminaries on the Dixon-Anderson density} \label{sec:DA}

We start by verifying that $\Lambda^{N,\theta}$ solves an appropriate partial differential equation (PDE) classically in the interior $\overset{\circ}{\mathcal{K}}$ of $\mathcal{K}$. To this end, we define the differential operators
\begin{eqnarray}
&& \mathcal{A}^\theta_N f = \sum_{i=1}^N c_i^{N,\theta}\,\partial_{x_i^N}f + \frac{1}{2}\sum_{i=1}^N\partial_{x_i^N x_i^N}f, \\
&& \mathcal{A}_{N-1}^{\theta,*}f = - \sum_{i=1}^{N-1} \partial_{x_i^{N-1}}(c_i^{N-1,\theta} f) + \frac{1}{2}\sum_{i=1}^{N-1} \partial_{x_i^{N-1}x_i^{N-1}}f, \\
&& \mathcal{L}^{\theta,*}_{1,...,N-1}f = - \sum_{1\leq i\leq k \leq N-1} \partial_{x_i^k}(b_i^{k,\theta} f) + \frac{1}{2}\sum_{1 \leq i \leq k \leq N-1} \partial_{x_i^k x_i^k} f, 
\end{eqnarray}
where 
\begin{equation}\label{c and b}
c_i^{N,\theta}(x):=\theta\sum_{j\neq i} \frac{1}{x^N_i-x^N_j}\quad\text{and}\quad
b_i^{k,\theta}(x) := (\theta-1)\bigg(\sum_{j=1}^{k-1}\,\frac{1}{x_i^k - x_j^{k-1}}
-\sum_{j \neq i }\,\frac{1}{x_i^k - x_j^k}\bigg).
\end{equation}

\begin{lemma}\label{DADElemma}
For any $\theta>0$, the following PDEs hold classically in $\overset{\circ}{\mathcal{K}}$:
\begin{eqnarray}
&& \mathcal{A}_{N-1}^{\theta,*}\lambda^{N,\theta} = \mathcal{A}^\theta_N \lambda^{N,\theta}, \\
&& \mathcal{L}^{\theta,*}_{1,...,N-1}\Lambda^{N,\theta} =\mathcal{A}^\theta_N \Lambda^{N,\theta}.
\end{eqnarray}
\end{lemma}

\noindent\textbf{Proof.} We fix a $\theta>0$ and omit $\theta$ in the superscripts. The two PDEs can be stated equivalently as
\begin{eqnarray}
&& \frac{1}{2}\sum_{i=1}^N \partial_{x_i^Nx_i^N}\lambda^N 
+ \sum_{i=1}^N c_i^N\,\partial_{x_i^N}\lambda^N 
= \frac{1}{2}\sum_{j=1}^{N-1}\partial_{x_j^{N-1}x_j^{N-1}}\lambda^N - \sum_{j=1}^{N-1}\partial_{x_j^{N-1}}(c_j^{N-1}\lambda^N), \label{singleinterlace} \\
&& \sum_{1 \leq i \leq k \leq N-1}\partial_{x_i^k}(b_i^k\Lambda^N) 
+ \sum_{i=1}^N\partial_{x_i^Nx_i^N} \Lambda^N + \sum_{i=1}^N c^N_i\,\partial_{x_i^N}\Lambda^N 
= \frac{1}{2} \Delta \Lambda^N. \label{DADE}
\end{eqnarray}
We start with the proof of \eqref{singleinterlace}. Let $L_{i,j} := (x_i^N - x_j^{N-1})^{-1}$ for $i=1,2,\ldots,N$ and $j=1,2,\ldots,N-1$, $B_{i,n}:=(x_i^N -x_n^N)^{-1}$ for $i=1,\,2,\,\ldots,\,N$ and $n=1,\,2,\,\ldots,\,N$, $n \neq i$, and $D_{j,m}:=(x_j^{N-1}-x_m^{N-1})^{-1}$ for $j=1,\,2,\,\ldots,\,N-1$ and $m=1,\,2,\,\ldots,\,N-1$, $m \neq j$. Then, we write the derivatives in \eqref{singleinterlace} as 
\begin{eqnarray*}
&& c_i^N\,\partial_{x_i^N}\lambda^N = \lambda^N\,\bigg(\theta(\theta-1)\Big(\sum_{n \neq i}B_{i,n}\Big)\Big(\sum_{j=1}^{N-1}L_{i,j}\Big) + \theta(1-2\theta)\Big(\sum_{n \neq i}B_{i,n}\Big)^2\bigg), \\
&& \partial_{x_i^Nx_i^N}\lambda^N = \lambda^N\,\bigg((\theta\!-\!1)\sum_{j=1}^{N-1}L_{i,j}+(1\!-\!2\theta)\sum_{n \neq i }B_{i,n}\bigg)^2
- \lambda^N\bigg((\theta\!-\!1)\sum_{j=1}^{N-1}L_{i,j}^2 + (1\!-\!2\theta) \sum_{n \neq i }B_{i,n}^2\bigg), \\
&& \partial_{x_j^{N-1}}(c_j^{N-1}\lambda^N) = \lambda^N\,\bigg(\theta\Big(\sum_{m \neq j}D_{j,m}\Big)^2 - \theta\sum_{m \neq j}D_{j,m}^2 - \theta(\theta-1)\Big(\sum_{m \neq j}D_{j,m}\Big)\Big(\sum_{i=1}^NL_{i,j}\Big)\bigg), \\
&&\partial_{x_j^{N-1}x_j^{N-1}}\lambda^N = \lambda^N\,\bigg(\sum_{m \neq j}D_{j,m} - (\theta-1)\sum_{i=1}^N L_{i,j}\bigg)^2 - \lambda^N\,\bigg(\sum_{m \neq j}D_{j,m}^2 + (\theta-1)\sum_{i=1}^N L_{i,j}^2\bigg).
\end{eqnarray*}

\smallskip

Next, following the argument on \cite[p.~254]{anderson_guionnet_zeitouni_2009}, we compute
\[
\sum_{i=1}^N \bigg(\Big(\sum_{n \neq i} B_{i,n}\Big)^2- \sum_{n \neq i} B_{i,n}^2\bigg) 
= \sum_{i\neq n\neq l\neq i} \frac{1}{x_i^N-x_n^N}\,\frac{1}{x_i^N-x_l^N} \qquad\qquad\qquad
\]
\[
\qquad\; =  \sum_{i \neq n\neq l\neq i} \frac{1}{x_l^N - x_n^N}\,\Big(\frac{1}{x_i^N - x_l^N} - \frac{1}{x_i^N - x_n^N}\Big) = -2\sum_{i \neq n\neq l\neq i} \frac{1}{x_i^N-x_n^N}\,\frac{1}{x_i^N-x_l^N}. 
\]
By comparing the second expression to the fourth one, we infer that
\[
\sum_{i=1}^N \bigg(\Big(\sum_{n \neq i} B_{i,n}\Big)^2- \sum_{n \neq i} B_{i,n}^2\bigg)=0.
\]
The same argument shows that
\[
\sum_{j=1}^{N-1}\bigg(\Big(\sum_{m \neq j} D_{j,m}\Big)^2- \sum_{m \neq j} D_{j,m}^2\bigg)=0.
\]
These allow us to reduce the two sides of \eqref{singleinterlace} to
\begin{eqnarray*}
&& \lambda^N\,\bigg(\frac{1}{2}(\theta-1)^2\sum_{i=1}^N\Big(\sum_{j=1}^{N-1}L_{i,j}\Big)^2 - (\theta-1)^2 \sum_{i=1}^N \Big(\sum_{j=1}^{N-1}L_{i,j}\Big)\Big(\sum_{n \neq i} B_{i,n}\Big) - \frac{1}{2}(\theta-1)\sum_{i=1}^N\sum_{j=1}^{N-1}L_{i,j}^2\bigg), \\
&& \lambda^N\,\bigg(\frac{1}{2}(\theta-1)^2\sum_{j=1}^{N-1}\Big(\sum_{i=1}^N L_{i,j}\Big)^2 + (\theta-1)^2\sum_{j=1}^{N-1} \Big(\sum_{m \neq j}D_{j,m}\Big)\Big(\sum_{i=1}^NL_{i,j}\Big) - \frac{1}{2}(\theta-1)\sum_{j=1}^{N-1}\sum_{i=1}^N L_{i,j}^2\bigg).
\end{eqnarray*}
Now, we expand the squares and move the cross terms, simplifying \eqref{singleinterlace} to
\[
\sum_{j=1}^{N-1}\sum_{i=1}^N\sum_{n \neq i}\Big(\frac{1}{2} L_{i,j}L_{n,j} + B_{i,n}L_{i,j}\Big) = \sum_{i=1}^N \sum_{j=1}^{N-1}\sum_{m \neq j}\Big(\frac{1}{2}L_{i,j}L_{i,m} - D_{j,m}L_{i,j} \Big).
\]
Thanks to $L_{i,j}L_{n,j}=B_{i,n}(L_{n,j}-L_{i,j})$ and $L_{i,j}L_{i,m}=D_{j,m}(L_{i,j}-L_{i,m})$, the equation reduces to
\[
\sum_{i=1}^N\sum_{n \neq i} B_{i,n}\sum_{j=1}^{N-1} (L_{i,j}+L_{n,j}) = -\sum_{j=1}^{N-1}\sum_{m \neq j} D_{j,m}\sum_{i=1}^N (L_{i,j}+L_{i,m}).
\]
Since $B_{i,n}=-B_{n,i}$ and $D_{m,j}=-D_{j,m}$, both sides are equal to $0$, and \eqref{singleinterlace} readily follows.

\medskip 

We prove \eqref{DADE} by induction over $N$. Since $\Lambda^2=\lambda^2$ and $\mathcal{L}^{\theta,*}_1=\mathcal{A}_1^{\theta,*}$, the $N=2$ case follows by equation \eqref{singleinterlace}. Now, assume equation \eqref{DADE} holds for some $N\ge 2$. We start with the first term in~\eqref{DADE} for $N+1$. Using the identities $\partial_{x_i^N}\Lambda^N = \Lambda^N(b_i^N - c_i^N)$ and $\Lambda^{N+1}=\Lambda^N\lambda^{N+1}$ we find that
\begin{equation*}
\begin{split}
&\,\sum_{1 \leq i \leq k \leq N}\partial_{x_i^k}(b_i^k \Lambda^N \lambda^{N+1}) = \lambda^{N+1}\sum_{1 \leq i \leq k \leq N-1}\partial_{x_i^k}(b_i^k \Lambda^N) 
+ \sum_{i=1}^N \partial_{x_i^N}(b_i^N \Lambda^N \lambda^{N+1}) \\
& =\lambda^{N+1}\sum_{1 \leq i \leq k \leq N-1}\partial_{x_i^k}(b_i^k \Lambda^N) 
+ \sum_{i=1}^N \partial_{x_i^N}(\lambda^{N+1}\partial_{x_i^N} \Lambda^N) 
+ \sum_{i=1}^N \partial_{x_i^N}(\lambda^{N+1}c_i^N\Lambda^N) \\
& = \lambda^{N+1}\Big(\sum_{1 \leq i \leq k \leq N-1}\partial_{x_i^k}(b_i^k \Lambda^N) + \sum_{i=1}^N\partial_{x_i^Nx_i^N}\Lambda^N + 
\sum_{i=1}^N c_i^N \partial_{x_i^N}\Lambda^N\Big) 
+ \nabla_{x^N} \lambda^{N+1}\cdot \nabla_{x^N} \Lambda^N \\
&\qquad\qquad\qquad\qquad\qquad\qquad\qquad\qquad\qquad\qquad\qquad\qquad\qquad\quad\;\;\;
+ \Lambda^N\sum_{i=1 }^N\partial_{x_i^N}(c_i^N \lambda^{N+1}) \\
& = \frac{1}{2}\lambda^{N+1}\Delta\Lambda^N + \nabla_{x^N} \lambda^{N+1}\cdot \nabla_{x^N} \Lambda^N+\Lambda^N\sum_{i=1 }^N\partial_{x_i^N}(c_i^N \lambda^{N+1}) \\
& = \frac{1}{2}\Delta\Lambda^{N+1} + \Lambda^N \sum_{i=1 }^N\partial_{x_i^N}(c_i^N \lambda^{N+1}) - \frac{1}{2}\Lambda^N \sum_{i=1}^{N+1} \partial_{x_i^{N+1}x_i^{N+1}}\lambda^{N+1} 
- \frac{1}{2}\Lambda^N\sum_{j=1}^N \partial_{x_j^Nx_j^N}\lambda^{N+1}.
\end{split}
\end{equation*}
Adding the other two terms on the left-hand side of \eqref{DADE} for $N+1$ we obtain
\begin{equation*}
\begin{split}
& \,\frac{1}{2}\Delta\Lambda^{N+1} + \Lambda^N \Big(\frac{1}{2}\sum_{i=1}^{N+1}\partial_{x_i^{N+1}x_i^{N+1}}\lambda^{N+1}
\!+\! \sum_{i=1}^{N+1}c_i^{N+1} \partial_{x_i^{N+1}}\lambda^{N+1} 
\!+\! \sum_{i=1}^N \partial_{x_i^N}(c_i^N \lambda^{N+1}) 
\!-\! \frac{1}{2}\sum_{j=1}^N \partial_{x_j^Nx_j^N}\lambda^{N+1}\Big) \\
& = \frac{1}{2}\Delta\Lambda^{N+1},
\end{split}
\end{equation*}
where we have used equation \eqref{singleinterlace} for $N+1$. \qed

\medskip

The following crucial lemma, of interest in its own right, allows to control the inverses of spacings between points on consecutive levels in terms of those on the same level. 

\begin{lemma}\label{diagtoflat}
For all $k\in\N$ and $x^{k+1}\in\mathcal{W}_{k+1}$, let 
\[
\mathcal{R}(x^{k+1}):=\{z \in \mathbb{R}^k:\,x^{k+1}_i \le z_i \le x^{k+1}_{i+1},\,i=1,\,2,\,\ldots,\,k\}.
\]
Then, for all $0 \leq \alpha,\alpha_1,\alpha_2 < \theta$ and $a,b\in\{1,\,2,\,\ldots,\,k\}$, $a\neq b$:
\begin{equation*}
\begin{split}
& \,\int_{\mathcal{R}(x^{k+1})} (x_{a+1}^{k+1}-x_a^k)^{-\alpha}\,\lambda^{k+1}(x^{k+1},x^k)\,\mathrm{d}x^k 
\le \frac{C_1(k,\alpha,\theta)}{\theta-\alpha}\,(x_{a+1}^{k+1}-x_a^{k+1})^{-\alpha}, \\
& \,\int_{\mathcal{R}(x^{k+1})} (x_{a+1}^{k+1}-x_a^k)^{-\alpha_1}\,(x_{b+1}^{k+1} - x_b^k)^{-\alpha_2}\,\lambda^{k+1}(x^{k+1},x^k)\,\mathrm{d}x^k \\
& \le \frac{C_2(k,\alpha_1,\alpha_2,\theta)}{(\theta-\alpha_1)(\theta-\alpha_2)}\, 
(x_{a+1}^{k+1} - x_a^{k+1})^{-\alpha_1}\,(x_{b+1}^{k+1} - x_b^{k+1})^{-\alpha_2},
\end{split}
\end{equation*}  
where $C_1(k,\alpha,\theta) := \frac{\Gamma(\theta-\alpha+1)\,\Gamma((k+1)\theta)}{\Gamma(\theta)\,\Gamma((k+1)\theta-\alpha)}$, $C_2(k,\alpha_1,\alpha_2,\theta):= \frac{\Gamma(\theta-\alpha_1+1)\,\Gamma(\theta-\alpha_2+1)\,\Gamma((k+1)\theta)}{\Gamma(\theta)^2\,\Gamma((k+1)\theta-\alpha_1-\alpha_2)}$.
\end{lemma}

\noindent\textbf{Proof.} We only prove the second inequality, as the first can then be recovered by setting $\alpha_2=0$ in the second. The Dixon-Anderson identity (cf.~ \cite[Proposition 4.2.1]{8a873437-53d6-348d-9b7a-19c7327502a1}) states that for $s_1,\,s_2,\,\ldots,\,s_{k+1}>0$:
\begin{equation}
\int_{\mathcal{R}(x^{k+1})} \prod_{i=1}^k\prod_{j=1}^{k+1} |z_i - x^{k+1}_j|^{s_j-1} 
\prod_{1\leq i < n \leq k}(z_n-z_i)\,\mathrm{d}z 
= \frac{\prod_{j=1}^{k+1}\Gamma(s_j)}{\Gamma(\sum_{j=1}^{k+1}s_j)}\prod_{1 \leq j < m \leq k+1}(x^{k+1}_m-x^{k+1}_j)^{s_m+s_j-1}.
\end{equation}
Moreover, denoting by $\delta_{\cdot,\cdot}$ the Kronecker delta, we can write 
\begin{equation*}
\begin{split}
&\,\int_{\mathcal{R}(x^{k+1})} \!(x_{a+1}^{k+1}\!-\! x_a^k)^{-\alpha_1}\,(x_{b+1}^{k+1} \!-\! x_b^k)^{-\alpha_2}\,\lambda^{k+1}(x^{k+1},x^k)\,\mathrm{d}x^k =\frac{\Gamma((k+1)\theta)}{\Gamma(\theta)^{k+1}}
\prod_{1\leq j < m \leq k+1}\! (x_m^{k+1} \!-\! x_j^{k+1})^{1-2\theta} \\
& \cdot\!\int_{\mathcal{R}(x^{k+1})}\!\prod_{1\leq i < n \leq k} (x_n^k\!-\! x_i^k) \prod_{i=1}^{{k}}\prod_{j=1}^{k+1} |x_j^{k+1} \!-\! x_i^{k}|^{\theta-1-\alpha_1 \delta_{j,a+1}-\alpha_2\delta_{j,b+1}}\prod_{\overset{i=1}{i\neq a}}^k |x_{a+1}^{k+1}\!-\! x_i^k|^{\alpha_1} 
\prod_{\overset{i=1}{i\neq b}}^k |x_{b+1}^{k+1}\!-\! x_i^k|^{\alpha_2}\,\mathrm{d}x^k.
\end{split}
\end{equation*}
Thanks to the interlacing inequalities, $0 \le x_{a+1}^{k+1}-x_i^k \le x_{a+1}^{k+1}-x_{i}^{k+1} $, $i=1,\,2,\,\ldots,\,a-1$, and $0 \le x_i^k - x_{a+1}^{k+1}\le x_{i+1}^{k+1}-x_{a+1}^{k+1}$, $i=a+1,\,a+2,\,\ldots,\,k$. The same inequalities hold also with $a$ replaced by $b$. Therefore, the above can be upper bounded by
\begin{equation}\label{interineq}
\begin{split}
& \frac{\Gamma((k+1)\theta)}{\Gamma(\theta)^{k+1}}\prod_{1\leq j < m \leq k+1} (x_m^{k+1} - x_j^{k+1})^{1-2\theta} \prod_{j=1,j \neq a, j \neq a+1}^{k+1} |x_{a+1}^{k+1}-x_j^{k+1}|^{\alpha_1}\prod_{j=1, j \neq b, j \neq b+1}^{k+1}|x_{b+1}^{k+1}-x_j^{k+1}|^{\alpha_2} \\
& \qquad\qquad\qquad\qquad\qquad\qquad\;\,
\cdot  \int_{\mathcal{R}(x^{k+1})}\prod_{1\leq i < n \leq k} (x_n^k-x_i^k)\, \prod_{i=1}^{{k}}\prod_{j=1}^{k+1} |x_j^{k+1} - x_i^{k}|^{\theta-1-\alpha_1 \delta_{j,a+1}-\alpha_2\delta_{j,b+1}}\,\mathrm{d}x^k \\
& = \frac{\Gamma((k+1)\theta)}{\Gamma(\theta)^{k+1}}
\prod_{1\leq j < m \leq k+1} (x_m^{k+1} \!-\! x_j^{k+1})^{1-2\theta} 
\prod_{j=1, j \neq a, j \neq a+1}^{k+1}|x_{a+1}^{k+1}\!-\! x_j^{k+1}|^{\alpha_1}
\prod_{j=1, j \neq b, j \neq b+1}^{k+1}|x_{b+1}^{k+1}\!-\! x_j^{k+1}|^{\alpha_2}
\end{split}
\end{equation}
\begin{equation*}
\begin{split}
\qquad\quad\;\;\,
\cdot\frac{\Gamma(\theta)^{k-1}\,\Gamma(\theta-\alpha_1)\,\Gamma(\theta-\alpha_2)}
{\Gamma((k+1)\theta-\alpha_1-\alpha_2)} 
\prod_{1\leq j < m \leq k+1}(x_m^{k+1}-x_j^{k+1})
^{2\theta -\alpha_1(\delta_{j,a+1}+\delta_{m,a+1})-\alpha_2(\delta_{j,b+1}+\delta_{m,b+1})-1} \\
=\frac{\Gamma(\theta-\alpha_1)\,\Gamma(\theta-\alpha_2)\,\Gamma((k+1)\theta)}
{\Gamma(\theta)^2\,\Gamma((k+1)\theta - \alpha_1-\alpha_2)}
\,(x_{a+1}^{k+1}-x_{a}^{k+1})^{-\alpha_1}
\,(x_{b+1}^{k+1}-x_{b}^{k+1})^{-\alpha_2}. \qquad\qquad\qquad\qquad\qquad\qquad\;
\end{split}
\end{equation*}
To get the promised constant, we use that $(\theta-\alpha_i)\,\Gamma(\theta-\alpha_i)=\Gamma(\theta-\alpha_i+1)$, $i=1,2$. \qed

\medskip

The next lemma lets us decouple the interactions in $\Lambda^N$ near $\partial\mathcal{K}$.~For $k=1,\,2,\,\ldots,\,N$ and $1\le a\le b\le k$, it uses the notation $x^k_{a,...,b}$ for the vector $(x^k_i)_{i=a}^b$.  

\begin{lemma}\label{decouple}
Fix any $\theta\ge1$ and pick an $\varepsilon \in (0,1]$. For $x^{k+1} \in \mathcal{W}_{k+1}$ such that $x_i^{k+1} + 2 \varepsilon \leq x_{i+1}^{k+1}$, $i=1,\,2,\,\ldots,\,k$, define 
\[
\mathcal{R}_{\varepsilon}(x^{k+1}) = \{x^k \in \mathcal{R}(x^{k+1}):\,x_i^{k+1} + 2\varepsilon \le x_i^k + \varepsilon \le x_{i+1}^{k+1},\,i=1,\,2,\,\ldots,\,k\}. 
\]
Then, for all $x^k \in \mathcal{R}_{\varepsilon}(x^{k+1})$ satisfying $x_a^k+\varepsilon= x_{a+1}^{k+1}$, we have that 
\[
\lambda^{k+1}(x^{k+1},x^k) \le C \varepsilon^{\theta -1}(x_{a+1}^{k+1}-x_a^{k+1})^{-\theta}\,\lambda^a(x_{1,...,a}^{k+1},x_{1,...,a-1}^k)
\,\lambda^{k-a+1}(x_{a+1,...,k+1}^{k+1},x_{a+1,...,k}^k).
 \]
Here, $C$ depends on $\theta,k$ and $a$ but not on $\varepsilon$, and $\lambda_1$ is understood to be identically equal to $1$.
\end{lemma}

\noindent\textbf{Proof.} We start by noting that
\[
\prod_{\stackrel{i=1}{i \neq a}}^k |x_i^k - x_a^k| 
=\prod_{i=1}^{a-1}(x_{a+1}^{k+1}-\varepsilon-x_i^k)
\,\prod_{i=a+1}^k (x_i^k - x_{a+1}^{k+1}+\varepsilon).
\]
Since $x^k \in \mathcal{R}_{\varepsilon}(x^{k+1})$, it holds $x_{i}^k - x_{a+1}^{k+1}+ \varepsilon \leq 2(x_i^k - x_{a+1}^{k+1})$ for $i=a+1,\,a+2,\,\ldots,\,k$. Therefore, by the same argument used to get the upper bound \eqref{interineq}, the above is at most
\begin{align}\label{in1}
2^{k-a}\prod_{\stackrel{i=1}{i \neq a}}^k |x_i^k - x_{a+1}^{k+1}| \leq 2^{k-a} \prod_{j=1,j\neq a, j \neq a+1}^{k+1}|x_j^{k+1}-x_{a+1}^{k+1}|.
\end{align}
In addition, the requirement $x^k \in \mathcal{R}(x^{k+1})$ implies that
\begin{align}\label{in2}
\prod_{i=1}^{a-1}\prod_{n=a+1}^k (x_n^k-x_i^k) 
\le \prod_{j=1}^{a-1}\prod_{m=a+2}^{k+1} (x_m^{k+1}-x_j^{k+1}).
\end{align}
Putting \eqref{in1} and \eqref{in2} together, we infer that
\begin{equation}\label{foineq}
\begin{split}
& \prod_{1 \leq i < n \leq k} (x_n^k - x_i^k) \leq 2^{k-a} \prod_{1 \leq i < n \leq a-1}(x_n^k - x_i^k)\,\prod_{a+1 \leq i < n \leq k} (x_n^k - x_i^k) \\
& \qquad\qquad\qquad\qquad\;\;
\cdot\prod_{j=1}^{a-1}|x_j^{k+1}-x_{a+1}^{k+1}|\,
\prod_{j=a+2}^{k+1} |x_j^{k+1}-x_{a+1}^{k+1}|
\,\prod_{j=1}^{a-1}\prod_{m=a+2}^{k+1}(x_m^{k+1}-x_j^{k+1}) \\
& \leq  2^{k-a}  (x_{a+1}^{k+1}-x_a^{k+1})^{-1}
\prod_{1 \leq i < n \leq a-1} (x_n^k - x_i^k)\,
\prod_{a+1 \leq i < n \leq k}(x_n^k - x_i^k)\,
\prod_{j=1}^a \prod_{m=a+1}^{k+1} (x_m^{k+1}-x_j^{k+1}),
\end{split}
\end{equation}
where we applied $|x_j^{k+1}-x_{a+1}^{k+1}|\le x_j^{k+1}-x_{a}^{k+1}$, $j=a+2,\,a+3,\,\ldots,\,k+1$ in the last step. 

\medskip

For the double product in \eqref{DA}, we follow the derivation of \eqref{in1}, and then employ $x^{k+1} \in \mathcal{W}_{k+1}$:
\[
\prod_{j=1}^{k+1}|x_a^k - x_j^{k+1}|^{\theta-1} 
\le \varepsilon^{\theta -1}2^{(\theta-1)(k-a)}
\prod_{\stackrel{j=1}{j\neq a+1}}^{k+1}|x_{a+1}^{k+1}-x_j^{k+1}|^{\theta-1} 
\qquad\qquad\qquad\qquad\qquad\qquad\qquad\qquad\;\;
\]
\begin{eqnarray}
\qquad\qquad\qquad\quad\;\;\;
&\le&\!\!\!\! \varepsilon^{\theta -1}2^{(\theta-1)(k-a)}\prod_{j=1}^a |x_{a+1}^{k+1}-x_j^{k+1}|^{\theta -1}
\prod_{j=a+2}^{k+1}|x_{a}^{k+1}-x_j^{k+1}|^{\theta -1} \label{2nddiagineq} \\
\qquad\qquad\qquad\quad\;\;\;
&=&\!\!\!\! \varepsilon^{\theta -1}2^{(\theta-1)(k-a)}(x_{a+1}^{k+1}-x_a^{k+1})^{-\theta+1}\prod_{j=1}^a |x_{a+1}^{k+1}-x_j^{k+1}|^{\theta -1}\prod_{j=a+1}^{k+1}|x_{a}^{k+1}-x_j^{k+1}|^{\theta -1}.
\nonumber
\end{eqnarray}
Due to $x^k \in \mathcal{R}(x^{k+1})$, we have for $i=1,\,2,\,\ldots,\,a-1$:
\[
\prod_{j=1}^{k+1}|x_j^{k+1}-x_i^k|^{\theta-1} \leq \prod_{j=1}^a|x_j^{k+1}-x_i^k|^{\theta-1} \prod_{j=a+1}^{k+1}|x_j^{k+1}-x_i^{k+1}|^{\theta -1}.
\]
Therefore, it holds
\begin{equation}\label{2ndlowineq}
\begin{split}
\prod_{i=1}^{a-1}\prod_{j=1}^{k+1}|x_j^{k+1}-x_i^k|^{\theta-1}
& \le \prod_{i=1}^{a-1}\prod_{j=1}^a|x_j^{k+1}-x_i^k|^{\theta-1}\, \prod_{m=1}^{a-1}\prod_{j=a+1}^{k+1} |x_j^{k+1}-x_m^{k+1}|^{\theta-1} \\
& = \prod_{i=1}^{a-1}\prod_{j=1}^a|x_j^{k+1}\!-\! x_i^k|^{\theta-1} \,\prod_{m=1}^{a}\prod_{j=a+1}^{k+1} |x_j^{k+1}\!-\! x_m^{k+1}|^{\theta-1}  
\,\Big(\prod_{j=a+1}^{k+1} |x_j^{k+1}\!-\! x_a^{k+1}|^{\theta-1} \Big)^{-1}.
\end{split}
\end{equation}
A similar argument shows that
\begin{equation}\label{3rdlowineq}
\begin{split}
\prod_{i=a+1}^{k}\prod_{j=1}^{k+1}|x_{j}^{k+1}-x_i^{k}|^{\theta-1} 
\le \prod_{i=a+1}^k\prod_{j=a+1}^{k+1}| x_{j}^{k+1}-x_i^k|^{\theta-1}\prod_{m=a+2}^{k+1}\prod_{j=1}^a|x_{j}^{k+1}-x_{m}^{k+1}|^{\theta-1}
\qquad\qquad\qquad\qquad\quad \\
= \prod_{i=a+1}^k\prod_{j=a+1}^{k+1}| x_{j}^{k+1}\!-\! x_i^k|^{\theta-1}
\prod_{m=a+1}^{k+1}\prod_{j=1}^a|x_{j}^{k+1}\!-\! x_{m}^{k+1}|^{\theta-1}
\Big( \prod_{j=1}^a |x_{j}^{k+1}\!-\! x_{a+1}^{k+1}|^{\theta-1} \Big)^{-1}.
\end{split}
\end{equation}

\smallskip

Putting the inequalities \eqref{2nddiagineq}, \eqref{2ndlowineq} and \eqref{3rdlowineq} together we arrive at
\begin{equation}\label{2ndineq}
\begin{split}
\prod_{i=1}^k\prod_{j=1}^{k+1} |x_i^{k}-x_j^{k+1}|^{\theta-1} = \prod_{j=1}^{k+1}|x_a^k - x_j^{k+1}|^{\theta-1} \prod_{i=1}^{a-1}\prod_{j=1}^{k+1}|x_i^k - x_j^{k+1}|^{\theta-1}\prod_{i=a+1}^k\prod_{j=1}^{k+1} |x_i^k - x_j^{k+1}|^{\theta-1} \\
\leq  \varepsilon^{\theta -1}2^{(\theta-1)(k-a)}(x_{a+1}^{k+1}-x_a^{k+1})^{-\theta+1}
\prod_{j=1}^a |x_{a+1}^{k+1}-x_j^{k+1}|^{\theta -1}
\prod_{j=a+1}^{k+1}|x_{a}^{k+1}-x_j^{k+1}|^{\theta -1} \quad \\
\cdot\prod_{i=1}^{a-1}\prod_{j=1}^a|x_j^{k+1}-x_i^k|^{\theta-1} \prod_{m=1}^{a}\prod_{j=a+1}^{k+1} |x_j^{k+1}-x_m^{k+1}|^{\theta-1} 
\Big(\prod_{j=a+1}^{k+1} |x_j^{k+1}-x_a^{k+1}|^{\theta-1} \Big)^{-1} \;\;\, \\
\cdot\prod_{i=a+1}^k\prod_{j=a+1}^{k+1}| x_{j}^{k+1}-x_i^k|^{\theta-1} \prod_{m=a+1}^{k+1}\prod_{j=1}^a|x_{j}^{k+1}-x_{m}^{k+1}|^{\theta-1}
\Big( \prod_{j=1}^a |x_{j}^{k+1}-x_{a+1}^{k+1}|^{\theta-1} \Big)^{-1} \\
=\varepsilon^{\theta -1}2^{(\theta-1)(k-a)}(x_{a+1}^{k+1}-x_a^{k+1})^{-\theta+1} \prod_{m=1}^a\prod_{j=a+1}^{k+1} |x_{j}^{k+1}-x_{m}^{k+1}|^{2\theta-2} \qquad\qquad\qquad\quad\;\;\; \\
\cdot\prod_{i=1}^{a-1}\prod_{j=1}^a|x_j^{k+1}-x_i^k|^{\theta-1} \prod_{i=a+1}^k\prod_{j=a+1}^{k+1}| x_{j}^{k+1}-x_i^k|^{\theta-1}.
\qquad\qquad\qquad\qquad\qquad\qquad\quad\;
\end{split}
\end{equation}

\smallskip

Thanks to the identity
\begin{equation}\label{3rdineq}
\begin{split}
\prod_{1\leq j < m \leq k+1} (x_m^{k+1}-x_j^{k+1})^{1-2\theta} = &\,\prod_{j=1}^a\prod_{m=a+1}^{k+1} |x_m^{k+1}-x_j^{k+1}|^{1-2\theta} \\
&\cdot\prod_{1 \leq j < m \leq a} (x_m^{k+1}-x_j^{k+1})^{1-2\theta}
\prod_{a+1 \leq j < m \leq k+1} (x_m^{k+1}-x_j^{k+1})^{1-2\theta},
\end{split}
\end{equation}
it remains to combine \eqref{foineq}, \eqref{2ndineq} and \eqref{3rdineq} in order to obtain the lemma. \qed

\medskip

Finally, we recall the intertwining between $\theta$-DBMs of different dimensions established in \cite[Theorem 1.1, Remark 1.2]{ramanan2016intertwinings}.

\begin{proposition} \label{intertwinings}
Fix any $\theta>0$. For $k\in\N$, let $(P_t^k)_{t\ge0}$ be the transition semigroup of a $\theta$-DBM of dimension $k$, and let $L^{k+1}$ be the transition operator associated with $\lambda^{k+1}$. Then, 
\[
\Big(\prod_{r=k+1}^N L^r\Big)\,P^k_t = P_t^N\,\Big(\prod_{r=k+1}^N L^r\Big),\quad t\ge0,\quad N\ge k\ge 1,
\]
where both sides are understood as transition operators taking bounded measurable functions on $\mathcal{W}_k$ to bounded measurable functions on $\mathcal{W}_N$.
\end{proposition}

\section{Superposition principle} \label{sec:super}

Take any $\theta>1$, $N\in\N$ and $T>0$. Let $\mu_0$ be a probability measure on $\mathcal{W}_N$ and let $(\mu_t)_{t\in[0,T]}$ be the fixed time distributions of a $\theta$-DBM started from $\mu_0$. The main goal of this section is to prove that the curve $(\nu_t)_{t \in [0,T]} := (\mu_t \times \Lambda^{N,\theta})_{t \in [0,T]}$  of probability measures can be lifted to a solution of the martingale problem associated with the SDE \eqref{MLSDE}, an instance of the ``superposition principle''. We start by giving a brief summary of the  superposition principle in the form of \cite{10.1214/16-EJP4453}. 

\subsection{Summary of the superposition principle}

A curve $(\nu_t)_{t \in [0,T]}$ of probability measures on~$\R^d$ is called Borel if $t \rightarrow \nu_t(A)$ is Borel measurable for every Borel $A \subseteq \mathbb{R}^d$. Given any measurable $b\!:[0,T]\times\R^d\to\R^d$ and symmetric positive semidefinite $a\!:[0,T]\times\R^d\to\R^{d\times d}$, let 
\begin{equation}
\mathcal{L} := \mathcal{L}_{a,b,t} := \frac{1}{2}\sum_{i,j=1}^d a_{ij}(t,\cdot)\,\partial_{x_i x_j}
+\sum_{i=1}^d b_i(t,\cdot)\,\partial_{x_i}.
\end{equation}
Write $\mathcal{A}_c$ for the space of $C^{1,2}$ functions compactly supported in $(0,T)\times\mathbb{R}^d$, and $\mathcal{A}$ for the space of bounded, continuous functions on $[0,T]\times\mathbb{R}^d$ with bounded first derivative in $t$ and bounded first and second derivatives in $x$ on $(0,T)\times \mathbb{R}^d$.

\begin{definition}\label{PDE}
A Borel $(\nu_t)_{t \in [0,T]}$ is a weak solution of the Fokker-Planck equation $\partial_t\nu_t = \mathcal{L}^* \nu_t$ if 
\begin{align}\label{FPint}
\int_0^T \int_{\mathbb{R}^d} |a(t,x)|+|b(t,x)|\,\nu_t(\mathrm{d}x)\,\mathrm{d}t < \infty,
\end{align}
and for every $f \in \mathcal{A}_c$, 
\begin{align}\label{FP}
\int_0^T \int_{\mathbb{R}^d} \partial_t f(t,x) + \mathcal{L} f(t,x)\,\nu_t(\mathrm{d}x)\,\mathrm{d}t = 0.
\end{align}
\end{definition}

\begin{definition}
Let $\eta$ be a probability measure on $C([0,T],\mathbb{R}^d)$, and let $e_t$, $t\in[0,T]$ denote the coordinate mappings. Then, $\eta$ is said to solve the martingale problem associated with $a$ and $b$ if
\begin{equation}
\int_{C([0,T],\mathbb{R}^d)} \int_0^T  |a|(t,\cdot) \circ e_t + |b|(t,\cdot) \circ e_t \,\mathrm{d}t\,\mathrm{d}\eta < \infty,
\end{equation}
and for each $f \in \mathcal{A}$, the following process is a martingale under $\eta$ with respect to the natural filtration:
\[
f(t,\cdot) \circ e_t - \int_0^t \partial_s f(s,\cdot) \circ e_s + \mathcal{L}f(s,\cdot)\circ e_s \,\mathrm{d}s,\quad t\in[0,T].
\]
 \end{definition}

\smallskip

The following ``superposition principle'' can be found in \cite[Theorem 2.5]{10.1214/16-EJP4453}.

\begin{proposition}\label{tsuper}
If a Borel weakly continuous $(\nu_t)_{t \in [0,T]}$ is a weak solution of the Fokker-Planck equation $\partial_t\nu_t = \mathcal{L}^* \nu_t$, then there exists a solution $\eta$ of the corresponding martingale problem such that $\eta\circ e_t^{-1} = \nu_t$, $t \in [0,T]$.
\end{proposition}

\subsection{Superposition principle in our context} \label{super context}

In this subsection, we verify \eqref{FPint} and \eqref{FP} for $(\nu_t)_{t \in [0,T]} := (\mu_t \times \Lambda^{N,\theta})_{t \in [0,T]}$ and $\mathcal{L}:=\mathcal{L}_{1,...,N} := \frac{1}{2}\Delta + b \cdot\nabla$, where 
\[
b_i^k(x):=\sum_{j=1}^{k-1} \frac{\theta -1 }{x_i^k - x_j^{k-1}} - \sum_{j \neq i} \frac{\theta -1}{x_i^k - x_j^k},\quad 1\le i\le k\le N. 
\]
The proof is broken into six steps. In Step 1, we show \eqref{FPint}. In Steps 2, 3 and 4, we establish \eqref{FP} for a restricted class of initial distributions $\mu_0$ and of test functions $f$. In Steps 5 and 6, we extend~\eqref{FP} to all $f \in \mathcal{A}_c$ and to all probability measures $\mu_0$ on $\mathcal{W}_N$, respectively. 

\medskip

\noindent\textbf{Step 1.} In order to show \eqref{FPint}, it suffices to check 
\[
\int_0^T \int_{\mathcal{K}}|b_i^k(x)|\,\nu_t(\mathrm{d}x)\,\mathrm{d}t < \infty
\]
for all $1 \leq i \leq k \leq N$, $k \ge 2$. With $C$ being a finite positive constant that  changes from line to line, 
\[ 
|b_i^k(x)| = \bigg|\sum_{j=1}^{k-1} \frac{\theta -1 }{x_i^k - x_j^{k-1}} - \sum_{j \neq i} \frac{\theta -1}{x_i^k - x_j^k}\bigg| \leq C\bigg(\frac{1}{x_i^k - x_{i-1}^{k-1}} + \frac{1}{x_i^{k-1}-x_i^k}\bigg).
\]
For $x^N \in \mathcal{W}_N$ and $m=1,\,2,\,\ldots,\,N-1$, we define
\begin{eqnarray*}
&& \mathcal{K}(x^N) = \big\{(x_i^k)_{1 \leq i \leq k \leq N-1}:\, x_i^{k+1} \le x_i^{k} \le x_{i+1}^{k+1},\,1 \leq i \leq k \leq N-1 \big\}, \\
&& \mathcal{K}^m(x^N) = \big\{(x_i^k)_{1 \leq i \leq k,\,m\leq k \leq N-1}:\, x_i^{k+1} \leq x_i^{k} \leq x_{i+1}^{k+1},\,1 \leq i \leq k,\, m \leq k \leq N-1 \big\}. 
\end{eqnarray*}
Then, by Lemma \ref{diagtoflat}, 
\begin{eqnarray} 
&& \int_{\mathcal{K}} |b_i^k(x)|\,\nu_t(\mathrm{d}x) 
= \int_{\mathcal{W}_N} \int_{\mathcal{K}(x^N)} |b_i^k(x)|\,\Lambda^{N,\theta}(x)\,
\mathrm{d}x^1\,\ldots\,\mathrm{d}x^{N-1}\,\mu_t(\mathrm{d}x^N) \nonumber \\
&& \qquad\;\;\;
\le C\int_{\mathcal{W}_N} \int_{\mathcal{K}^k(x^N)}\bigg(\frac{1}{x_i^k - x_{i-1}^{k}} + \frac{1}{x_{i+1}^{k}-x_i^k}\bigg)\,\prod_{m=k+1}^N \lambda^{m,\theta}(x^m,x^{m-1})\, \mathrm{d}x^k\,\ldots\,\mathrm{d}x^{N-1}\,\mu_t(\mathrm{d}x^N) \nonumber \\
&& \qquad\;\;\;
= C\,\mathbb{E}\bigg[\frac{1}{X_i^k(t) - X_{i-1}^{k}(t)} + \frac{1}{X_{i+1}^{k}(t)-X_i^k(t)}\bigg], \label{int1}
\end{eqnarray}
where $X^k(t)$ is obtained by sampling first a $\theta$-DBM $X^N$ with the initial distribution~$\mu_0$, and then $X^k(t)$ according to the Dixon-Anderson density. By Lemma \ref{intertwinings}, a copy of $X^k(t)$ can be obtained by sampling first $X^N(0)$ according to $\mu_0$, then $X^k(0)$ according to the Dixon-Anderson density, and finally $X^k(t)$ according to the transition operator of the $\theta$-DBM. Thus, the integral from $0$ to $T$ of the upper bound \eqref{int1} is finite by Lemma \ref{DBMint}(a).

\medskip  

\noindent\textbf{Step 2.} Now, assume additionally that $\mu_0(\mathrm{d}x^N)=p_0(x^N)\,\mathrm{d}x^N$ for some $p_0\in L^1(\mathcal{W}_N)$ satisfying 
\[
\mathrm{supp}(p_0)\subset\,\stackrel{\circ}{\mathcal{W}}_N \text{compactly and }
\int_{\mathcal{W}_N} |\log p_0(x^N)|\,p_0(x^N)\,\mathrm{d}x^N < \infty.
\]
Then, by Lemma \ref{reg}, $\mu_t$ admits a density $p_t\in W^1_1(\mathcal{W}_N)$ for $t\in(0,T]$. Also, let $f\in C^{1,2}_c((0,T)\times\mathcal{K})$ be such that for some $\varepsilon_0>0$, $f(t,x)=0$ if $\mathrm{dist}(x^N,\partial{\mathcal{W}_N})<\varepsilon_0$. Further, define the operators
\begin{eqnarray}
&& \mathcal{L}_{1,...,N-1} = \frac{1}{2}\sum_{1\leq i \leq k \leq N-1} \partial_{x_i^kx_i^k}
+ \sum_{1\leq i \leq k \leq N-1} b_i^k\,\partial_{x_i^k}, \\
&& \mathcal{L}_{N} = \frac{1}{2}\sum_{i=1}^N \partial_{x_i^N x_i^N}
+ \sum_{i=1}^N b_i^N\,\partial_{x_i^N}.
\end{eqnarray}
In Steps 2, 3 and 4, we study respectively
\begin{eqnarray}
&& \int_0^T \int_{\mathcal{K}} \partial_t f(t,x)\,p_t(x^N)\,\Lambda^{N,\theta}(x)\,\mathrm{d}x\,\mathrm{d}t, \label{timederiv} \\
&& \int_0^T \int_{\mathcal{K}} (\mathcal{L}_Nf)(t,x)\,p_t(x^N)\,\Lambda^{N,\theta}(x)\,\mathrm{d}x\,\mathrm{d}t, \label{topderiv} \\
&& \int_0^T \int_{\mathcal{K}} (\mathcal{L}_{1,...,N-1}f)(t,x)\,p_t(x^N)\,\Lambda^{N,\theta}(x)\,\mathrm{d}x\,\mathrm{d}t. \label{bottomderiv}
\end{eqnarray}

\smallskip

While dealing with \eqref{timederiv} and \eqref{topderiv}, we drop the superscript of $N$ in $b^N$ and $c^N$. The function
\[
h(t,x^N):= \int_{\mathcal{K}(x^N)} f(t,x)\,\Lambda^{N,\theta}(x)\,\mathrm{d}x^1\,\ldots\,\mathrm{d}x^{N-1}
\]
belongs to $C_c^{1,1}((0,T)\times\stackrel{\circ}{\mathcal{W}}_N)$ by Lemma \ref{diagtoflat}, the identity $\nabla_{x^N}\Lambda^{N,\theta} = (b - c)\Lambda^{N,\theta}$, and the boundedness of $\partial_t f$ and $\nabla_{x^N}f$. Therefore, by Lemma \ref{TIBP},
\begin{equation*}
\begin{split}
\int_0^T\int_{\mathcal{K}} \partial_t f(t,x)\,p_t(x^N)\,\Lambda^{N,\theta}(x)\,\mathrm{d}x\, \mathrm{d}t 
& = \int_0^T \int_{\mathcal{W}_N} \partial_t h(t,x^N)\, p_t(x^N)\, \mathrm{d}x^N \,\mathrm{d}t \\
& = \int_0^T \int_{\mathcal{W}_N} -c \cdot (\nabla_{x^N}h) \, p_t 
+ \frac{1}{2}\nabla_{x^N}h \cdot \nabla_{x^N} p_t \, \mathrm{d}x^N \,\mathrm{d}t.
\end{split}
\end{equation*}
Finally, $|\nabla_{x^N} \Lambda^{N,\theta}|$ is integrable over $\mathcal{K}(x^N)$ for any $x^N \in \stackrel{\circ}{\mathcal{W}}_N$ by Lemma \ref{diagtoflat}. This, the boundedness of $f$ and $|\nabla_{x^N} f|$, and $\Lambda^{N,\theta}(x^N,\cdot)\equiv0$ on the codimension one part of $\partial \mathcal{K}(x^N)$ yield
\[ 
\nabla_{x^N}\int_{\mathcal{K}(x^N)} f(t,x)\,\Lambda^{N,\theta}(x)\,\mathrm{d}x^1\,\ldots\,\mathrm{d}x^{N-1}
= \int_{\mathcal{K}(x^N)} \nabla_{x^N}(f \Lambda^{N,\theta})\,\mathrm{d}x^1\,\ldots\,\mathrm{d}x^{N-1}.
\]
All in all, the term in \eqref{timederiv} equals to
\begin{align}\label{ftimederiv}
\int_{0}^T \int_{\mathcal{K}} -c \cdot \nabla_{x^N}(f\Lambda^{N,\theta})\, p_t + \frac{1}{2}\nabla_{x^N}(f\Lambda^{N,\theta})\cdot\nabla_{x^N}p_t \, \mathrm{d}x \, \mathrm{d}t .
\end{align}

\smallskip 

\noindent\textbf{Step 3.} We turn to \eqref{topderiv}. Since $\nabla_{x^N}\Lambda^{N,\theta} = (b - c)\Lambda^{N,\theta}$, we have that
\[ 
\int_0^T \int_{\mathcal{K}} (\mathcal{L}_Nf)\, p_t \, \Lambda^{N,\theta}\,\mathrm{d}x \,\mathrm{d}t 
= \int_0^T \int_{\mathcal{K}} \big(c \cdot \nabla_{x^N}f + \frac{1}{2}\Delta_{x^N} f\big)\, p_t \, \Lambda^{N,\theta}\,\mathrm{d}x \,\mathrm{d}t 
+ \int_0^T \int_{\mathcal{K}} (\nabla_{x^N} f \cdot \nabla_{x^N} \Lambda^{N,\theta})\, p_t \,\mathrm{d}x \,\mathrm{d}t.
\]
We seek to integrate by parts. To this end, for $\varepsilon>0$, we let 
\[
\mathcal{K}_{\varepsilon} := \{x \in \mathcal{K}:\, x_i^N + 2\varepsilon \le x_i^{N-1} +\varepsilon \le x_{i+1}^{N},\,i=1,\,2,\,\ldots,\,N-1\}, 
\]
and let $\mathcal{K}^m_{\varepsilon}$ be the $(x_i^k)_{1 \leq i \leq k,\, m \leq k \leq N}$ projection of $\mathcal{K}_{\varepsilon}$. Since $\Delta_{x^N}f$ is bounded, 
\[
\lim_{\varepsilon\downarrow0}\,\int_0^T \int_{\mathcal{K} \setminus \mathcal{K}_{\varepsilon}} |\Delta_{x^N}f|\, p_t \,\Lambda^{N,\theta}\,\mathrm{d}x \,\mathrm{d}t=0.
\]
Write $\kappa(x) = (\kappa_i^{k}(x))_{1 \leq i \leq k \leq N}$ for the outward normal to $ \partial\mathcal{K}_{\varepsilon}$ at $x \in \partial \mathcal{K}_{\varepsilon}$, and set 
\begin{eqnarray*}
&& D_{i,+,\varepsilon} = \{x\in \partial \mathcal{K}_{\varepsilon}:\, x_{i-1}^{N-1} + \varepsilon = x_i^N\},\quad i=2,\,3,\,\ldots,\,N, \\ 
&& D_{i,-,\varepsilon} = \{ x\in \partial \mathcal{K}_{\varepsilon}:\, x_i^{N-1} - \varepsilon = x_i^N\},\quad i =1,\,2,\,\ldots,\,N-1. 
\end{eqnarray*}
Then, $\kappa_i^N\neq 0$ only on $D_{i,+,\varepsilon}$ and $D_{i,-,\varepsilon}$. Integrating by parts, we find that
\begin{equation}\notag
\begin{split}
\int_0^T\int_{\mathcal{K}_{\varepsilon}} \partial_{x_i^Nx_i^N}f \, p_t \,\Lambda^{N,\theta}\,\mathrm{d}x \,\mathrm{d}t = 
&\,-\int_0^T\int_{\mathcal{K}_{\varepsilon}} \partial_{x_i^N}f\,\partial_{x_i^N}(p_t\Lambda^{N,\theta})\,\mathrm{d}x \,\mathrm{d}t \\
&\,+\int_0^T\int_{D_{i,+,\varepsilon}} \partial_{x_i^N}f \, p_t \,\Lambda^{N,\theta}\, \kappa_i^N\,\mathrm{d}S \,\mathrm{d}t
+\int_0^T\int_{D_{i,-,\varepsilon}} \partial_{x_i^N}f \, p_t \,\Lambda^{N,\theta}\,\kappa_i^N \,\mathrm{d}S \,\mathrm{d}t.
\end{split}
\end{equation}

\smallskip

We only consider the $D_{i,+,\varepsilon}$ term, as the $D_{i,-,\varepsilon}$ term can be studied similarly. Let $0<t_1 < t_2<T$ be such that $f(t,x)=0$ unless $t \in [t_1,t_2]$. Estimating $\partial_{x_i^N}f$ and $\kappa_i^N$ by constants, integrating out $x^1,\,\ldots,\,x^{N-2}$, using Lemma \ref{decouple}, and integrating out $x^{N-1}_1,\,\ldots,\,x^{N-1}_{i-2},\,x^{N-1}_i,\,\ldots,\,x^{N-1}_{N-1}$ we infer that  
\begin{equation*}
\begin{split}
\bigg |\int_0^T\int_{D_{i,+,\varepsilon}} \partial_{x_i^N}f\,p_t\,\Lambda^{N,\theta}\, \kappa_i^N \,\mathrm{d}S \,\mathrm{d}t\bigg| 
& \le C\varepsilon^{\theta-1}\int_{t_1}^{t_2}\int_{\mathcal{W}_N}(x_i^N-x_{i-1}^N)^{-\theta}\, p_t(x^N)\, \mathrm{d}x^N \,\mathrm{d}t \\
& = C \varepsilon^{\theta-1} \int_{t_1}^{t_2} \mathbb{E}\big[(X_i^N(t) - X_{i-1}^N(t))^{-\theta}\big] \,\mathrm{d}t,
\end{split}
\end{equation*}
where $X^N$ is a $\theta$-DBM with the initial distribution $\mu_0$. The latter integral is finite by Lemma \ref{DBMint}(c), and so the overall expression goes to $0$ with $\varepsilon$.

\medskip

Finally, we argue that $\partial_{x_i^N}f\,\partial_{x_i^N}(p_t\Lambda^{N,\theta})$ is integrable over $(0,T)\times\mathcal{K}$. Bounding $\partial_{x_i^N}f$ by a constant and applying Lemma \ref{diagtoflat} we see that, on the one hand,
\begin{equation*}
\begin{split}
&\int_0^T\int_{\mathcal{K}} |\partial_{x_i^N}f\,\partial_{x_i^N}\Lambda^{N,\theta}|\, p_t \,\mathrm{d}x \,\mathrm{d}t 
\le C \int_{t_1}^{t_2} \int_{\mathcal{W}_N} \int_{\mathcal{R}(x^N)} |b_i^N-c_i^N|\,\lambda^{N,\theta}\, p_t \,\mathrm{d}x^{N-1} \,\mathrm{d}x^{N} \,\mathrm{d}t \\
&\qquad\qquad\qquad\qquad\quad
\le C \int_{t_1}^{t_2} \int_{\mathcal{W}_N} \int_{\mathcal{R}(x^N)} \big((x_i^N - x_{i-1}^{N-1})^{-1} + (x_{i}^{N-1}-x_i^N)^{-1} \big) \,
\lambda^{N,\theta} \, p_t \,\mathrm{d}x^{N-1} \,\mathrm{d}x^N \,\mathrm{d}t \\
&\qquad\qquad\qquad\qquad\quad
\le C\int_{t_1}^{t_2} \int_{\mathcal{W}_N} 
\big((x_i^N - x_{i-1}^N)^{-1} + (x_{i+1}^N - x_i^N)^{-1}\big)\, p_t \,\mathrm{d}x^N \,\mathrm{d}t \\ 
&\qquad\qquad\qquad\qquad\quad
= C\int_{t_1}^{t_2} \mathbb{E}\big[(X_i^N(t)-X_{i-1}^N(t))^{-1} + (X_{i+1}^N(t) - X_i^N(t))^{-1} \big] \,\mathrm{d}t,
\end{split}
\end{equation*}
where the last integral is finite by Lemma \ref{DBMint}(c). On the other hand, Hölder's inequality implies
\begin{equation*}
\begin{split}
\int_{0}^T\int_{\mathcal{K}}|\partial_{x_i^N}f \,\partial_{x_i^N}p_t|\,\Lambda^{N,\theta}\,\mathrm{d}x \,\mathrm{d}t  
& \le C \int_0^T \int_{\mathcal{W}_N}|\partial_{x_i^N}p_t| \,\mathrm{d}x^N \,\mathrm{d}t \\
& \le C \bigg(\int_0^T \int_{\mathcal{W}_N}p_t \,\mathrm{d}x^N \,\mathrm{d}t \bigg)^{1/2}\,
\bigg(\int_0^T \int_{\mathcal{W}_N}\frac{|\nabla_{x^N}p_t|^2}{p_t} \,\mathrm{d}x^N \,\mathrm{d}t \bigg)^{1/2},
\end{split}
\end{equation*}
where the last double integral is finite by Lemma \ref{reg}. All in all, we have shown that
\[
\int_0^T \int_{\mathcal{K}} (\Delta_{x^N}f) \, p_t \,\Lambda^{N,\theta} \,\mathrm{d}x \,\mathrm{d}t = - \int_{0}^T \int_{\mathcal{K}}\nabla_{x^N}f \cdot \nabla_{x^N}(p_t\Lambda^{N,\theta}) \,\mathrm{d}x \,\mathrm{d}t.
\]
Thus, the term in \eqref{topderiv} equals to
\begin{align}\label{ftopderiv}
\int_0^T \int_{\mathcal{K}} (c \cdot \nabla_{x^N}f)\, p_t \,\Lambda^{N,\theta} 
- \frac{1}{2}\nabla_{x^N}f \cdot \nabla_{x^N}(p_t\Lambda^{N,\theta}) 
+ (\nabla_{x^N}f \cdot \nabla_{x^N} \Lambda^{N,\theta})\, p_t \, \mathrm{d}x \,\mathrm{d}t.
\end{align}

\smallskip

\noindent\textbf{Step 4.} In this (most involved) step, we deal with the term in \eqref{bottomderiv}. By Lemma \ref{DADElemma}, 
\[
\int_{\mathcal{K}_{\varepsilon}} \sum_{1\leq i \leq k \leq N-1} 
\Big(-fp_t \, \partial_{x_i^k}(b_i^k \Lambda^{N,\theta}) 
+ \frac{1}{2}fp_t \,\partial_{x_i^k x^k_i} \Lambda^{N,\theta}\Big)\,\mathrm{d}x 
= \int_{\mathcal{K}_{\varepsilon}} f p_t \, c^N \cdot \nabla_{x^N} \Lambda^{N,\theta} 
+ \frac{1}{2}f p_t\,\Delta_{x^N}\Lambda^{N,\theta}\,\mathrm{d}x.
\]
Now, we integrate by parts on both sides. Recalling that $p_t$ is a function of $x^N$ only, we find that
\begin{equation}\label{miracle}
\begin{split}
& \,\int_{\mathcal{K}_{\varepsilon}}  \sum_{1\leq i \leq k \leq N-1} \Big(b_i^k\Lambda^{N,\theta}
\, p_t \,\partial_{x_i^k}f  + \frac{1}{2}\Lambda^{N,\theta}\, p_t \,\partial_{x_i^kx^k_i} f \Big)\,\mathrm{d}x \\
& = \int_{\mathcal{K}_{\varepsilon}} fp_t \, c^N \cdot \nabla_{x^N} \Lambda^{N,\theta} 
- \frac{1}{2}\nabla_{x^N} (f p_t) \cdot \nabla_{x^N}\Lambda^{N,\theta}\,\mathrm{d}x 
+ R_1(t) + R_2(t),
\end{split}
\end{equation}
where the boundary terms are given by 
\begin{equation*}
\begin{split}
& R_1(t) = \int_{\partial \mathcal{K}_{\varepsilon}} fp_t\,\Big(\sum_{1\leq i \leq k \leq N-1} \Big(b_i^k  \Lambda^{N,\theta} - \frac{1}{2} \partial_{x_i^k}\Lambda^{N,\theta}\Big)\,\kappa_i^k
+ \frac{1}{2} \sum_{i=1}^N  \partial_{x_i^N} \Lambda^{N,\theta}\,\kappa_i^N\Big)\,\mathrm{d}S, \\
& R_2(t) =  \frac{1}{2}\int_{\partial \mathcal{K}_{\varepsilon}} \sum_{1 \leq i \leq k \leq N-1} p_t \Lambda^{N,\theta}\,\partial_{x_i^k}f \,\kappa_i^k \,\mathrm{d}S,
\end{split}
\end{equation*}
and $\kappa_i^k$ is the $(k,i)$ component of the outward normal to $\partial\mathcal{K}_\varepsilon$.

\medskip

We first aim to show that the integrals of $|R_2(t)|$ and of $|R_1(t)|$ with respect to $t$ tend to $0$ with~$\varepsilon$. Note that $\partial\mathcal{K}_\varepsilon\subset\bigcup _{1\leq i \leq k \leq N-1} (D_{k,i,+,\varepsilon}\cup D_{k,i,-,\varepsilon})$, where  
\[
D_{k,i,+,\varepsilon}:=\{x\in\partial\mathcal{K}_{\varepsilon}:\, x_i^k + \varepsilon = x_{i+1}^{k+1} \}
\quad\text{and}\quad
D_{k,i,-,\varepsilon}:=\{x\in\partial\mathcal{K}_{\varepsilon}:\, x_i^k - \varepsilon = x_i^{k+1}\}. 
\]
Also, for $m=1,\,2,\,\ldots,\, N$, let $\mathcal{K}^m$ be the $(x_i^k)_{1\le i \le k \le N,\, k \ge m}$ projection of $\mathcal{K}$.
Bounding $\partial_{x_j^m}f$ and $\kappa^m_j$ by constants, integrating out $x^1,\,\ldots,\, x^{k-1}$, using Lemma \ref{decouple}, and finally integrating out $x^k_1,\,\ldots,\,x^k_{i-1},\,x^k_{i+1},\,\ldots,x^k_k$ we see that
\begin{equation}\label{decoupleexample}
\begin{split}
&\,\int_{D_{k,i,+,\varepsilon}} \sum_{1 \leq j \leq m \leq N-1} |p_t\Lambda^{N,\theta}\, \partial_{x_j^m}f\,\kappa_j^m| \,\mathrm{d}S \\
&\le C\,\mathbf{1}_{[t_1,t_2]}(t)\, \varepsilon^{\theta-1}\int_{\mathcal{K}^{k+1}}(x_{i+1}^{k+1}-x_i^{k+1})^{-\theta} \prod_{m=k+2}^{N}\lambda^{m,\theta}(x^m,x^{m-1})\, p_t(x^N)\,\mathrm{d}x^{k+1}\,\ldots\,\mathrm{d}x^N \\
&= C\,\mathbf{1}_{[t_1,t_2]}(t)\,\varepsilon^{\theta-1}\, \mathbb{E}\big[(X_{i+1}^{k+1}(t)-X_i^{k+1}(t))^{-\theta} \big],
\end{split}
\end{equation}
where $X^{k+1}$ is a $\theta$-DBM with the initial density
\[ 
\int_{\mathcal{K}^{k+2}} \mu_0(\mathrm{d}x^N)\,\prod_{m=k+2}^N \lambda^{m,\theta}(x^m,x^{m-1})\, \mathrm{d}x^{k+2}\,\ldots\,\mathrm{d}x^{N-1}, 
\]
and we have made use of Lemma \ref{intertwinings}. A similar argument shows that
\[
\int_{D_{k,i,-,\varepsilon}} \sum_{1 \leq j \leq m \leq N-1} |p_t\Lambda^{N,\theta}\, \partial_{x_j^m}f \,\kappa_j^m| \,\mathrm{d}S  
\le C\,\mathbf{1}_{[t_1,t_2]}(t)\,\varepsilon^{\theta-1}\, \mathbb{E}\big[(X_{i}^{k+1}(t)-X_{i-1}^{k+1}(t))^{-\theta} \big].
\]
In conclusion,
\[
\int_{0}^T |R_2(t)| \,\mathrm{d}t \le C\varepsilon^{\theta-1}\int_{t_1}^{t_2}\sum_{1 \leq i \leq k \leq N-1} \mathbb{E}\big[(X_{i+1}^{k+1}(t)-X_{i}^{k+1}(t))^{-\theta}\big] \,\mathrm{d}t,
\]
which goes to $0$ with $\varepsilon$ thanks to Lemma \ref{DBMint}(c). 

\medskip

We turn to the more complicated $R_1(t)$. Consider $D_{k,i,+,\varepsilon}$ for fixed $1\le i\le k\le N-2$. On $D_{k,i,+,\varepsilon}$, it holds $\kappa_i^k \equiv -\kappa_{i+1}^{k+1}$, and $\kappa^m_j\equiv0$ for all other $(m,j)$. Therefore,
\begin{equation*}
\begin{split}
&\,\bigg|\int_{D_{k,i,+,\varepsilon}} fp_t\,\Big(\sum_{1\leq j \leq m \leq N-1} \Big(b_j^m  \Lambda^{N,\theta} - \frac{1}{2} \partial_{x_j^m}\Lambda^{N,\theta}\Big)\,\kappa_j^m
+ \frac{1}{2} \sum_{j=1}^N  \partial_{x_j^N} \Lambda^{N,\theta}\,\kappa_j^N\Big)\,\mathrm{d}S\bigg| \\
& \le C\,\mathbf{1}_{[t_1,t_2]}(t) \int_{D_{k,i,+,\varepsilon}} p_t\,\Big|b_i^k\Lambda^{N,\theta} - \frac{1}{2}\partial_{x_i^k}\Lambda^{N,\theta} - b_{i+1}^{k+1}\Lambda^{N,\theta} 
+ \frac{1}{2}\partial_{x_{i+1}^{k+1}}\Lambda^{N,\theta}\Big|\,\mathrm{d}S.
\end{split}
\end{equation*}
Recalling the definition of $b^k_i$ in \eqref{c and b}, as well as 
\[
\partial_{x_i^k}\Lambda^{N,\theta}(x) 
= \Lambda^{N,\theta}(x)\,\bigg((\theta-1)\sum_{j=1}^{k-1}\frac{1}{x_i^k - x_j^{k-1}}+2(1-\theta)\sum_{\stackrel{j=1}{j \neq i}}^k\frac{1}{x_i^k - x_j^k}+ (\theta-1)\sum_{j=1}^{k+1}\frac{1}{x_i^k - x_j^{k+1}}\bigg),
\]
we note that the terms involving $(x_{i+1}^{k+1}-x_i^k)^{-1}$ in $-\frac{1}{2}\partial_{x_i^k}\Lambda^{N,\theta}$, $-b_{i+1}^{k+1}\Lambda^{N,\theta}$, and $\frac{1}{2}\partial_{x_{i+1}^{k+1}}\Lambda^{N,\theta}$ cancel out exactly! Thus, we have the estimate (with obvious adjustments if $i=1$ or $i=k$):
\begin{equation}\label{BAD}
\begin{split}
& \,\mathbf{1}_{[t_1,t_2]}(t) \int_{D_{k,i,+,\varepsilon}} p_t\,\Big|b_i^k\Lambda^{N,\theta} - \frac{1}{2}\partial_{x_i^k}\Lambda^{N,\theta} - b_{i+1}^{k+1}\Lambda^{N,\theta} 
+ \frac{1}{2}\partial_{x_{i+1}^{k+1}}\Lambda^{N,\theta}\Big|\,\mathrm{d}S \\
& \le C\,\mathbf{1}_{[t_1,t_2]}(t) \int_{D_{k,i,+,\varepsilon}} p_t(x^N)\,\Lambda^{N,\theta}(x)
\bigg(\frac{1}{x_i^{k-1}-x_i^k}+\frac{1}{x_i^k-x_{i-1}^{k-1}} + \frac{1}{x_i^k - x_i^{k+1}} \\
& \qquad\qquad\qquad\qquad\qquad\qquad\qquad\qquad\;\;\;
+\frac{1}{x_{i+1}^{k}-x_{i+1}^{k+1}} 
+\frac{1}{x_{i+1}^{k+1}-x_{i+1}^{k+2}} 
+\frac{1}{x_{i+2}^{k+2}-x_{i+1}^{k+1}}\bigg)\,\mathrm{d}S.
\end{split}
\end{equation}

\smallskip

We deal with the terms in \eqref{BAD} one by one in the order of their appearance. Letting $D^m_{k,i,+,\varepsilon}$ be the respective $(x^n_j)_{1\le j\le n\le N,\,n\ge m}$ projections of $D_{k,i,+,\varepsilon}$ for $m=1,\,2,\,\ldots,\,k$, integrating out $x^1,\,\ldots,\,x^{k-2}$, and using Lemma \ref{diagtoflat} we obtain
\begin{equation*}
\begin{split}
&\,\mathbf{1}_{[t_1,t_2]}(t)\int_{D_{k,i,+,\varepsilon}^{k-1}} \frac{p_t(x^N)}{x_i^{k-1}-x_i^k} \,\prod_{m=k}^N \lambda^{m,\theta}(x^m,x^{m-1})\,\mathrm{d}S \\
&\le C\,\mathbf{1}_{[t_1,t_2]}(t)\int_{D_{k,i,+,\varepsilon}^k} \frac{p_t(x^N)}{x_{i+1}^k - x_i^k}\, \prod_{m=k+1}^N \lambda^{m,\theta}(x^m,x^{m-1})\,\mathrm{d}S \\
&\le C\,\mathbf{1}_{[t_1,t_2]}(t)\int_{D_{k,i,+,\varepsilon}^k} \frac{p_t(x^N)}{x_{i+1}^k - x_{i+1}^{k+1}}\,\prod_{m=k+1}^N \lambda^{m,\theta}(x^m,x^{m-1})\,\mathrm{d}S.
\end{split}
\end{equation*}
Note that the following argument also takes care of the fourth term in \eqref{BAD}. Applying Lemma \ref{decouple} as in the proof of inequality \eqref{decoupleexample}, and then Lemma \ref{diagtoflat} again, we find the upper bound
\begin{equation}\notag
\begin{split}
&\,C\,\mathbf{1}_{[t_1,t_2]}(t)\,\varepsilon^{\theta-1}\int_{\mathcal{K}^{k+1}} p_t(x^N)\, (x_{i+1}^{k+1}-x_i^{k+1})^{-\theta}(x_{i+2}^{k+1} - x_{i+1}^{k+1})^{-1}\prod_{m=k+2}^N \lambda^{m,\theta}(x^m,x^{m-1})\,\mathrm{d}x^{k+1}\,\ldots\,\mathrm{d}x^N \\
&= C\,\mathbf{1}_{[t_1,t_2]}(t)\,\varepsilon^{\theta-1}\, \mathbb{E}\big[(X_{i+1}^{k+1}(t)-X_{i}^{k+1}(t))^{-\theta}(X_{i+2}^{k+1}(t)-X_{i+1}^{k+1}(t))^{-1}\big].
\end{split}
\end{equation}
We have used Lemma \ref{intertwinings} to claim the latter equality. Integrating in $t$ and employing Lemma \ref{DBMint}(c), we see that the contribution from the first term in \eqref{BAD} goes to $0$ with $\varepsilon$. The second term in \eqref{BAD} can be upper bounded in the same manner by 
\[
C\, \mathbf{1}_{[t_1,t_2]}(t)\,\varepsilon^{\theta-1}\,
\mathbb{E}\big[(X_{i+1}^{k+1}(t)-X_{i}^{k+1}(t))^{-\theta}(X_{i}^{k+1}(t)-X_{i-1}^{k+1}(t))^{-1}\big].
\]

\smallskip

For the third term in \eqref{BAD}, we apply $(x_i^k - x_i^{k+1})^{-1}\leq 2(x_{i+1}^{k+1}-x_i^{k+1})^{-1}$, $x \in D_{k,i,+,\varepsilon}$, and then Lemma \ref{decouple}  as in the proof of inequality \eqref{decoupleexample} to obtain the upper bound 
\[
C\,\mathbf{1}_{[t_1,t_2]}(t)\,\varepsilon^{\theta-1}\,\int_{\mathcal{K}^{k+1}} p_t(x^N)\, (x_{i+1}^{k+1}-x_i^{k+1})^{-\theta -1}\prod_{m=k+2}^N \lambda^{m,\theta}(x^m,x^{m-1})\, \mathrm{d}x^{k+1}\,\ldots\,\mathrm{d}x^N.
\]
In view of Lemma \ref{intertwinings}, this equals to 
\[
C\,\mathbf{1}_{[t_1,t_2]}(t)\,\varepsilon^{\theta-1}\, \mathbb{E}\big[(X_{i+1}^{k+1}(t)-X_i^{k+1}(t))^{-\theta-1}\big].
\]
Integrating in $t$ and recalling Lemma \ref{DBMint}(c) we see that the contribution from the third term in~\eqref{BAD} goes to $0$ with $\varepsilon$. 

\medskip

For the fifth term in \eqref{BAD}, we rely on Lemma \ref{decouple} to get
\begin{equation}\notag
\begin{split}
&\,\mathbf{1}_{[t_1,t_2]}(t) \int_{D_{k,i,+,\varepsilon}} p_t(x^N)\,\Lambda^{N,\theta}(x)\, (x_{i+1}^{k+1}-x_{i+1}^{k+2})^{-1}\,\mathrm{d}S \\
&\leq C\mathbf{1}_{[t_1,t_2]}(t)\varepsilon^{\theta-1} \int_{\mathcal{K}_{\varepsilon}^{k+1}} p_t(x^N)(x_{i+1}^{k+1}\!-\! x_i^{k+1})^{-\theta}(x_{i+1}^{k+1}\!-\! x_{i+1}^{k+2})^{-1} \prod_{m=k+2}^N \lambda^{m,\theta} (x^m,x^{m-1})\,\mathrm{d}x^{k+1}\,\ldots\,\mathrm{d}x^N \\
&\leq C\mathbf{1}_{[t_1,t_2]}(t)\varepsilon^{\theta-1} \int_{\mathcal{K}_{\varepsilon}^{k+1}} p_t(x^N)(x_{i+1}^{k+2}\!-\! x_i^{k+1})^{-\theta}(x_{i+1}^{k+1}\!-\! x_{i+1}^{k+2})^{-1} \prod_{m=k+2}^N \lambda^{m,\theta}(x^m,x^{m-1}) \,\mathrm{d}x^{k+1}\,\ldots\,\mathrm{d}x^N.
\end{split}
\end{equation}
Now, let $\alpha\in(0,\theta-1)$ and estimate the latter expression from above by
\[
C \mathbf{1}_{[t_1,t_2]}(t)\varepsilon^{\theta-1-\alpha}\int_{\mathcal{K}_{\varepsilon}^{k+1}} \! p_t(x^N)(x_{i+1}^{k+2}-x_i^{k+1})^{-\theta+\alpha}(x_{i+1}^{k+1}-x_{i+1}^{k+2})^{-1} \prod_{m=k+2}^N \!\lambda^{m,\theta}(x^m,x^{m-1})\,\mathrm{d}x^{k+1}\ldots\,\mathrm{d}x^N.
\]
Thanks to Lemmas \ref{diagtoflat} and \ref{intertwinings}, this can be bounded further by 
\[
C\,\mathbf{1}_{[t_1,t_2]}(t)\,\varepsilon^{\theta-1-\alpha}\, \mathbb{E}\big[(X_{i+1}^{k+2}(t)-X_{i}^{k+2}(t))^{-\theta+\alpha}(X_{i+2}^{k+2}(t)-X_{i+1}^{k+2}(t))^{-1}\big].
\]
Integrating in $t$ and using Lemma \ref{DBMint}(c) we infer that the contribution from the fifth term in \eqref{BAD} goes to $0$ with $\varepsilon$. A similar argument also takes care of the sixth term in \eqref{BAD}. Moreover, the~$D_{k,i,-,\varepsilon}$ terms can be dealt with analogously. 

\medskip

In the case $k = N-1$, we note that
\begin{equation}\notag
\begin{split}
&\,\bigg| \sum_{i=1}^{N-1}\int_{D_{N-1,i,+,\varepsilon}} fp_t\,\Big(\sum_{1\leq j \leq m \leq N-1} \Big(b_j^m  \Lambda^{N,\theta} - \frac{1}{2} \partial_{x_j^m}\Lambda^{N,\theta}\Big)\,\kappa_j^m
+ \frac{1}{2} \sum_{j=1}^N  \partial_{x_j^N} \Lambda^{N,\theta}\,\kappa_j^N\Big)\,\mathrm{d}S\bigg| \\
&\leq C\,\mathbf{1}_{[t_1,t_2]}(t)\sum_{i=1}^{N-1}\int_{D_{N-1,i,+,\varepsilon}} p_t\,\Big|\Big(b_i^{N-1} \Lambda^{N,\theta} - \frac{1}{2}\partial_{x_i^{N-1}}\Lambda^{N,\theta}\Big)\,\kappa_i^{N-1} + \frac{1}{2}\partial_{x_{i+1}^N}\Lambda^{N,\theta}\,\kappa_{i+1}^{N}\Big|\,\mathrm{d}S \\
&\leq C\,\mathbf{1}_{[t_1,t_2]}(t)\sum_{i=1}^{N-1}\int_{D_{N-1,i,+,\varepsilon}} p_t \,\Big|b_i^{N-1} \Lambda^{N,\theta} - \frac{1}{2}\partial_{x_i^{N-1}}\Lambda^{N,\theta} - \frac{1}{2}\partial_{x_{i+1}^N}\Lambda^{N,\theta}\Big| \,\mathrm{d}S.
\end{split}
\end{equation}
Again, the crucial singularity $(x_{i+1}^N - x_i^{N-1})^{-1}$ cancels out on $D_{N-1,i,+\varepsilon}$, resulting in the upper bound
\[
C\,\mathbf{1}_{[t_1,t_2]}(t)\int_{D_{N-1,i,+,\varepsilon}} p_t(x^N)\,\Lambda^{N,\theta}(x) \bigg(\frac{1}{x_i^{N-1}\!-\! x_{i-1}^{N-2}}+\frac{1}{x_{i}^{N-2}\!-\! x_{i}^{N-1}} + \frac{1}{x_{i}^{N-1}\!-\! x_i^{N}} + \frac{1}{x_{i+1}^{N-1}\!-\! x_{i+1}^N}\bigg)\,\mathrm{d}S.
\]
It remains to proceed as with the first four terms in \eqref{BAD}. 

\medskip

We proved that the contributions of the remainder terms in \eqref{miracle} go to $0$ with $\varepsilon$. In order to get
\begin{equation}\label{IBP}
\begin{split}
&\,\int_{0}^T\int_{\mathcal{K}}  \sum_{1\leq i \leq k \leq N-1} \Big(b_i^k \,\Lambda^{N,\theta} p_t \,\partial_{x_i^k}f  + \frac{1}{2}\Lambda^{N,\theta}\,p_t\,\partial_{x_i^kx_i^k}f \Big)\,\mathrm{d}x \,\mathrm{d}t \\
& = \int_{0}^T\int_{\mathcal{K}} f p_t \, c^N \cdot \nabla_{x^N} \Lambda^{N,\theta} - \frac{1}{2}\nabla_{x^N} (f p_t) \cdot \nabla_{x^N}\Lambda^{N,\theta}\,\mathrm{d}x \,\mathrm{d}t,
\end{split}
\end{equation}
it remains to check the integrability on the right-hand side. (Note that the integrability on the left-hand side follows from the boundedness of $\partial_{x_i^k}f$, $\partial_{x_i^kx_i^k}f$, and Step 1.) For the first term, $\nabla_{x^N}\Lambda^{N,\theta} = (b^N - c^N)\Lambda^{N,\theta}$,
Lemma \ref{diagtoflat}, and the Cauchy-Schwarz inequality yield
\begin{equation}\notag
\begin{split}
\int_{\mathcal{K}}  |f| \,|c^N \cdot \nabla_{x^N} \Lambda^{N,\theta}| \, p_t \,\mathrm{d}x 
&\le C \,\mathbf{1}_{[t_1,t_2]}(t) \sum_{i=1}^N \bigg(\int_{\mathcal{K}} p_t \,|c_i^N|\,|b_i^N|\, \Lambda^{N,\theta} \,\mathrm{d}x +  \int_{\mathcal{K}} p_t \,|c_i^N|^2\,\Lambda^{N,\theta} \,\mathrm{d}x\bigg) \\
&\leq C \,\mathbf{1}_{[t_1,t_2]}(t)\sum_{i=1}^N \, \int_{\mathcal{W}_N} p_t(x^N)\, \bigg(\frac{1}{(x_i^N - x_{i-1}^N)^2} + \frac{1}{(x_{i+1}^N - x_{i}^N)^2}\bigg)\,\mathrm{d}x^N,
\end{split}
\end{equation}
which is integrable in $t$ by Lemma \ref{DBMint}(c). 

\medskip

For the second term on the right-hand side of \eqref{IBP}, we use Lemmas \ref{diagtoflat} and \ref{DBMint}(c) to find
\[
\sum_{i=1}^N \,\int_0^T \int_{\mathcal{K}} p_t\,|\partial_{x_i^N}f\,\partial_{x_i^N }\Lambda^{N,\theta}| \,\mathrm{d}x \,\mathrm{d}t
\leq C\sum_{i=1}^N \,\int_{t_1}^{t_2} \int_{\mathcal{W}_N} p_t(x^N)\,\bigg(\frac{1}{x_i^N - x_{i-1}^N} + \frac{1}{x_{i+1}^N - x_{i}^N} \bigg) \,\mathrm{d}x^N \,\mathrm{d}t < \infty.
\]
Finally, by $\nabla_{x^N}\Lambda^{N,\theta} = (b^N - c^N)\Lambda^{N,\theta}$ and Lemma \ref{diagtoflat}, 
\begin{equation*}
\begin{split}
\sum_{i=1}^N \,\int_{0}^T \int_{\mathcal{K}} |f|\,|\partial_{x_i^N}p_t \,\partial_{x_i^N}\Lambda^{N,\theta}| \,\mathrm{d}x \,\mathrm{d}t &\leq C\sum_{i=1}^N\,\int_{t_1}^{t_2} \int_{\mathcal{K}} |\partial_{x_i^N}p_t|\,|b_i^N|\,\Lambda^{N,\theta} + |\partial_{x_i^N}p_t|\,|c_i^N|\,\Lambda^{N,\theta} \,\mathrm{d}x \,\mathrm{d}t \\
&\leq C \sum_{i=1}^N \,\int_{t_1}^{t_2} \int_{\mathcal{W}_N} |\partial_{x_i^N}p_t(x^N)|\,\bigg(\frac{1}{x_i^N - x_{i-1}^N} + \frac{1}{x_{i+1}^N - x_{i}^N} \bigg) \,\mathrm{d}x^N \,\mathrm{d}t.
\end{split}
\end{equation*}
In view of the Cauchy-Schwarz inequality, the latter expression is at most
\[
C \sum_{i=1}^N \bigg(\int_{t_1}^{t_2} \int_{\mathcal{W}_N} p_t(x^N) \bigg(\frac{1}{(x_i^N - x_{i-1}^N)^2} + \frac{1}{(x_{i+1}^N - x_{i}^N)^2}\bigg)\, \mathrm{d}x^N \,\mathrm{d}t\bigg)^{1/2}\,
\bigg(\int_{t_1}^{t_2} \int_{\mathcal{W}_N} \frac{|\nabla p_t|^2}{p_t} \,\mathrm{d}x^N \,\mathrm{d}t\bigg)^{1/2},
\]
which is finite by Lemmas \ref{DBMint}(c) and \ref{reg}. This finishes the proof of equation \eqref{IBP}. Lastly, we combine \eqref{ftimederiv}, \eqref{ftopderiv} and \eqref{IBP} to arrive at
\begin{equation*}
\begin{split}
&\int_0^T \int_{\mathcal{K}} (\partial_t f + \mathcal{L}_{1,...,N}f)\, p_t \,\Lambda^{N,\theta}\,\mathrm{d}x \,\mathrm{d}t
= \int_0^T \int_{\mathcal{K}}- c^N \cdot \nabla_{x^N} (f \Lambda^{N,\theta})\, p_t + \frac{1}{2}\nabla_{x^N}(f\Lambda^{N,\theta})\cdot\nabla_{x^N} p_t \,\mathrm{d}x \,\mathrm{d}t \\
&\qquad\qquad\qquad\quad\;\;
+ \int_0^T \int_{\mathcal{K}} (c^N \cdot \nabla_{x^N} f)\, p_t \,\Lambda^{N,\theta}
- \frac{1}{2}\nabla_{x^N}f \cdot \nabla_{x^N} (p_t \Lambda^{N,\theta}) 
+ (\nabla_{x^N} f \cdot \nabla_{x^N}\Lambda^{N,\theta})\, p_t \,\mathrm{d}x \,\mathrm{d}t \\
&\qquad\qquad\qquad\quad\;\;
+\int_0^T \int_{\mathcal{K}} fp_t \, c^N \cdot \nabla_{x^N} \Lambda^{N,\theta} 
- \frac{1}{2}\nabla_{x^N}(fp_t)\cdot\nabla_{x^N}\Lambda^{N,\theta} \,\mathrm{d}x \,\mathrm{d}t=0.
\end{split}
\end{equation*}

\smallskip

\noindent\textbf{Step 5.} Next, we remove the assumption that for some $\varepsilon_0>0$, $f(t,x)=0$ if $\mathrm{dist}(x^N,\partial\mathcal{W}_N)<\varepsilon_0$. Given an $\varepsilon > 0$, let $\zeta_\varepsilon\!:[0,\infty)\rightarrow [0,1]$ be a $C^2([0,\infty))$-function such that $\zeta_\varepsilon(y)=0$ for $y\le\varepsilon$, $\zeta_\varepsilon(y)=1$ for $y\ge 2\varepsilon$, $|\zeta'_\varepsilon|\le C\varepsilon^{-1}\,\mathbf{1}_{\{y \ge \varepsilon\}}$, and $|\zeta''_\varepsilon|\le C\varepsilon^{-2}\,\mathbf{1}_{\{y \ge \varepsilon\}}$. Take $\phi_\varepsilon(x)= \prod_{i=1}^{N-1} \zeta_{\varepsilon}(x_{i+1}^N - x_i^N)$. 

\begin{lemma}\label{approx1}
For every $f \in \mathcal{A}_c$, there exists a constant $C_{f,\theta}<\infty$ (independent of $\mu_0$) so that 
\begin{equation}
\bigg|\int_0^T \int_{\mathcal{K}} \partial_t (f-f\phi_{\varepsilon}) + \mathcal{L}_{1,...,N}(f-f\phi_{\varepsilon})\,\mathrm{d}\nu_t \,\mathrm{d}t\bigg|
\le C_{f,\theta}\,\varepsilon^{\theta-1},\quad\varepsilon\in(0,1).
\end{equation}
\end{lemma}

\noindent\textbf{Proof of Lemma \ref{approx1}.} Let $X^N$ be a $\theta$-DBM with the initial distribution $\mu_0$, and let $Z^N$ be a $\theta$-DBM started from $0\in\mathcal{W}_N$. Write $\mathrm{dist}$ for the notion of distance on $\mathbb{R}^N$ induced by the infinity norm. Then, $\mathrm{dist}(x^N,\partial \mathcal{W}_N) \le \min_{1\le i\le N-1} |x_{i+1}^N-x_i^N| \le 2\,\mathrm{dist}(x^N,\partial \mathcal{W}_N) $. By Lemmas~\ref{comparisons}(b) and~\ref{DBMint}(b),
\begin{equation}\label{probestimate}
\begin{split}
\mu_t\big(\{x^N\in\mathcal{W}_N:\,\mathrm{dist}(x^N,\partial\mathcal{W}_N) \le 2\varepsilon\}\big)
&\le \sum_{i=1}^{N-1} \mathbb{P}(|X_{i+1}^N(t) - X_i^N(t)|\leq 4\varepsilon)\\
&\le 4^{2\theta}\varepsilon^{2\theta}\sum_{i=1}^{N-1}\mathbb{E}\big[(Z_{i+1}^N(t) - Z_i^N(t))^{-2\theta}\big] \le C_\theta\,\varepsilon^{2\theta}t^{-\theta}.
\end{split}
\end{equation}
Let $[t_1,t_2]\subset(0,T)$ satisfy $f(t,x) = 0$, $t \notin [t_1,t_2]$ and note that $|\partial_t(f-f\phi_{\varepsilon})|\leq C_f\,\mathbf{1}_{\{\mathrm{dist}(x^N, \partial \mathcal{W}_N) \leq 2\varepsilon\}}$. Therefore, it holds 
\[
\bigg|\int_0^T \int_{\mathcal{K}} \partial_t (f-f\phi_{\varepsilon})\,\mathrm{d}\nu_t\,\mathrm{d}t\bigg|
\le C_f \int_{t_1}^{t_2} \mu_t\big(\{x^N\in\mathcal{W}_N:\,\mathrm{dist}(x^N,\partial\mathcal{W}_N) \le 2\varepsilon\}\big)\,\mathrm{d}t 
\leq C_{f,\theta}\,\varepsilon^{2\theta}. 
\]
Arguing similarly, we find that 
\[
\bigg|\int_0^T \int_{\mathcal{K}} \sum_{1\leq i \leq k \leq N-1} \partial_{x_i^kx_i^k}(f-f\phi_{\varepsilon})\,\mathrm{d}\nu_t\,\mathrm{d}t\bigg| \le 
C_{f,\theta}\,\varepsilon^{2\theta-2}.
\]

\smallskip

In addition, Lemma \ref{diagtoflat} and the Cauchy-Schwarz inequality yield
\begin{equation*}
\begin{split}
&\,\bigg|\int_0^T \int_{\mathcal{K}} \sum_{1\leq i \leq k \leq N} b_i^k \, \partial_{x_i^k}(f-f\phi_{\varepsilon})\,\mathrm{d}\nu_t \,\mathrm{d}t \bigg| \\
& \le C_f \,\varepsilon^{-1}\sum_{1 \leq i \leq k \leq N} \,\int_{t_1}^{t_2} \int_{\mathcal{K}} 
|b_i^k|\,\mathbf{1}_{\{\mathrm{dist}(x^N,\partial \mathcal{W}_N) \le 2\varepsilon\}}
\,\mathrm{d}\nu_t \,\mathrm{d}t \\
 &\le C_{f,\theta}\,\varepsilon^{-1}\sum_{1 \leq i \leq k \leq N}\,
 \int_{t_1}^{t_2} \int_{\mathcal{K}^k} 
 \bigg(\frac{1}{x_i^k - x_{i-1}^{k}}+\frac{1}{x_{i+1}^k - x_i^k}\bigg) 
 \,\mathbf{1}_{\{\mathrm{dist}(x^N,\partial \mathcal{W}_N) \le 2\varepsilon\}}\,
 \prod_{m=k+1}^N \lambda^{m,\theta}(x^m,x^{m-1}) \\
&\qquad\qquad\qquad\qquad\qquad\qquad\qquad\qquad\qquad\qquad\qquad\qquad\qquad\qquad\qquad\quad\;
 \mu_t(\mathrm{d}x^N)\,\mathrm{d}x^k \,\ldots\,\mathrm{d}x^{N-1} \\
& \leq C_{f.\theta}\,\varepsilon^{-1} \sum_{1\leq i \leq k \leq N}\, \int_{t_1}^{t_2} \bigg(\int_{\mathcal{K}^k} \bigg(\frac{1}{(x_i^k - x_{i-1}^{k})^2}+\frac{1}{(x_{i+1}^k - x_i^k)^2} \bigg)\prod_{m=k+1}^N \lambda^{m,\theta}(x^m,x^{m-1}) \\ 
&\qquad\qquad\qquad\qquad\qquad\qquad\quad\;
\mu_t(\mathrm{d}x^N)
\,\mathrm{d}x^k \,\ldots\,\mathrm{d}x^{N-1}\bigg)^{1/2} 
\mu_t\big(\{x^N\in\mathcal{W}_N \!:\mathrm{dist}(x^N,\mathcal{W}_N)\le 2\varepsilon\}\big)^{1/2} \,\mathrm{d}t.
\end{split}
\end{equation*}
Letting $Z^k$ be a $\theta$-DBM started from $0\in\mathcal{W}_k$, and using Lemma \ref{intertwinings}, Lemma \ref{comparisons}(b), and \eqref{probestimate}, the latter expression can be upper bounded by
 \[
C_{f,\theta}\,\varepsilon^{\theta-1} \sum_{1\leq i \leq k \leq N} \,\int_{t_1}^{t_2} \mathbb{E}\big[(Z_i^{k}(t)-Z_{i-1}^k(t))^{-2} + (Z_{i+1}^k(t) - Z_i^k(t))^{-2}\big]^{1/2}\, t^{-\theta/2}\,\mathrm{d}t \le C_{f,\theta}\,\varepsilon^{\theta-1},
\]
where the integral is finite by Lemma \ref{DBMint}(b). It remains to notice that 
\begin{equation*}
\begin{split}
&\,\bigg|\int_0^T \int_{\mathcal{K}}\frac{1}{2}\Delta_{x^N} (f-f\phi_\varepsilon)\,\mathrm{d}\nu_t\,\mathrm{d}t \bigg| \\
& \leq C_f \int_{t_1}^{t_2}\int_{\mathcal{W}_N} 1-\phi_\varepsilon \,\mathrm{d}\mu_t \,\mathrm{d}t 
+ C_f \sum_{i=1}^N \,\int_{t_1}^{t_2}\int_{\mathcal{W}_N} |\partial_{x_i^N}\phi_\varepsilon| \,\mathrm{d}\mu_t \,\mathrm{d}t 
+ C_f \sum_{i=1}^N \, \int_{t_1}^{t_2}\int_{\mathcal{W}_N}|\partial_{x_i^Nx_i^N}\phi_\varepsilon|\,\mathrm{d}\mu_t \,\mathrm{d}t. 
\end{split}
\end{equation*}
As before, the three summands can be bounded by $C_{f,\theta}\varepsilon^{2\theta}$, $C_{f,\theta}\varepsilon^{2\theta-1}$, and $C_{f,\theta}\varepsilon^{2\theta-2}$, respectively.
Combining the four obtained estimates we arrive at the lemma. \qed

\medskip

For any $f \in \mathcal{A}_c$, we can now rely on Lemma \ref{approx1} when taking the limit $\varepsilon\downarrow0$ in 
\[
\int_0^T\int_{\mathcal{K}}(\partial_t + \mathcal{L}_{1,...,N})(\phi_\varepsilon f) 
\,\mathrm{d}\nu_t\,\mathrm{d}t=0
\]
and conclude that
\[
\int_0^T\int_{\mathcal{K}}(\partial_t + \mathcal{L}_{1,...,N})f 
\,\mathrm{d}\nu_t\,\mathrm{d}t=0.
\]

\smallskip

\noindent\textbf{Step 6.} In this final step, we remove the extra assumption on $\mu_0$ imposed at the beginning of Step~2. 

\begin{lemma}\label{edgeprob}
Fix an $\varepsilon\in(0,1)$. Then, for any $1 \le p < \theta$, there exists a constant $C_{\varepsilon,p}< \infty$ such that for all $x^N \in\mathcal{W}_N$ with $\mathrm{dist}(x^N,\partial \mathcal{W}_N)\ge\varepsilon$ and all $1\leq i \leq k \leq N-1$,
\[
\int_{\mathcal{K}(x^N)} \! (x_i^k - x_i^{k+1})^{-p}
\,\Lambda^{N,\theta}(x^N,\mathrm{d}x^1,\ldots,\mathrm{d}x^{N-1}) \,\vee\,
\int_{\mathcal{K}(x^N)} \! (x_{i+1}^{k+1} - x_i^{k})^{-p}
\,\Lambda^{N,\theta}(x^N,\mathrm{d}x^1,\ldots,\mathrm{d}x^{N-1})
\le C_{\varepsilon,p}.
\]
Moreover, with $\mathcal{K}_{\delta}(x^N) := \{(x_i^k)_{1\le i \le k \le N-1}:\, 
x_i^{k+1}+2\delta\le x_i^k +\delta\le x_{i+1}^{k+1},\,1\le i\le k\le N-1\}$ and $L_{x^N}(\mathrm{d}x^1,\ldots,\mathrm{d}x^{N-1}):=\Lambda^{N,\theta}(x^N,\mathrm{d}x^1,\ldots,\mathrm{d}x^{N-1})$, there exists a constant $C_\varepsilon<\infty$ such that for all $x^N \in \mathcal{W}_N$ with $\mathrm{dist}(x^N,\partial \mathcal{W}_N)\ge \varepsilon$,
\[
L_{x^N}(\mathcal{K}(x^N)\setminus\mathcal{K}_{\delta}(x^N))\le C_\varepsilon\delta.
\]
\end{lemma}

\noindent\textbf{Proof of Lemma \ref{edgeprob}.} We just prove the first claim, as the second claim then follows via Markov's inequality and a union bound. For the first claim, we proceed by backward induction over $k$. For $k=N-1$, Lemma \ref{diagtoflat} implies that
\[
\int_{\mathcal{K}^{N-1}(x^N)}(x_i^{N-1} - x_i^{N})^{-p}\,\lambda^{N,\theta}(x^N,x^{N-1})\,\mathrm{d}x^{N-1}
\le C_p(x_{i+1}^N - x_i^N)^{-p}\leq C_p \varepsilon^{-p} =: C_{\varepsilon,p}.
\]
For $k<N-1$, Lemma \ref{diagtoflat} and $(x^1,\ldots,x^{N-1}) \in \mathcal{K}(x^N)$ give
\begin{equation*}
\begin{split}
&\,\int_{\mathcal{K}(x^N)}(x_i^{k} - x_i^{k+1})^{-p}\,\Lambda^{N,\theta}(x^N,\mathrm{d}x^1,\ldots,\mathrm{d}x^{N-1}) \\
&\le C_p \,\int_{\mathcal{K}^{k+1}(x^N)}(x_{i+1}^{k+1} - x_i^{k+1})^{-p}\prod_{m=k+2}^N\lambda^{m,\theta}(x^m,x^{m-1}) \,\mathrm{d}x^{k+1}\,\ldots\,\mathrm{d}x^N \\
&\le C_p \,\int_{\mathcal{K}^{k+1}(x^N)}(x_{i+1}^{k+2} - x_{i}^{k+1})^{-p}
\prod_{m=k+2}^N\lambda^{m,\theta}(x^m,x^{m-1}) \,\mathrm{d}x^{k+1}\,\ldots\,\mathrm{d}x^N
\le C_{\varepsilon,p},
\end{split}
\end{equation*}
where, in the last line, we have used the induction hypothesis. The first claim readily follows. \qed

\medskip

To remove the extra assumption on $\mu_0$, we aim to approximate $\mu_0$ by $\mu^n_0$, $n\in\N$ satisfying the extra assumption. Writing $(\mu^n_t)_{t\in[0,T]}$ and $(\nu^n_t)_{t\in[0,T]}$ for the associated flows of probability measures, we then need to check that  
\begin{equation}\label{meas approx}
\int_0^T \!\int_{\mathcal{K}} (\partial_t+\mathcal{L}_{1,...,N})f \,\mathrm{d}(\nu^n_t-\nu_t)\,\mathrm{d}t
= \int_0^T \!\int_{\mathcal{W}_N} \!\int_{\mathcal{K}(x^N)} (\partial_t+\mathcal{L}_{1,...,N})f \,
\Lambda^{N,\theta}(x^N,\mathrm{d}x^1,\ldots,\mathrm{d}x^{N-1})\,
\mathrm{d}(\mu^n_t-\mu_t)\,\mathrm{d}t 
\end{equation}
goes to $0$ as $n\to\infty$. Hereby, due to Lemma \ref{approx1}, we may assume that for some $\varepsilon>0$, $f(t,x)=0$ if $\mathrm{dist}(x^N, \partial \mathcal{W}_N) < \varepsilon$. Our first claim is that 
\[
h(t,x^N) := \int_{\mathcal{K}(x^N)} (\partial_t+\mathcal{L}_{1,...,N})f \,
\Lambda^{N,\theta}(x^N,\mathrm{d}x^1,\ldots,\mathrm{d}x^{N-1})
\]  
is continuous in $x^N \!\in\!\mathcal{W}_N$. If $\mathrm{dist}(x^N,\partial \mathcal{W}_N)<\varepsilon$, we have $h\equiv0$ in a neighborhood of $x^N$.
 
\medskip
 
Now, consider $x^N \in \mathcal{W}_N$ such that $\mathrm{dist}(x^N,\partial \mathcal{W}_N)\ge\varepsilon$ and let $x^{N,n} \rightarrow x^{N}$ with $\mathrm{dist}(x^{N,n},x^N) < \varepsilon/2$ for all $n$. For $\delta\in(0,\varepsilon/4)$, we bound the derivatives of $f$ by a constant and use Lemma \ref{edgeprob} to obtain
\begin{equation*}
\begin{split}
& \,\bigg| \int_{\mathcal{K}(x^{N,n})\setminus \mathcal{K}_{\delta}(x^{N,n})} 
\Big(\partial_t f+ \frac{1}{2}\Delta_x f\Big)(t,x^{N,n},x^1,\ldots,x^{N-1}) \,
\Lambda^{N,\theta}(x^{N,n},\mathrm{d}x^1,\ldots,\mathrm{d}x^{N-1}) \bigg| \\
& \le C_f L_{x^{N,n}}(\mathcal{K}(x^{N,n})\setminus \mathcal{K}_{\delta}(x^{N,n})) 
\le C_{f,\frac{\varepsilon}{2}}\,\delta.
\end{split}
\end{equation*}
For the same reason, 
\[
\bigg| \int_{\mathcal{K}(x^N)\setminus \mathcal{K}_{\delta}(x^N)} 
\Big(\partial_t f+ \frac{1}{2}\Delta_x f\Big)(t,x) \,
\Lambda^{N,\theta}(x^N,\mathrm{d}x^1,\ldots,\mathrm{d}x^{N-1}) \bigg|
\le C_{f,\varepsilon}\delta.
\]
Moreover, by the Bounded Convergence Theorem, 
\begin{equation*}
\begin{split}
& \lim_{n\to\infty}\,\int_{\mathbb{R}^{\frac{(N-1)N}{2}}} \mathbf{1}_{\mathcal{K}(x^{N,n})\cap \mathcal{K}_{\delta}(x^{N,n})}\, 
\Big(\partial_t f+\frac{1}{2}\Delta_x f\Big)(t,x^{N,n},x^1,\ldots,x^{N-1})\, \Lambda^{N,\theta}(x^{N,n},x^1,\ldots,x^{N-1}) \\
&\qquad\qquad\qquad\;
-\mathbf{1}_{\mathcal{K}(x^{N})\cap \mathcal{K}_{\delta}(x^{N})} 
\Big(\partial_t f + \frac{1}{2}\Delta_x f\Big)(t,x)\,\Lambda^{N,\theta}(x^N,x^1,\ldots,x^{N-1}) 
\,\mathrm{d}x^1\,\ldots\,\mathrm{d}x^{N-1} = 0.
\end{split}
\end{equation*}
Therefore, $\int_{\mathcal{K}(x^{N})} \big(\partial_t f + \frac{1}{2}\Delta_x f\big)(t,x)\, \Lambda^{N,\theta}(x^N,\mathrm{d}x^1,\ldots,\mathrm{d}x^{N-1})$ is continuous at $x^N$.

\medskip

Fix $p\in(1,\theta)$, and let $q$ satisfy $p^{-1}+q^{-1}=1$. Then, bounding $\partial_{x_i^k}f$ by a constant and employing Hölder's inequality we have, with the convention $x^{k,n}=x^k$ for $k=1,\,2,\,\ldots,\, N-1$, 
\begin{equation*}
\begin{split}
&\,\bigg|\int_{\mathcal{K}(x^{N,n})\setminus \mathcal{K}_{\delta}(x^{N,n})} b_i^k(x^{k,n},x^{k-1})\, \partial_{x_i^{k}}f(t,x^{N,n},x^1,\ldots,x^{N-1})\, \Lambda^{N,\theta}(x^{N,n},\mathrm{d}x^1,\ldots,\mathrm{d}x^{N-1})\bigg| \\
&\le C_f\bigg(\int_{\mathbb{R}^{\frac{(N-1)N}{2}}} \mathbf{1}_{\mathcal{K}(x^{N,n})\setminus \mathcal{K}_{\delta}(x^{N,n})}\, \Lambda^{N,\theta}(x^{N,n},\mathrm{d}x^1,\ldots,\mathrm{d}x^{N-1})\bigg)^{1/q} \\
&\quad\cdot\bigg(\int_{\mathbb{R}^{\frac{(N-1)N}{2}}} \big((x_{i}^{k-1} - x_{i}^{k,n})^{-p} + (x_{i}^{k,n}-x_{i-1}^{k-1})^{-p} \big)\,\Lambda^{N,\theta}(x^{N,n},\mathrm{d}x^1,\ldots,\mathrm{d}x^{N-1})\bigg)^{1/p}
\le C_{f,\frac{\varepsilon}{2},p}\,\delta^{1/q},
\end{split}
\end{equation*}
where the last inequality is due to Lemmas \ref{edgeprob} and \ref{diagtoflat}. The same argument also shows that
\begin{equation}\label{h bounded}
\bigg|\int_{\mathcal{K}(x^{N})\setminus \mathcal{K}_{\delta}(x^{N})}b_i^k(x^{k},x^{k-1})\, \partial_{x_i^{k}}f(t,x)\,\Lambda^{N,\theta}(x^{N},\mathrm{d}x^1,\ldots,\mathrm{d}x^{N-1})\bigg|
\leq C_{f,\varepsilon,p}\,\delta^{1/q}.
\end{equation}
Consequently, 
\begin{equation*}
\begin{split}
&\bigg| \int_{\mathcal{K}(x^{N,n})}b_i^k(x^{k,n},x^{k-1})\, \partial_{x_i^{k}}f(t,x^{N,n},x^1,\ldots,x^{N-1})\,\Lambda^{N,\theta}(x^{N,n},\mathrm{d}x^1,\ldots,\mathrm{d}x^{N-1}) \\
&\, -\int_{\mathcal{K}(x^{N})} b_i^k(x^k,x^{k-1})\,\partial_{x_i^{k}}f(t,x)\, \Lambda^{N,\theta}(x^N,\mathrm{d}x^1,\ldots,\mathrm{d}x^{N-1})\bigg| \\
&\leq  C_{f,\frac{\varepsilon}{2},p}\,\delta^{1/q} + C_{f,\varepsilon,p}\,\delta^{1/q} \\
&\quad +\int_{\mathbb{R}^{\frac{(N-1)N}{2}}} \Big|\mathbf{1}_{\mathcal{K}(x^{N,n})\cap\mathcal{K}_{\delta}(x^{N,n})}\, b_i^k(x^{k,n},x^{k-1}) \,\partial_{x_i^{k}}f(t,x^{N,n},x^1,\ldots,x^{N-1})\, \Lambda^{N,\theta}(x^{N,n},x^1,\ldots,x^{N-1}) \\ 
&\qquad\qquad\qquad\;
- \mathbf{1}_{\mathcal{K}(x^{N})\cap\mathcal{K}_{\delta}(x^{N})}\, b_i^k(x^k,x^{k-1})\, \partial_{x_i^{k}}f(t,x)\,\Lambda^{N,\theta}(x^{N},x^1,\ldots,x^{N-1})\Big| \,\mathrm{d}x^1\,\ldots\,\mathrm{d}x^{N-1}. 
\end{split}
\end{equation*}
The last integral goes to $0$ as $n\to\infty$ by the Bounded Convergence Theorem. Taking $n\to\infty$ and~then $\delta\downarrow0$, we see that $h(t,x^N)$ is continuous in $x^N\in\mathcal{W}_N$. Since $h$ is bounded (recall $f\in\mathcal{A}_c$ and \eqref{h bounded}), 
\[
\int_{\mathcal{K}} (\partial_t+\mathcal{L}_{1,...,N})f \,\mathrm{d}(\nu^n_t-\nu_t)
= \int_{\mathcal{W}_N} \int_{\mathcal{K}(x^N)} (\partial_t+\mathcal{L}_{1,...,N})f \,
\Lambda^{N,\theta}(x^N,\mathrm{d}x^1,\ldots,\mathrm{d}x^{N-1})\,
\mathrm{d}(\mu^n_t-\mu_t)
\]
goes to $0$ in the limit $n\to\infty$, as long as the weak convergence $\mu^n_t\to\mu_t$ holds.

\medskip

Next, we observe that thanks to estimate \eqref{int1} and Lemma \ref{comparisons}(b), 
\begin{equation*} 
\begin{split}
\int_{\mathcal{W}_N} \int_{\mathcal{K}(x^N)} (\partial_t+\mathcal{L}_{1,...,N})f \,
\Lambda^{N,\theta}(x^N,\mathrm{d}x^1,\ldots,\mathrm{d}x^{N-1})\,
\mathrm{d}\mu_t 
&\le C + C\sum_{k=2}^N \sum_{i=1}^{k-1} 
\mathbb{E}\bigg[\frac{1}{X_{i+1}^k(t) - X_i^{k}(t)}\bigg] \\
&\le C + C\sum_{k=2}^N \sum_{i=1}^{k-1} 
\mathbb{E}\bigg[\frac{1}{Z_{i+1}^k(t) - Z_i^{k}(t)}\bigg],
\end{split}
\end{equation*}
where each $X^k$ is a $\theta$-DBM with initial density
$\int_{\mathcal{K}^{k+1}}\! \prod_{m=k+1}^N \! \lambda^{m,\theta}(x^m,x^{m-1})\,\mathrm{d}x^{k+1}\ldots\mathrm{d}x^{N-1}\mu_0(\mathrm{d}x^N)$ and each $Z^k$ is a $\theta$-DBM started from $0\in\mathcal{W}_k$. The latter upper bound is integrable over $[0,T]$ thanks to Lemma \ref{DBMint}(b), so that by the Dominated Covergence Theorem, the quantities in \eqref{meas approx} go to $0$ in the limit $n\to\infty$, as long as the weak convergence $\mu^n_t\to\mu_t$ holds for all $t\in[0,T]$. 

\medskip

In view of Lemma \ref{wass}, the Fokker-Planck equation \eqref{FP} is fulfilled for every probability measure $\mu_0$ on $\mathcal{W}_N$ such that $\mu_0$ can be approximated in the $1$-Wasserstein distance by probability measures $\mu_0^n$ on $\mathcal{W}_N$ for which \eqref{FP} has been established. If $\mu_0$ is a Dirac mass at $x\in\,\stackrel{\circ}{\mathcal{W}}_N$, a desired approximation is given by the convolutions with a scaled bump function, which satisfy the extra assumption imposed at the beginning of Step 2. If $\mu_0$ is a Dirac mass at $x \in \partial{\mathcal{W}_N}$, a desired approximation is given by the Dirac masses at $\stackrel{\circ}{\mathcal{W}}_N \,\ni x^n$, $x^n_i:=x_i+\frac{i}{n}$. Thus, \eqref{FP} holds for any $\mu_0$ given by a finite convex combination of Dirac masses in $\mathcal{W}_N$. Since every probability measure $\mu_0$ on $\mathcal{W}_N$ can be approximated in the $1$-Wasserstein distance by such convex combinations, the proof of \eqref{FP} is finished. \qed  

\section{Proof of Theorem \ref{MainTheorem}} \label{sec:conseq}

Section \ref{super context} allows us to apply Proposition \ref{tsuper}, hence concluding the existence of a solution $\eta$ to the martingale problem associated with $ \mathcal{L}_{1,...,N}$ such that $\eta\circ e_t^{-1}=\nu_t$, $t\in[0,T]$. A standard localization argument shows that $\eta$ is a solution to the local martingale problem in the sense of \cite[Chapter 5, Definition 4.5]{karatzas1991brownian}. By \cite[Chapter 5, Proposition 4.6]{karatzas1991brownian},  $\eta$ induces a weak solution to the SDE \eqref{MLSDE}. The remaining statements of Theorem \ref{MainTheorem} are restated for convenience in the next proposition. 

\begin{proposition}\label{proofofthm1}
Let $\eta$ be a solution to the martingale problem associated with $ \mathcal{L}_{1,...,N}$ such that $\eta \circ e_t^{-1} = \nu_t$, $t\ge0$, and let $X = (X_i^k)_{1 \leq i \leq k \leq N}$ be distributed according to $\eta$. Then:
\begin{enumerate}[(a)]
\item $X$ is a Markov process.
\item If $\eta'$ also satisfies the assumptions of the proposition, then $\eta' = \eta$.
\item Each $X^k$ is a $\theta$-DBM. 
\end{enumerate}
\end{proposition}

\noindent\textbf{Proof. (a).} The proof closely follows \cite[Chapter 5, Lemma 4.19 and Theorem 4.20]{karatzas1991brownian}, so we omit the details that are the same and emphasize the differences. Fix a $t_0>0$, and for $\omega \in C\big([0,t_0],\mathbb{R}^{\frac{N(N+1)}{2}}\big)$, let $Q_\omega$ be the regular conditional distribution $\eta$ given $\mathcal{G}_{t_0} := \sigma(s\mapsto\omega(s),\,s\in[0,t_0])$. Write $\theta_{t_0}$ for the time shift operator, and define $P_{\omega}= Q_{\omega}\circ \theta_{t_0}^{-1}$.
For given $0\!\le\! s \!<\! t \!<\!\infty$ and $f \in C_b^2\big(\mathbb{R}^{\frac{N(N+1)}{2}}\big)$, set
\[
Z(y) = f(y(t))-f(y(s)) - \int_{s}^t(\mathcal{L}_{1,...,N}f)(y(r)) \,\mathrm{d}r.
\]
Note that $\mathbb{E}[\int |Z(y)| \, P_\omega(\mathrm{d}y)] = \mathbb{E}[|Z \circ \theta_{t_0}|] < \infty$, so  $\int |Z(y)| \, P_\omega(\mathrm{d}y)< \infty$ almost surely. By repeating \cite[Chapter 5, proof of Lemma 4.19]{karatzas1991brownian}, we find an $\eta$-null event such that for all $\omega$ in its complement, $P_\omega$ solves the martingale problem associated with $\mathcal{L}_{1,...,N}$. Now, off the $\eta$-null event that $\omega(t_0) \in \partial \mathcal{K}$, the martingale problem started from $\omega(t_0)$ has a unique solution $\eta^{\omega(t_0)}$ by \cite[Theorem 1.2]{gorin2014multilevel}. Therefore, for any Borel $A \subset C\big([0,\infty),\mathbb{R}^{\frac{N(N+1)}{2}}\big)$ and any $\omega$ off an $\eta$-null event,
\[
\eta(\theta_{t_0}^{-1}A | \mathcal{G}_{t_0}) = Q_{\omega}(\theta_{t_0}^{-1}A) = P_{\omega}(A) 
= \eta^{\omega(t_0)}(A).
\]
This amounts to the Markov property of $X$. 

\medskip

\noindent\textbf{(b).} For any $t_0>0$, the restrictions of $\eta'$, $\eta$ to $C\big([t_0,\infty),\R^{\frac{N(N+1)}{2}}\big)$ solve the martingale problem associated with $\mathcal{L}_{1,...,N}$ starting from $\nu_{t_0}$. Since $\nu_{t_0}(\partial\mathcal{K})=0$, these restrictions are the same by \cite[Theorem 1.2]{gorin2014multilevel}. This implies $\eta'=\eta$ because both are measures on $C\big([0,\infty),\R^{\frac{N(N+1)}{2}}\big)$.   

\medskip

\noindent\textbf{(c).} We first show that $X^N$ is a $\theta$-DBM. Our goal is to use \cite[Theorem 2]{3e0b7045-d14d-39db-8a2a-09b90e489549}. To this end, we define a transition kernel from $\mathcal{W}$ to $\mathcal{K}$ by 
\[
\overline{\Lambda}(x^N,\mathrm{d}y^1,\ldots,\mathrm{d}y^N) = \delta_{x^N}(\mathrm{d}y^N)\,\Lambda^{N,\theta}(x^N,\mathrm{d}y^1,\ldots,\mathrm{d}y^{N-1}).
\] 
Let $\Phi\!: C_b(\mathcal{W})\to C_b(\mathcal{K})$ be given by $(\Phi f)(x) = f(x^N)$. Clearly, $\overline{\Lambda}\Phi = I$, the identity on $C_b(\mathcal{W})$. Write $(P_t)_{t\ge0}$ for the transition kernels of $X$. Set $Q_t = \overline{\Lambda}P_t \Phi$. Then, for any $f\in C_b(\mathcal{W})$ and $x^N \in \mathcal{W}$,
\[
(Q_tf)(x^N) = \int_{\mathcal{K}}\int_{\mathcal{K}} f(z^N) \, 
P_t(y,\mathrm{d}z)\,\overline{\Lambda}(x^N,\mathrm{d}y)
= \mathbb{E}[f(Y^N(t))] 
= \int_{\mathcal{W}} f(y^N)\, p_t(x^N,\mathrm{d}y^N),
\]
where $Y$ is constructed via superposition from $\mu_0=\delta_{x^N}$, and $p_t$ is the transition kernel of the $\theta$-DBM. Therefore, $Q_t$ is the transition kernel of the $\theta$-DBM. 

\medskip

To apply \cite[Theorem 2]{3e0b7045-d14d-39db-8a2a-09b90e489549}, it remains to show that $\overline{\Lambda}P_t = Q_t\overline{\Lambda}$. This follows from the superposition construction, since for any $f\in C_b(\mathcal{K})$ and $x^N\in\mathcal{W}$,
\begin{equation*}
\begin{split}
(\overline{\Lambda}P_tf)(x^N) 
& =\int_{\mathcal{K}(x^N)} \int_{\mathcal{K}} f(y) \, P_t(x,\mathrm{d}y)
\, \Lambda^{N,\theta}(x^N,\mathrm{d}x^1,\ldots,\mathrm{d}x^{N-1})\\
& = \mathbb{E}[f(Y(t))] = \int_{\mathcal{W}} \int_{\mathcal{K}(y^N)} f(y)
\,\Lambda^{N,\theta}(y^N,\mathrm{d}y^1,\ldots,\mathrm{d}y^{N-1})\, p_t(x^N,\mathrm{d}y^N) 
= (Q_t \overline{\Lambda}f)(x^N).
\end{split}
\end{equation*}
Hence, $X^N$ is Markov with the transition kernels $(Q_t)_{t\ge0}$, and thus a $\theta$-DBM, by \cite[Theorem 2]{3e0b7045-d14d-39db-8a2a-09b90e489549}. Using the same argument for $(X^m_i)_{1\le i\le m\le k}$, which satisfies the assumptions of the proposition with $\mathrm{d}\nu_t = \big(\int_{\mathcal{K}^{k+1}} \prod_{m=k+1}^N \lambda^{m,\theta}(x^m,x^{m-1})\, \mathrm{d}x^{k+1}\,\ldots\,\mathrm{d}x^{N-1}\,\mu_t(\mathrm{d}x^N)\big)\,\prod_{m=2}^k \lambda^{m,\theta}(x^m,x^{m-1})\,\mathrm{d}x^1\,\ldots\,\mathrm{d}x^k$, $t\ge0$, we see that $X^k$ is a $\theta$-DBM for any $k=1,\,2,\,\ldots,\,N-1$. \qed

\section{Proof of Theorem \ref{convergence}} \label{sec:Thm2}

\subsection{Preliminaries on the Submartingale Problem}

Let us recall the submartingale problem and the stochastic differential equations with reflection (SDERs) from \cite{kang2014submartingale}. Throughout the rest of this subsection, $b\!: \mathbb{R}^J \rightarrow \mathbb{R}^J$ and $\sigma\!: \mathbb{R}^J \rightarrow \mathbb{R}^{J\times J}$ are locally bounded measurable functions, and we let $\mathcal{C} := C([0,\infty),\mathbb{R}^J)$. Now, we define domains for which the submartingale problem is most helpful.

\begin{definition}
We say an open $\emptyset\neq G \subset \mathbb{R}^J$ endowed with a multi-valued vector field $\upsilon$ on $ \overline{G}$ is \textit{piecewise} $C^2$ \textit{with continuous reflection} if:
\begin{enumerate}
\item $G = \bigcap_{i \in \mathcal{I}} G_i$, where $\mathcal{I}$ is a finite index set and each $G_i$ has a $C^2$ boundary. (For $x\in\partial G_i$, we let $n^i(x)$ be the unit inward normal to $\partial G_i$ at $x$. Also, we let $\mathcal{I}(x) := \{ i \in \mathcal{I} : x \in \partial G_i \}$.)
\item The reflection field $\upsilon : \overline{G} \rightarrow 2^{\mathbb{R}^J}$ is given by $\upsilon(x) = \{0\}$ for $x\in G$, and 
\[
\upsilon(x) = \bigg\{ \sum_{i \in \mathcal{I}(x)}s_i \upsilon^i(x):\, s_i \geq 0,\, i\in\mathcal{I}(x) \bigg\}
\]
for $x \in \partial G$ where each $\upsilon^i$ is a continuous unit vector field on $\partial G_i$ with $\langle n^i, \upsilon^i \rangle >0$.
\end{enumerate}
\end{definition}

\smallskip 

Define $n(x)$ as the set of inward normal vectors at $x \in \partial G$. We introduce
\begin{align}\label{goodset}
\mathcal{U} := \big\{x \in \partial G: \;\exists\,n \in n(x)\,\forall \upsilon\in\upsilon(x)\setminus\{0\}:\, \langle n , \upsilon \rangle >0 \big\}
\end{align}
and the \textit{exceptional set} $\mathcal{V} := \partial G \setminus \mathcal{U}$. 
Next, we recall the definition of the extended Skorokhod problem (ESP) which is used to define solutions of SDERs.

\begin{definition}
Fix $(G,\upsilon)$ and $ \psi \in \mathcal{C}$ with $\psi(0) \in \overline{G}$. Then, the pair $(\phi,\eta) \in \mathcal{C}^2$ is said to solve the ESP if $\phi(0) = \psi(0)$, and for all $t \geq s \geq 0$, we have:
\begin{enumerate}
\item $\phi(t) = \psi(t) + \eta(t)$.
\item $\phi(t) \in \overline{G}$.
\item $\eta(t) - \eta(s) \in \overline{\text{co}}(\cup_{r \in [s,t]}\,\upsilon(\phi(r)))$, where $\overline{\text{co}}$ stands for the closure of the convex hull.
\end{enumerate}
\end{definition}

\begin{definition}\label{def:SDER}
Given a probability measure $\nu_0$ on $\overline{G}$, a solution to an SDER associated with $(G,\upsilon), b$, $\sigma$ and $\nu_0$ is a triplet $(\Omega,\mathcal{F}, \{\mathcal{F}_t \} ), \mathbb{P}_{\nu_0}, (X,B)$, where $(\Omega,\mathcal{F}, \{\mathcal{F}_t \} )$ is a filtered probability space supporting a probability measure $\mathbb{P}_{\nu_0}$ and a pair $(X,B)$ of continuous processes, such that: 
\begin{enumerate}
\item For every Borel $A\subset\mathbb{R}^J$, $\mathbb{P}_{\nu_0}(X(0)\in A)=\nu_0(A)$.
\item $(B_t,\mathcal{F}_t)$ is a $J$-dimensional standard Brownian motion under $\mathbb{P}_{\nu_0}$.
\item $\mathbb{P}_{\nu_0}\big(\int_0^t |b(X(s))|\,\mathrm{d}s 
+ \int_0^t |\sigma(X(s))|^2 \,\mathrm{d}s < \infty\big)=1$ for all $t \geq 0$.
\item There exists a continuous $\mathcal{F}_t$-adapted process $Y$ such that $\mathbb{P}_{\nu_0}$-almost surely, $(X,Y)$ solves the ESP associated with $(G,\upsilon)$ and 
\[
Z(t) := X(0) + \int_0^t b(X(s)) \,\mathrm{d}s + \int_0^t \sigma(X(s)) \,\mathrm{d}B_s,\quad t\ge0.
\]
\item $\mathbb{P}_{\nu_0}$-almost surely, the set $\{t\ge0\!: X(t) \in \partial G \}$ has zero Lebesgue measure.
\end{enumerate}
\end{definition}

\smallskip

To define the submartingale problem, we introduce the space $C_c^2(\overline{G}) \oplus \mathbb{R}$ of functions in $C_c^2(\overline{G})$ plus a constant, as well as the test functions space
\begin{equation}
\mathcal{H} := \big\{h \in C_c^2(\overline{G}) \oplus \mathbb{R}:\, \nabla h\equiv0 \text{ on a neighborhood of } \mathcal{V}, \text{ and } \langle \nabla h(x), \upsilon \rangle \ge 0,\,\upsilon \in \upsilon(x),\, x \in \partial G\big\}.
\end{equation}
Also, assuming throughout that $a := \sigma\sigma^\top$ is uniformly elliptic, we let
\[
\mathcal{L} :=  \frac{1}{2}\sum_{i,j=1}^J a_{ij}\,\partial_{x_i x_j} + \sum_{i=1}^J b_i\,\partial_{x_i}.
\]

\begin{definition}\label{submart}
Given a probability measure $\nu_0$ on $\overline{G}$, a probability measure $\mathbb{Q}_{\nu_0}$ on $(\mathcal{C},\mathcal{B}(\mathcal{C}))$ is said to solve the \textit{submartingale problem} with initial distribution $\nu_0$ associated with $(G,\upsilon), b, \sigma$ if:
\begin{enumerate}
\item For every Borel $A \subset \mathbb{R}^J$, $\mathbb{Q}_{\nu_0}(\omega(0)\in A)
=\nu_0(A)$.
\item $\mathbb{Q}_{\nu_0}(\omega(t)\in \overline{G},\,t\ge0)=1$.
\item For every $h \in \mathcal{H}$, the following process is a $\mathbb{Q}_{\nu_0}$-submartingale:
\[ 
h(\omega(t)) - \int_0^t\mathcal{L}h(\omega(s)) \,\mathrm{d}s. 
\]
\item $\mathbb{Q}_{\nu_0}$-almost surely, $\int_0^{\infty} \mathbf{1}_{\mathcal{V}}(\omega(s)) \,\mathrm{d}s = 0$.
\end{enumerate}
\end{definition}

\smallskip

We need the implication ``submartingale problem $\Rightarrow$ SDER'', in the sense of the next proposition. It is proven as \cite[Theorem 3]{kang2014submartingale}, and while the latter is stated for deterministic initial conditions, the same proof works also for random initial conditions.

\begin{proposition}\label{correspondence}
Suppose $(G,\upsilon)$ is piecewise $C^2$ with continuous reflection, $\mathcal{V}$ is a finite union of closed connected sets, and $\mathbb{Q}_{\nu_0}$ solves the submartingale problem with initial distribution $\nu_0$ associated with $(G,\upsilon)$, $b$, $\sigma$. Let $X(\omega,t):=\omega(t)$ be the canonical process on $\mathcal{C}$, write $\{\mathcal{F}_t\}$ for the filtration generated by the projections $t\mapsto\omega(t)$, and consider the $\mathcal{F}_t$-stopping time
\[
\tau_{\mathcal{V}}:=\inf\{t\geq 0:\,\omega(t) \in \mathcal{V} \}.
\]
Then, there exists a continuous process $B$ on $(\mathcal{C},\mathcal{B}(\mathcal{C}),\{\mathcal{F}_t\})$ with the property that $(\mathcal{C},\mathcal{B}(\mathcal{C}),\{\mathcal{F}_t\})$, $\mathbb{Q}_{\nu_0}$, $(X(\cdot \wedge \tau_{\mathcal{V}}),B(\cdot \wedge \tau_{\mathcal{V}}))$ is a solution to the SDER associated with $(G,\upsilon)$, $b$, $\sigma$ and $\nu_0$.
\end{proposition}

\subsection{Proof of Theorem \ref{convergence}}

We break the proof of Theorem \ref{convergence} into three steps. First, we identify the exceptional set in the submartingale problem. Second, we prove tightness of the initial distributions and conclude the tightness of the processes. Third, we show that any limit point must be the Warren process.

\medskip

\noindent\textbf{Step 1.} We begin by noting that
\[
\stackrel{\circ}{\mathcal{K}}\; = \bigcap_{1 \leq i \leq k \leq N-1} (G_{k,i,+}\cap G_{k,i,-}),
\]
where $G_{k,i,+} := \big\{x \in \mathbb{R}^{\frac{N(N+1)}{2}}\!: x_i^k < x_{i+1}^{k+1} \big\}$ and $G_{k,i,-} := \big\{x \in \mathbb{R}^{\frac{N(N+1)}{2}} : x_i^k > x_{i}^{k+1}\big\}$. We define the reflection vectors $\upsilon^{k,i,+}(x)= (\upsilon_{m,j}^{k,i,+}(x))_{1 \leq j \leq m \leq N}$ and $\upsilon^{k,i,-}(x)= (\upsilon_{m,j}^{k,i,-}(x))_{1 \leq j \leq m \leq N}$ as
\[
\upsilon_{m,j}^{k,i,+}(x) =  \delta_{m,k+1}\,\delta_{j,i+1}\,\mathbf{1}_{\partial G_{k,i,+}}(x) 
\quad\text{and}\quad
\upsilon_{m,j}^{k,i,-}(x) =  -\delta_{m,k+1}\,\delta_{j,i}\,\mathbf{1}_{\partial G_{k,i,-}}(x), 
\]
where $\delta_{\cdot,\cdot}$ stands for the Kronecker delta. 

\begin{proposition}\label{vset}
The exceptional set $\mathcal{V}$ in the setup of the preceding paragraph is given by
\begin{equation}\label{Vchar}
\mathcal{V} = \{x \in \partial \mathcal{K}:\, \exists\,1 \leq i \leq k-1 \leq N-2 \text{ such that } x_i^k = x_{i+1}^k\}.
\end{equation}
\end{proposition}

\noindent\textbf{Proof of Proposition \ref{vset}.} Let $\mathcal{I}_+(x) := \{(k,i):\,1 \leq i \leq k \leq N-1,\, x \in \partial G_{k,i,+} \}$ and $\mathcal{I}_-(x) := \{(k,i):\, 1 \leq i \leq k \leq N-1,\, x \in \partial G_{k,i,-} \}$. Write $n^{k,i,+}(x)$ and $n^{k,i,-}(x)$ for the inward unit normal vectors to $\partial G_{k,i,+}$ at $x\in\partial G_{k,i,+}$ and to $\partial G_{k,i,-}$ at $x\in\partial G_{k,i,-}$, respectively. The set of inward normal vectors at $x \in \partial \mathcal{K}$ is then
\[
n(x) = \bigg\{\sum_{(k,i) \in \mathcal{I}_{+}(x)} s^{k,i,+} n^{k,i,+}(x) 
+\sum_{(k,i) \in \mathcal{I}_{-}(x)} s^{k,i,-} n^{k,i,-}(x):\, s^{k,i,\pm}\ge0,\,(k,i)\in  \mathcal{I}_\pm(x)\bigg\}.
\]
We prove the proposition by induction over $N$. For $N=2$, it is equivalent to $\mathcal{U} = \partial \mathcal{K}$. Fix an $x \in \partial \mathcal{K}$. If $x \in \partial G_{1,1,+}$ and $x \notin \partial G_{1,1,-}$, any $\upsilon \in \upsilon(x)\setminus\{0\}$ satisfies $\upsilon_2^2>0$ and $\upsilon_1^1=\upsilon_1^2=0$, so for $n$ given by $n_2^2 = -n_1^1 = 1$ and $n_1^2=0$, it holds $n \in n(x)$ and $\langle n, \upsilon \rangle > 0$ for all $\upsilon \in \upsilon(x)\setminus \{0 \}$. For $x \in \partial G_{1,1,-}$ and $x \notin \partial G_{1,1,+}$, the argument is similar. If $x \in \partial G_{1,1,+} \cap \partial G_{1,1,-}$, let $n_1^1 = 0$, $n_1^2 = - n_2^2 = -1$. Then, $n \in n(x)$. Any $\upsilon \in \upsilon(x)\setminus\{ 0\}$ satisfies $\upsilon_1^2\leq 0$ and $\upsilon_2^2\geq 0$, where one of the inequalities is strict. Therefore, $\langle n, \upsilon \rangle > 0$. The base case is complete.

\medskip 

For $N >2$, we assume the proposition is true for $N-1$. Consider an $x \in \partial \mathcal{K}$ with $x_1^k < \cdots < x_k^k$, $2 \leq k \leq N-1$. By the induction hypothesis, there is an $n \in n(x)$ such that $n_i^N=0$, $i=1,\,2,\,\ldots,\, N$,
\begin{align}\label{plusineq}
\langle n, \upsilon^{k,i,+} \rangle > 0 \text{ for all } (k,i) \in \mathcal{I}_+(x) \text{ with } k < N-1,
\end{align}
and
\begin{align}\label{minusineq}
\langle n, \upsilon^{k,i,-} \rangle > 0 \text{ for all } (k,i) \in \mathcal{I}_-(x) \text{ with } k < N-1.
\end{align}
If $(N-1,1) \in \mathcal{I}_-(x)$, we define $p_1^N = -p_1^{N-1} = -1$ and $p_i^k=0$ for all other $(k,i)$. Then, for an $\varepsilon_1 > 0$ small enough, \eqref{plusineq} and \eqref{minusineq} are maintained for $\widetilde{n} := n + \varepsilon_1 p$, $\widetilde{n} \in n(x)$, and $\langle \widetilde{n}, \upsilon^{N-1,1,-} \rangle > 0$. We drop the tilde and refer to $\widetilde{n}$ simply as $n$.

\medskip

If $(N-1,1) \in \mathcal{I}_+(x)$, then our assumption on $x$ implies $(N-1,2) \notin \mathcal{I}_-(x)$. Set $p_2^N = -p_1^{N-1}=1$ and $p_i^k = 0$ for all other $(k,i)$. For an $\varepsilon_2>0$ small enough, \eqref{plusineq} and \eqref{minusineq} are maintained for $\widetilde{n} := n + \varepsilon_2 p$, $\widetilde{n} \in n(x)$, and $\langle \widetilde{n}, \upsilon^{N-1,1,+} \rangle > 0$. Note that the $(N,1)$ component of $n$ stayed the same. If instead $(N-1,2) \in \mathcal{I}_-(x)$, then $(N-1,1) \notin \mathcal{I}_+(x)$. Define $p_2^N = - p_2^{N-1}=-1$ and $p_i^k = 0$ for all other $(k,i)$. Then, for an $\varepsilon_2>0$ small enough, \eqref{plusineq} and \eqref{minusineq} are maintained for $\widetilde{n} := n + \varepsilon_2 p$, $\widetilde{n} \in n(x)$, and $\langle \widetilde{n}, \upsilon^{N-1,2,-} \rangle > 0$. Note that the $(N,1)$ component of $n$ stayed the same. Proceeding in this way, we find a vector $n \in n(x)$ such that \eqref{plusineq} and \eqref{minusineq} hold up to $k=N-1$. Thus, $x \in \mathcal{U}$.

\medskip

Conversely, let $x \in \mathcal{K}$ be such that $x_i^k = x_{i+1}^k$ for some $1 \leq i \leq k-1 \leq N-2$. Then, $(k,i) \in \mathcal{I}_+(x)$ and $(k,i+1) \in \mathcal{I}_{-}(x)$. Let $n \in n(x)$ be arbitrary. If $n_{i+1}^{k+1}=0$, then $\langle n, \upsilon^{k,i,+} \rangle = 0$. If $n_{i+1}^{k+1}<0$, then $\langle n, \upsilon^{k,i,+} \rangle < 0$. If $n_{i+1}^{k+1}>0$, then $\langle n, \upsilon^{k,i+1,-} \rangle < 0$. Therefore, $x \notin \mathcal{U}$. \qed

\medskip

\noindent\textbf{Step 2.} We turn our attention to tightness. 

\begin{lemma}\label{tight}
In the setting of Theorem \ref{convergence}, it holds $\mu_0^\theta\times\Lambda^{N,\theta}\to \mu_0 \times \Lambda^{N,1}$ weakly as $\theta\downarrow1$.
\end{lemma}

\noindent\textbf{Proof of Lemma \ref{tight}.} We prove the following statement, whose iterative application yields the lemma. For any $2 \leq m \leq N$ and any probability measures $\{\nu_0^\theta\}_{\theta>1}$ on $\mathcal{K}^m$ such that $\nu_0^\theta\to\nu_0$ weakly as $\theta\downarrow1$ and that $\nu_0(\partial \mathcal{K}^m)=0$, one has $\nu_0^{\theta} \times \lambda^{m,\theta}\to\nu_0 \times \lambda^{m,1}$ weakly as $\theta\downarrow1$. Since $\lambda^{m,\theta}(x^m,\cdot)$ is supported in $[x_1^m,x_m^m]^{m-1}$, tightness of $\{\nu_0^\theta\times \lambda^{m,\theta}\}_{\theta\downarrow1}$ follows directly from that of $\{\nu_0^\theta\}_{\theta\downarrow1}$. Hence, it suffices to check that for any $f\in C_b(\mathcal{K}^{m-1})$ compactly supported in the interior of $\mathcal{K}^{m-1}$,
\begin{equation}\label{lemma5.7claim}
\begin{split}
&\,\lim_{\theta\downarrow 1}\; \int_{\mathcal{K}^m}\int_{\mathcal{R}(x^m)}
f(x^{m-1},\ldots,x^N)\,\lambda^{m,\theta}(x^m,x^{m-1})\,\mathrm{d}x^{m-1}\,\nu_0^\theta(\mathrm{d}x^m,\ldots,\mathrm{d}x^N) \\
&=\int_{\mathcal{K}^m}\int_{\mathcal{R}(x^m)}
f(x^{m-1},\ldots,x^N)\,\lambda^{m,1}(x^m,x^{m-1})\,\mathrm{d}x^{m-1}\,\nu_0(\mathrm{d}x^m,\ldots,\mathrm{d}x^N).
\end{split}
\end{equation}

\smallskip

We begin with the estimate 
\begin{equation*}
\begin{split}
&\bigg|\int_{\mathcal{K}^m} \int_{\mathcal{R}(x^m)}
f(x^{m-1},\ldots,x^N)\,\lambda^{m,\theta}(x^m,x^{m-1})\,\mathrm{d}x^{m-1}
\,\nu_0^\theta(\mathrm{d}x^m,\ldots,\mathrm{d}x^N) \\
& -\int_{\mathcal{K}^m}\int_{\mathcal{R}(x^m)}
f(x^{m-1},\ldots,x^N)\,\lambda^{m,1}(x^m,x^{m-1})\,\mathrm{d}x^{m-1}
\,\nu_0^\theta(\mathrm{d}x^m,\ldots,\mathrm{d}x^N)\bigg| \\
& \le C_f \int_{\mathcal{K}^m} \int_{\mathcal{R}(x^m)} \mathbf{1}_{\text{supp}(f)}(x^{m-1},\ldots,x^{N})\, |\lambda^{m,\theta}-\lambda^{m,1}|(x^m,x^{m-1})\,\mathrm{d}x^{m-1}
\,\nu_0^\theta(\mathrm{d}x^m,\ldots,\mathrm{d}x^N).
\end{split}
\end{equation*}
Let $K>1$ be such that $f(x^{m-1},\ldots,x^{N})=0$ if $x_1^N<-K$ or $x_N^N>K$, and let $\varepsilon\in(0,1)$ be such that $f(x^{m-1},\ldots,x^{N})=0$ if $x_i^{k-1} - x_i^{k} < \varepsilon$ or $x_{i+1}^k - x_i^{k-1} < \varepsilon$ for some $1 \leq i \leq k -1\leq N-1$, $k \geq m$. Notice that for any $n\in\N$, any $M>0$ and any $z_1,\ldots,z_n,w_1,\ldots,w_n\in[-M,M]$,
\[
\bigg|\prod_{i=1}^n z_i - \prod_{i=1}^n w_i\bigg| \leq M^{n-1}\sum_{i=1}^n |z_i - w_i|.
\]
Since $|x_i^{m-1}-x_j^m|^{\theta-1}\le 2K$ and $|x_j^m - x_i^m|^{1-2\theta} \le \varepsilon^{-2}$ on $\text{supp}(f)$ for all $\theta\ge1$ close enough to $1$, the quantity  $|\lambda^{m,\theta} - \lambda^{m,1}|$ is upper bounded on $\text{supp}(f)$ by
\begin{equation*}
\begin{split}
&\,M^{a_m} \bigg(\bigg|\frac{\Gamma(m\theta)}{\Gamma(\theta)^m}-\Gamma(m)\bigg|+
\sum_{i=1}^{m-1}\sum_{j=1}^m \big||x_i^{m-1}-x_j^{m}|^{\theta-1}-1\big| + \sum_{1 \leq i < j \leq m} \big|(x_j^m - x_i^m)^{1-2\theta}- (x_j^m - x_i^m)^{-1}\big|\bigg) \\
&\le M^{a_m} \bigg(\bigg|\frac{\Gamma(m\theta)}{\Gamma(\theta)^m}-\Gamma(m)\bigg|+\frac{3m(m-1)}{2\varepsilon}\,
\big(|\varepsilon^{\theta-1}\!-\!1| \vee |(2K)^{\theta-1}\!-\!1| \vee |\varepsilon^{2-2\theta}\!-\!1| \vee |(2K)^{2-2\theta}\!-\!1|\big)\!\bigg),
\end{split}
\end{equation*}
where $M := \max(2\Gamma(m),2K,\varepsilon^{-2})$, $a_m :=\frac{(m-1)(m-2)}{2}+ m(m-1) + \frac{m(m-1)}{2}$. This goes to $0$ as $\theta\downarrow1$.

\medskip

Moreover, the function $(x^m,\ldots,x^{N}) \mapsto  \int_{\mathcal{R}(x^m)}  f(x^{m-1},\ldots,x^{N})\,\lambda^{m,1}(x^m,x^{m-1})\,\mathrm{d}x^{m-1}$ belongs to $C_b(\mathcal{K}_m)$, so that the weak convergence $\nu^\theta_0\to\nu_0$ implies
\begin{equation}
\begin{split}
&\,\lim_{\theta \downarrow 1} \int_{\mathcal{K}^m}  \int_{\mathcal{R}(x^m)}  f(x^{m-1},\ldots,x^{N})\,\lambda^{m,1}(x^m,x^{m-1})\,\mathrm{d}x^{m-1}
\,\nu_0^\theta(\mathrm{d}x^m,\ldots,\mathrm{d}x^N) \\
&= \int_{\mathcal{K}^m} \int_{\mathcal{R}(x^m)}  f(x^{m-1},\ldots,x^{N})\,\lambda^{m,1}(x^m,x^{m-1})\,\mathrm{d}x^{m-1}
\,\nu_0(\mathrm{d}x^m,\ldots,\mathrm{d}x^N).
\end{split}
\end{equation}
Together with the findings of the preceding paragraph, this yields \eqref{lemma5.7claim}. \qed

\medskip

By Lemma \ref{tight}, the initial distributions
$\int_{\mathcal{K}^{k+1}} \mu_0^{\theta}(\mathrm{d}x^N)\,\prod_{m=k+1}^N\lambda^{m,\theta}(x^m,x^{m-1})\,\mathrm{d}x^{k+1}\,\ldots\,\mathrm{d}x^{N-1}$ converge weakly to $\int_{\mathcal{K}^{k+1}} \mu_0(\mathrm{d}x^N)\,\prod_{m=k+1}^N\lambda^{m,1}(x^m,x^{m-1})\,\mathrm{d}x^{k+1}\,\ldots\,\mathrm{d}x^{N-1}$ for all $k=1,\,2,\,\ldots,\,N$. Lemma \ref{DBMtight} then implies the process level convergence $X^{k,\theta}\to X^{k,1}$ for all $k=1,\,2,\,\ldots,\, N$. 
The tightness of $\{X^\theta\}_{\theta\downarrow 1}$ readily follows.

\medskip

\noindent\textbf{Step 3.} To identify the limit points of $\{X^\theta\}_{\theta\downarrow 1}$, we apply the Skorokhod Representation Theorem (see, e.g., \cite[Theorem 6.7]{billing}), which allows us to assume that the convergence $X^\theta\to X$ holds almost surely along a sequence of $\theta\downarrow1$. We now check that $X$ satisfies the properties (1)--(4) of Definition~\ref{submart}. Due to Lemma \ref{tight}, $X(0)$ is distributed according to $\mu_0\times\Lambda^{N,1}$, so property (1) holds. Property (2) is clear since $\mathbb{P}(X^{\theta}(t) \in \mathcal{K},\,t\ge0)=1 $ and $X^\theta\to X$ almost surely. For property (4), we recall from the preceding paragraph that each $X^k$ is a $1$-DBM. Thus, $\mathbb{P}(X^k(t) \in \,\stackrel{\circ}{\mathcal{W}}_k,\,t>0)=1$ (see, e.g., \cite[Proposition 4.3.5]{anderson_guionnet_zeitouni_2009}). This implies property (4) via the characterization of $\mathcal{V}$ in \eqref{Vchar}. In fact, $(\mu_0\times\Lambda^{N,1})(\partial\mathcal{K})=0$ yields $\mathbb{P}(X^k(t) \in \,\stackrel{\circ}{\mathcal{W}}_k,\,t\ge0)=1$, $k=1,\,2,\,\ldots,\,N$, and so 
\[
\tau_{\mathcal{V}} = \inf\{t \geq 0 : X(t) \in \mathcal{V} \} = \infty. 
\]

\smallskip

To see property (3), let $t>s\ge0$, let $F$ be a continuous bounded functional from $C\big([0,s],\mathbb{R}^{\frac{N(N+1)}{2}}\big)$ to $[0,\infty)$, and let $h \in \mathcal{H}$. Itô's Lemma reveals that
\begin{equation}\label{Step3Ito}
\begin{split}
&\,\mathbb{E}\bigg[F(X_{[0,s]}^{\theta})\,\bigg(h(X^{\theta}(t))-h(X^{\theta}(s)) -  \frac{1}{2}\int_s^t \Delta h(X^{\theta}(r))\,\mathrm{d}r\bigg)\bigg] \\
&= \mathbb{E}\bigg[F(X_{[0,s]}^{\theta})\sum_{1 \leq i \leq k \leq N}(\theta-1)\int_s^t \partial_{x_i^k}h(X^{\theta}(r))\,\widetilde{b}_i^k(X^{\theta}(r))\,\mathrm{d}r\bigg],
\end{split}
\end{equation}
where $(\theta-1)\widetilde{b}_i^k = b_i^{k,\theta}$. We decompose $\widetilde{b}_i^k$ into two parts:
\begin{equation*}
\begin{split}
\widetilde{b}_i^k(x) &= \sum_{j=1}^{k-1}\frac{1}{x_i^k - x_j^{k-1}} - \sum_{j\neq i} \frac{1}{x_i^k - x_j^k} \\
&= \sum_{j=1}^{i-1} \bigg(\frac{1}{x_i^k - x_j^{k-1}}-\frac{1}{x_i^k-x_j^k}\bigg) - \sum_{j=i}^{k-1}\bigg(\frac{1}{x_j^{k-1}-x_i^k}-\frac{1}{x_{j+1}^k - x_i^k} \bigg) 
=: \widetilde{b}_i^{k,+} + \widetilde{b}_i^{k,-}.
\end{split}
\end{equation*}
Note that $\widetilde{b}_i^{k,+}(x)\ge 0 \ge \widetilde{b}_i^{k,-}(x)$ for all $x \in\, \stackrel{\circ}{\mathcal{K}}$. Fix $1\leq i \leq k \leq N$ and $\delta > 0$. Then,
\begin{equation*}
\begin{split}
&\;\mathbb{E}\bigg[F(X_{[0,s]}^{\theta})\,(\theta-1)\int_s^t \partial_{x_i^k} h(X^{\theta}(r))\,\widetilde{b}_i^{k,+}(X^{\theta}(r))\,\mathrm{d}r\bigg] \\
&=\mathbb{E}\bigg[F(X_{[0,s]}^{\theta})\,(\theta-1)\int_{\{r \in [s,t]:\, |X_i^{k,\theta}(r)-X_{i-1}^{k-1,\theta}(r)|\ge \delta \}} \partial_{x_i^k}h(X^{\theta}(r)) \,\widetilde{b}_i^{k,+}(X^{\theta}(r))\,\mathrm{d}r\bigg] \\
&\quad +\mathbb{E}\bigg[F(X_{[0,s]}^{\theta})\,(\theta-1)\int_{\{r \in [s,t]:\, |X_i^{k,\theta}(r)-X_{i-1}^{k-1,\theta}(r)| < \delta \}} \partial_{x_i^k}h(X^{\theta}(r))\,\widetilde{b}_i^{k,+}(X^{\theta}(r))\,\mathrm{d}r\bigg].
\end{split}
\end{equation*}

\smallskip

The first term can be lower bounded by $-\frac{C_{F,h}(t-s)}{\delta}(\theta-1)$. Next, for $x \in \mathcal{K}$, we let $\Pi^+(x)$ be the projection of $x$ onto $\{x \in \mathcal{K}\!: x_i^{k} = x_{i-1}^{k-1} \}$. Note that if $|x_i^k-x_{i-1}^{k-1}|<\delta$, we only need to change a finite number of coordinates by at most $\delta$ in order for the resulting point to be in $\{x \in \mathcal{K}\!: x_i^{k} = x_{i-1}^{k-1} \}$. Therefore, $|x - \Pi^+(x)|\leq C\delta$ for such $x$.  Further, $h\in C_c^2(\mathcal{K})$ implies that the function $\partial_{x_i^k}h$ is Lipschitz. Combining Lemmas \ref{diagtoflat}, \ref{intertwinings} and \ref{DBMint}(a) we get for $\theta\in(1,2)$:
\begin{equation*}
\begin{split}
& \,\mathbb{E}\bigg[F(X_{[0,s]}^{\theta})\,(\theta-1)\int_{\{r \in [s,t]:\, |X_i^{k,\theta}(r)-X_{i-1}^{k-1,\theta}(r)|< \delta \}} \big|\partial_{x_i^k}h(X^{\theta}(r)) 
- \partial_{x_i^k}h(\Pi^{+}(X^{\theta}(r)))\big|
\,\widetilde{b}_i^{k,+}(X^{\theta}(r))\,\mathrm{d}r\bigg] \\
& \le C_{F,h}\,\delta \sup_{\theta\in(1,2)} \mathbb{E}\bigg[\int_s^t (\theta-1)\,\widetilde{b}_i^{k,+}(X^{\theta}(r))\,\mathrm{d}r\bigg] \leq C_{F,h}\,\delta.
\end{split}
\end{equation*}
Finally, since $h \in \mathcal{H}$ and $\widetilde{b}_i^{k,+} \geq 0$, we have that $\partial_{x_i^k}h(\Pi^{+}(X^\theta(r)))\,\widetilde{b}_i^{k,+}(X^\theta(r))\geq 0$. Thus, 
\[
\mathbb{E}\bigg[F(X_{[0,s]}^{\theta})\,(\theta-1)\sum_{1 \leq i \leq k \leq N}\int_s^t \partial_{x_i^k} h(X^{\theta}(r))\,\widetilde{b}_i^{k,+}(X^{\theta}(r))\,\mathrm{d}r\bigg]
\geq -C_{F,h}\,\delta - \frac{C_{F,h}\,(t-s)}{\delta}(\theta-1).
\]

\smallskip

By an analogous argument, one arrives at
\[
\mathbb{E}\bigg[F(X_{[0,s]}^{\theta})\,(\theta-1)\sum_{1 \leq i \leq k \leq N}\int_s^t \partial_{x_i^k} h(X^{\theta}(r))\,\widetilde{b}_i^{k,-}(X^{\theta}(r))\,\mathrm{d}r\bigg]
\geq -C_{F,h}\,\delta - \frac{C_{F,h}\,(t-s)}{\delta}(\theta-1).
\]
Plugging these into \eqref{Step3Ito}, we conclude 
\[
\mathbb{E}\bigg[F(X_{[0,s]}^{\theta})\,\bigg(h(X^{\theta}(t))-h(X^{\theta}(s)) -  \frac{1}{2}\int_s^t \Delta h(X^{\theta}(r))\,\mathrm{d}r\bigg)\bigg]
\geq -C_{F,h}\,\delta - \frac{C_{F,h}\,(t-s)}{\delta}(\theta-1).
\]
Passing to the limit $\theta\downarrow1$ then results in 
\[
\mathbb{E}\bigg[F(X_{[0,s]})\,\bigg(h(X(t))-h(X(s)) -  \frac{1}{2}\int_s^t \Delta h(X(r))\,\mathrm{d}r\bigg)\bigg]\geq -C_{F,h}\,\delta.
\]
Since $\delta>0$ was arbitrary, property (3) (``submartingale'' property) readily follows. 

\medskip

In view of Proposition \ref{correspondence} and $\tau_{\mathcal{V}}= \infty$, the canonical process under $\P\circ X^{-1}$, which we also label~$X$, solves the SDER associated with the Warren process $X^1$. In particular, there are $N(N+1)/2$-dimensional standard Brownian motion $B$ and process $Y$ such that $(X,Y)$ solves the ESP associated with $(\mathcal{K},\upsilon)$ and $Z(t):=X(0)+B(t)$, $t\ge0$ (see Definition \ref{def:SDER}(4)). We claim that $X$ actually solves the Skorokhod problem associated with $(\mathcal{K},\upsilon)$ and $Z$ by \cite[Theorem 2.9]{10.1214/EJP.v11-360}. Indeed, it is easy to see that $(\mathcal{K},\upsilon)$ satisfies \cite[Assumption 2.1]{10.1214/EJP.v11-360}; moreover, while $\mathcal{V}$ is defined differently in \cite[display~(2.15)]{10.1214/EJP.v11-360}, it clearly is a subset of our $\mathcal{V}$ (recall \eqref{Vchar}) for which $\tau_{\mathcal{V}}=\infty$. Since the Warren process is defined as the solution of the Skorokhod problem associated with $(\mathcal{K},\upsilon)$ and $Z$, the proof of Theorem~\ref{convergence} is complete. \qed

\medskip\smallskip

\bibliographystyle{plain}
\bibliography{bibliography.bib}

\smallskip

\appendix\section{Some Properties of Dyson Brownian Motion} \label{sec:app}

In this appendix, we collect the properties of $\theta$-DBM that are used throughout the paper. First, we recall that the $\theta$-DBM of dimension $N$ is the unique strong solution of the SDE \eqref{DBM} (see, e.g., \cite[Proposition 4.3.5]{anderson_guionnet_zeitouni_2009}). When $X(0)=0$, for each $t>0$, the random vector $X(t)$ admits a density $p_t$ on $\mathcal{W}$ which is given by  
\begin{equation}\label{eigv}
p_t(x) = \frac{1}{Z_{\theta,t}}\,\prod_{1\leq i < j \leq N} (x_j-x_i)^{2\theta}\,
\prod_{i=1}^N \exp{\bigg(\!-\frac{x_i^2}{2t}\bigg)},
\end{equation}
where $Z_{\theta,t}:=t^{\theta\frac{N(N-1)}{2}+\frac{N}{2}}(2\pi)^{\frac{N}{2}}
\prod_{i=1}^N\frac{\Gamma(i\theta)}{\Gamma(\theta)}$ (see, e.g., \cite[Proposition 2.9]{gorin2014multilevel}). When $X(0)\neq0$, we often apply the following comparison results to bound expectations with respect to $X$. 

\begin{lemma}\label{comparisons} 
Let $X$, $Y$ be two strong solutions of \eqref{DBM} with deterministic initial conditions $x$, $y$ in $\mathcal{W}$ and the same driving Brownian motion.
\begin{enumerate}[(a)]
\item Suppose that $y_i \le x_i$, $i=1,\,2,\,\ldots,\, N$. Then, almost surely, $Y_i(t)\le X_i(t)$, $t \ge 0$, $i=1,\,2,\,\ldots,\, N$.
\item Suppose that $y_j-y_i\le x_j-x_i$, $1\le i < j \le N$. Then, almost surely, $Y_j(t)-Y_i(t) \le X_j(t)-X_i(t)$, $t \ge 0$, $1\le i < j \le N$.
\item Almost surely, $|X(t)-Y(t)|_\infty \le 2N|x-y|_\infty$, and so $|X(t)-Y(t)|_\infty\wedge1 \le 2N(|x-y|_\infty\wedge 1)$, $t\ge0$.
\end{enumerate}
\end{lemma}

\noindent\textbf{Proof.} Parts (a) and (b) can be found in \cite[Lemma 2.10]{aggarwal2023edge}, so we only prove part (c). To this~end, we set $z^+_i = x_i + |x-y|_\infty$ and $z_i^- = x_i - |x-y|_\infty$. Then, $z_i^- \le x_i,y_i \le z_i^+$, $i=1,\,2,\,\ldots,\,N$. Let $Z^+$, $Z^-$ be strong solutions of \eqref{DBM} with the initial conditions $z^+$, $z^-$ driven by the same Brownian motion as $X$ and $Y$. By part (a), $Z_i^-(t)\le X_i(t),Y_i(t) \le Z_i^+(t)$, $t\ge0$, $i=1,\,2,\,\ldots,\,N$. Since $\P(Z^-(t),Z^+(t)\in\,\stackrel{\circ}{\mathcal{W}},\,t>0)=1$ and $\sum_{i=1}^N \sum_{j \neq i} \frac{1}{z_i-z_j}=0$, $z \in\,\stackrel{\circ}{\mathcal{W}}$,
\[
|X_j(t)-Y_j(t)| \le Z_j^{+}(t) - Z_j^-(t) \leq \sum _{i=1}^N (Z_i^+(t)-Z_i^-(t)) = \sum_{i=1}^N (Z_i^+(0)- Z_i^-(0)) = 2N |x-y|_{\infty},
\] 
$j=1,\,2,\,\ldots,\,N$, as desired. \qed

\medskip



Next, we give a bound on the $1$-Wasserstein distance between two strong solutions of \eqref{DBM}. 

\begin{lemma}\label{wass}
Let $X$ and $Y$ be two strong solutions of \eqref{DBM} with initial distributions $\mu_0$ and $\nu_0$, respectively, and let $\mu_t$ and $\nu_t$ be the laws of $X(t)$ and $Y(t)$, respectively. Writing $W_1(\cdot,\cdot)$ for the $1$-Wasserstein distance with respect to the metric $|x-y|_\infty\wedge 1$ on $\mathbb{R}^N$ we then have that
	\[
	W_1(\mu_t,\nu_t) \leq 2N W_1(\mu_0,\nu_0).
	\]
\end{lemma}

\noindent\textbf{Proof.} Consider any $(\mathcal{W}\times\mathcal{W})$-valued random vector $(\xi_1,\xi_2)$ with $\xi_1\sim\mu_0$ and $\xi_2\sim\nu_0$. Define $X^{\xi_1}$ and $X^{\xi_2}$ to be strong solutions of \eqref{DBM} with the initial conditions $\xi_1$ and $\xi_2$, respectively, driven by the same Brownian motion. Conditioning on $(\xi_1,\xi_2)$ and applying Lemma \ref{comparisons}(c) we find that almost surely, $\mathbb{E}[|X^{\xi_1}(t)-X^{\xi_2}(t)|_\infty\wedge 1\,|\,(\xi_1,\xi_2)] \le 2N(|\xi_1-\xi_2|_\infty\wedge1)$. Taking expectations and then the infimum over all couplings $(\xi_1,\xi_2)$ proves the lemma. \qed

\medskip

We repeatedly rely on the following inverse moment estimates for the $\theta$-DBM. 

\begin{lemma}\label{DBMint}
Let $\theta\ge\frac{1}{2}$ and let $X$ be a $\theta$-DBM. Then:
\begin{enumerate}[(a)]
\item For any $\theta_0 > \frac{1}{2}$ and any $p\in[1,2)$, there exists a $C_{\theta_0,p}<\infty$ such that as long as $\theta\in\big[\frac{1}{2},\theta_0\big]$,
\[
\mathbb{E}[(X_{i+1}(t)-X_i(t))^{-p}] \le C_{\theta_0,p}\, t^{-\frac{p}{2}},\quad t>0,\quad i=1,\,2,\,\ldots,\,N-1.
\]
\item For any $p\in[0,2\theta+1)$, there exists a $C_{\theta,p}<\infty$ such that
\[
\mathbb{E}[(X_{i+1}(t)-X_i(t))^{-p}] \leq C_{\theta,p}\, t^{-\frac{p}{2}},\quad t>0,\quad i=1,\,2,\,\ldots,\,N-1.
\]
\item For any $\alpha_1,\alpha_2 \in[0,2\theta+1)$, any $1 \le i \neq j \le N-1$, and any $t_2>t_1>0$, it holds
\[ 
\int_{t_1}^{t_2} \mathbb{E}[(X_{i+1}(t)-X_i(t))^{-\alpha_1}(X_{j+1}(t) - X_j(t))^{-\alpha_2}]\,\mathrm{d}t < \infty.
\]
\item If $\theta > \frac{1}{2}$ and for some $\varepsilon>0$, $\mathrm{dist}(X_0,\partial \mathcal{W})\ge\varepsilon$ with respect to $|\cdot|_\infty$ almost surely, then 
\[
\int_{0}^T \mathbb{E}[(X_{i+1}(t)-X_i(t))^{-2}]\,\mathrm{d}t < \infty,\quad T\in(0,\infty),\quad i=1,\,2,\,\ldots,\, N-1.
\]
\end{enumerate}
\end{lemma}

\noindent\textbf{Proof. (a), (b).} We start with part (a). By Lemma \ref{comparisons}(b), we may assume that $X(0)=0$. The normalization constant $Z_{\theta,t}$ in \eqref{eigv} is lower bounded by $\frac{1}{C_{\theta_0}}\, t^{\theta\frac{N(N-1)}{2}+\frac{N}{2}}$. Therefore,
\begin{equation}\label{bound1}
\begin{split}
& \mathbb{E}[(X_{i+1}(t)-X_i(t))^{-p}] \leq C_{\theta_0}\, t^{-\theta\frac{N(N-1)}{2}-\frac{N}{2}}\int_{\mathcal{W}} (x_{i+1}-x_i)^{-p} \prod_{1 \le j<m \le N} (x_m-x_j)^{2\theta} \,\prod_{j=1}^N \exp\bigg(\!-\frac{x_j^2}{2t}\bigg)\,\mathrm{d}x \\
&\qquad\qquad\qquad\qquad\quad\;\;\;\,
= C_{\theta_0}\, t^{-\frac{p}{2}}\int_{\mathcal{W}} (x_{i+1}-x_i)^{-p} \prod_{1 \le j < m \le N} (x_m-x_j)^{2\theta} \,\prod_{j=1}^N \exp\bigg(\!-\frac{x_j^2}{2}\bigg)\,\mathrm{d}x \\
&\qquad\qquad\qquad\qquad\quad\;\;\;\,
\le C_{\theta_0}\, t^{-\frac{p}{2}}\int_{\mathcal{W}} (x_{i+1}-x_i)^{2\theta-p}\, (x_N-x_1)^{\theta(N(N-1)-2)} \,\prod_{j=1}^N \exp\bigg(\!-\frac{x_j^2}{2}\bigg)\,\mathrm{d}x \\
& \qquad
=  C_{\theta_0}\, t^{-\frac{p}{2}}\int_{0}^{\infty}\ldots\,\int_0^{\infty}\int_{-\infty}^{\infty} y_{i+1}^{2\theta-p}\,\Big(\sum_{j=2}^N y_j\Big)^{\theta(N(N-1)-2)}\,
\prod_{j=1}^N \exp\bigg(\!-\frac{(\sum_{m=1}^j y_m)^2}{2}\bigg)\,\mathrm{d}y_1\,\mathrm{d}y_2\,\ldots\,\mathrm{d}y_N,
\end{split}
\end{equation}
where $y_1:=x_1$ and $y_j := x_j - x_{j-1}$, $j=2,\,3,\,\ldots,\, N$. 

\medskip

We aim to integrate in $y_1$ first. Note that
\begin{equation*}
\begin{split}
\sum_{j=1}^N \Big(\sum_{m=1}^j y_m\Big)^2 
& = Ny_1^2+2y_1\sum_{j=2}^N \sum_{m=2}^j y_m 
+ \sum_{j=2}^N \Big(\sum_{m=2}^j y_m\Big)^2 \\
& = N\bigg(y_1 + \frac{1}{N}\sum_{j=2}^N \sum_{m=2}^j y_m \bigg)^2 
+ \sum_{j=2}^N\Big(\sum_{m=2}^j y_m \Big)^2 
- \frac{1}{N}\Big(\sum_{j=2}^N \sum_{m=2}^j y_m\Big)^2 \\
& \ge N\bigg(y_1\!+\!\frac{1}{N}\sum_{j=2}^N \sum_{m=2}^j y_m \bigg)^2  
\!+\! \frac{1}{N} \sum_{j=2}^N\Big( \sum_{m=2}^j y_m\Big)^2 
\geq  N\bigg(y_1 \!+\! \frac{1}{N}\sum_{j=2}^N \sum_{m=2}^j y_m \bigg)^2  
\!+\! \frac{1}{N} \sum_{j=2}^N y_j^2.
\end{split}
\end{equation*}
Consequently,
\[
\int_{-\infty}^{\infty}\prod_{j=1}^N \exp\bigg(\!-\frac{(\sum_{m=1}^j y_m)^2}{2}\bigg)\,\mathrm{d}y_1 \leq C\prod_{j=2}^N \exp\bigg(\!-\frac{y_j^2}{2N} \bigg).
\]
So, the final expression in \eqref{bound1} can be upper bounded by
\[
C_{\theta_0}\, t^{-\frac{p}{2}} \sum_{j=2}^N \, \int_0^\infty \ldots\, \int_0^\infty y_{i+1}^{2\theta-p}\, y_j^{\theta(N(N-1)-2)}\,
\prod_{m=2}^N \exp\bigg(\!-\frac{y_m^2}{2N} \bigg)\,\mathrm{d}y_2\,\ldots\,\mathrm{d}y_N.
\]
Since $2\theta - p \ge 1 - p > -1$ and $\theta \le \theta_0$, the latter integral is bounded uniformly in $i=1,\,2,\,\ldots,\,N-1$ and $\theta\in\big[\frac{1}{2},\theta_0\big]$, yielding part (a).  Part (b) also readily follows.

\medskip

\noindent\textbf{(c).} To satisfy the assumptions of part (c), we must have $N \geq 3$. The same estimates as above yield
\begin{equation*}
\begin{split}
& \,\mathbb{E}[(X_{i+1}(t)-X_i(t))^{-\alpha_1}(X_{j+1}(t) - X_j(t))^{-\alpha_2}] \\
& \le C_\theta \, t^{-\frac{\alpha_1+\alpha_2}{2}}\sum_{m=2}^N \,
\int_0^\infty \ldots\, \int_0^\infty
y_{i+1}^{2\theta-\alpha_1}\, y_{j+1}^{2\theta-\alpha_2}\,
y_m^{\theta(N(N-1)-4)}\,
\prod_{n=2}^N \exp\bigg(\!-\frac{y_n^2}{2N} \bigg)\,\mathrm{d}y_2\,\ldots\,\mathrm{d}y_N.
\end{split}
\end{equation*}
Since $2\theta-\alpha_1,\,2\theta-\alpha_2 > -1$, the latter integral is finite. An integration in $t$ gives the result.

\medskip

\noindent\textbf{(d).} 
We start from
\[
\mathrm{d}(X_{i+1}(t)-X_i(t)) = \theta\Big(\sum_{j \neq i+1} \frac{1}{X_{i+1}(t)-X_j(t)} - \sum_{j\neq i} \frac{1}{X_i(t) - X_j(t)}\Big)\,\mathrm{d}t + \mathrm{d}(B_{i+1}(t)-B_i(t)),
\]
which, due to $X(s) \in \mathcal{W}$, $s \ge 0$, implies
\begin{equation*}
\begin{split}
&\,|(X_{i+1}(t) - X_i(t)) - (X_{i+1}(0)-X_i(0))| \\ 
&\le C_\theta\int_{0}^t \frac{1}{X_i(s)\!-\! X_{i-1}(s)}
+ \frac{1}{X_{i+1}(s)\!-\! X_i(s)} + \frac{1}{X_{i+2}(s)\!-\! X_{i+1}(s)}\,\mathrm{d}s 
+ |B_{i+1}(t)-B_{i}(t)| \\
& =: D_i(t) + |B_{i+1}(t)-B_{i}(t)|.
\end{split}
\end{equation*}
Next, we note that for the event
\[
A_t := \big\{\max_{1 \leq i \leq N-1} |(X_{i+1}(t) - X_i(t)) - (X_{i+1}(0)-X_i(0))| > \varepsilon \big\},
\]
it trivially holds
\[
\mathbb{P}(A_t) \leq \sum_{i=1}^{N-1} \Big(\mathbb{P}\big(D_i(t) > \varepsilon/2\big) 
+ \mathbb{P}\big(|B_{i+1}(t)-B_{i}(t)| > \varepsilon/2\big) \Big).
\]
By Chebyshev's inequality, the latter probability is at most $\frac{8t}{\varepsilon^2}$. For the former probability, we have 
\[
\mathbb{P}\big(D_i(t) > \varepsilon/2\big) \le \frac{2}{\varepsilon}\,\mathbb{E}[D_i(t)] 
\le \frac{2C_\theta}{\varepsilon} \int_{0}^t s^{-\frac{1}{2}}\,\mathrm{d}s 
= C_{\theta,\varepsilon}\,\sqrt{t}
\]
by part (b). Thus, $\mathbb{P}(A_t) \le C_{\theta,\varepsilon}\,\sqrt{t}$, $t\in(0,1]$.

\medskip

On the event $A_t^c$, $\mathrm{dist}(X(t),\partial \mathcal{W}) \ge \varepsilon/2$. Taking $t_0=T\wedge1$, applying part (b), letting $p\in(1,\theta+1/2)$ and $q^{-1} := 1-p^{-1}$, and using Hölder's inequality and part (b) again, we conclude that
\begin{equation*}\notag
\begin{split}
\int_0^T \mathbb{E}[(X_{i+1}(t)-X_i(t))^{-2}]\,\mathrm{d}t & = 
\int_0^{t_0} \mathbb{E}[(X_{i+1}(t)-X_i(t))^{-2}]\,\mathrm{d}t 
+ \int_{t_0}^T \mathbb{E}[(X_{i+1}(t)-X_i(t))^{-2}]\,\mathrm{d}t \\
&\le \int_0^{t_0} \mathbb{E}[(X_{i+1}(t)-X_i(t))^{-2}]\,\mathrm{d}t + C_{\theta,T} \\
&\le \int_0^{t_0}\mathbb{E}[(X_{i+1}(t)-X_i(t))^{-2}\,\mathbf{1}_{A_t}]\,\mathrm{d}t 
+ C_{\theta,T,\varepsilon} \\
&\le \int_0^{t_0} \mathbb{E}[(X_{i+1}(t)-X_i(t))^{-2p}]^{1/p}\,\mathbb{P}(A_t)^{1/q}\,\mathrm{d}t
+ C_{\theta,T,\varepsilon} \\
&\le C_{\theta,p,\varepsilon} \int_0^{t_0}t^{-1+\frac{1}{2q}}\,\mathrm{d}t 
+ C_{\theta,T,\varepsilon} < \infty. 
\qquad\qquad\qquad\qquad\qquad\qquad\quad\;\; \qed
\end{split}
\end{equation*}

The next lemma establishes the continuity of $\theta$-DBMs in $\theta$.

\begin{lemma}\label{DBMtight}
For $\theta\ge\frac{1}{2}$, let $X^\theta$ be a $\theta$-DBM with initial distribution $\mu_0^\theta$. If $\mu_0^\theta\to\mu_0^{\theta_0}$ weakly as $\theta\to\theta_0$, then $X^\theta$ converges in distribution to $X^{\theta_0}$, a $\theta_0$-DBM with initial distribution $\mu_0^{\theta_0}$, on $C([0,\infty),\mathbb{R}^N)$ with the topology of uniform convergence on compacts. 
\end{lemma}

\noindent\textbf{Proof.} Write $\mu^{\theta}$, $\mu^{\theta_0}$ for the laws of $X^\theta$, $X^{\theta_0}$, respectively, on $C([0,\infty),\mathbb{R}^N)$. It suffices to show that for every $T\in(0,\infty)$, the restriction $\mu_{[0,T]}^\theta$ of $\mu^\theta$ to $[0,T]$ converges weakly on $C([0,T],\mathbb{R}^N)$ to the restriction $\mu_{[0,T]}^{\theta_0}$ of $\mu^{\theta_0}$ to $[0,T]$. For the remainder of the proof, we fix such a $T$ and drop the subscripts $[0,T]$. To prove the tightness of $\{\mu^\theta\}_\theta$, we rely on \cite[Chapter 2, Theorem 4.10]{karatzas1991brownian}, and so need to verify 
\[
\lim_{K \to \infty} \, \sup_\theta\; \mu^\theta\big(\omega:\,|\omega(0)|_\infty > K\big) = 0
\]
and 
\[
\forall\,\varepsilon > 0:\quad \lim_{\delta \downarrow 0} \, \sup_\theta\; \mu^\theta\Big(\omega:\, \sup_{0 \le s,t \le T,\,|t-s|\le \delta}\,|\omega(t)-\omega(s)|_\infty > \varepsilon\Big)=0.
\]

\smallskip

The first claim is immediate since $\mu_0^\theta\to\mu_0^{\theta_0}$. For the second claim, denote by $(D^\theta_i)_{i=1}^N$ and $(c^\theta_i)_{i=1}^N$ the drift and the drift coefficient vectors of $X^\theta$. Then,
\begin{equation*}
\begin{split}
&\, \mu^\theta\Big(\omega:\, \sup_{0 \le s,t \le T,\,|t-s|\le \delta}\,|\omega(t)-\omega(s)|_\infty > \varepsilon\Big) \\
&\le  \sum_{i=1}^N \,\mathbb{P}\Big(\sup_{0 \le s,t \le T,\, |t-s|\le \delta} |D_i^\theta(t) - D_i^\theta(s)| > \frac{\varepsilon}{2}\Big) 
+ N\,\mathbb{P}\Big(\sup_{0 \le s,t \le T,\, |t-s|\le \delta} |B_1(t)-B_1(s)| > \frac{\varepsilon}{2}\Big).
\end{split}
\end{equation*}
The second term is independent of $\theta$ and goes to $0$ with $\delta$. Next, we fix a $p\in(1,2)$, and apply Lemma~\ref{DBMint}(a) to find for $\theta$ within $1$ of $\theta_0$,
\[
\mathbb{E}\bigg[\int_{0}^T |c_i^\theta(X^\theta(t))|^p \,\mathrm{d}t\bigg] 
\le C_{\theta_0,p}\,\mathbb{E}\bigg[\int_{0}^T (X_{i+1}^\theta(t)-X_i^\theta(t))^{-p} + (X_{i}^\theta(t)-X_{i-1}^\theta(t))^{-p}\,\mathrm{d}t\bigg] 
\le C_{\theta_0,p}\, T^{1-\frac{p}{2}} < \infty.
\]
Let $q^{-1}:=1-p^{-1}$. Then, using Hölder's inequality, the latter estimate, and Jensen's inequality,
\begin{equation*}
\begin{split}
&\mathbb{P}\Big(\sup_{0 \le s,t \le T,\,|t-s|\le \delta} |D_i^\theta(t) - D_i^\theta(s)| > \frac{\varepsilon}{2}\Big) 
\le \frac{2}{\varepsilon}\,\mathbb{E}\bigg[\sup_{0 \le s,t \le T,\, |t-s|\le \delta} 
\int_s^t |c_i^\theta(X^\theta(r))|\,\mathrm{d}r \bigg] \\
&\qquad\qquad\qquad\qquad\qquad
\le \frac{2}{\varepsilon}\,\mathbb{E}\bigg[\sup_{0 \le s,t \le T,\, |t-s|\le \delta}  
\bigg(\int_0^T |c_i^\theta(X^\theta(r))|^p \,\mathrm{d}r\bigg)^{1/p}\,
\bigg(\int_s^t 1 \,\mathrm{d}r \bigg)^{1/q} \bigg] 
\le \frac{C_{\theta_0,p,T}}{\varepsilon}\,\delta^{1/q}.
\end{split}
\end{equation*}
This establishes the tightness of $\{\mu^\theta\}_\theta$.

\medskip

It remains to check that any limit point of $\{\mu^\theta\}_\theta$ is equal to $\mu^{\theta_0}$. By the Skorokhod Representation Theorem (see, e.g., \cite[Theorem 6.7]{billing}), we may assume that, along the sequence in consideration, $X^\theta\to X$ almost surely, for some $X$. Let us show that $X$ solves the martingale problem associated with the $\theta_0$-DBM. Let $0 \le s < t \le T$, let $F\!:C([0,s],\mathbb{R}^N)\rightarrow \mathbb{R}$ be a bounded continuous functional, and let $f \in C_{c}^{\infty}(\mathbb{R}^N)$. We claim that 
\begin{align}\label{DBMMART}
\mathbb{E}\bigg[F(X_{[0,s]})\bigg(f(X(t))-f(X(s)) - \int_s^t \frac{1}{2}\Delta f(X(r))
+ \sum_{i=1}^N c_i^{\theta_0}(X(r))\,\partial_{x_i}f(X(r))\,\mathrm{d}r\bigg)\bigg]=0.
\end{align}
Implicit in the above statement is that 
\begin{equation}\label{integ}
\mathbb{E}\bigg[\int_{0}^T |c_i^{\theta_0}(X(r))| \,\mathrm{d}r\bigg]< \infty.
\end{equation}
To see \eqref{integ}, note that the summands in $c_i^{\theta_0}(X(r))$ are the almost sure limits of the summands in $c_i^{\theta_0}(X^\theta(r))$. Lemma \ref{DBMint}(a) with $p\in(1,2)$ renders the Vitali Convergence Theorem applicable, so the $L^1$-norms of the summands in $c_i^{\theta_0}(X(r))$ are at most $C_{\theta_0}\,r^{-\frac{1}{2}}$. This yields \eqref{integ}.  

\medskip

To obtain \eqref{DBMMART}, we begin with 
\[
\mathbb{E}\bigg[F(X_{[0,s]}^\theta)\bigg(f(X^\theta(t))-f(X^\theta(s)) 
- \int_s^t \frac{1}{2}\Delta f(X^\theta(r)) + \sum_{i=1}^N c_i^\theta(X^\theta(r))\, \partial_{x_i}f(X^\theta(r))\,\mathrm{d}r\bigg)\bigg]=0.
\]
Clearly, it holds 
\begin{eqnarray*}
&& \lim_{\theta\to\theta_0}\, \Big|\mathbb{E}\Big[F(X_{[0,s]}^\theta)\,\big(f(X^\theta(t))-f(X^\theta(s))\big) -F(X_{[0,s]})\,\big(f(X(t))-f(X(s))\big)\Big]\Big| = 0, \\
&& \lim_{\theta\to\theta_0}\,
\bigg|\mathbb{E}\bigg[F(X_{[0,s]}^\theta)\int_{s}^t \frac{1}{2}\Delta f(X^\theta(r))\,\mathrm{d}r 
-  F(X_{[0,s]})\int_{s}^t \frac{1}{2}\Delta f(X(r))\,\mathrm{d}r\bigg]\bigg| = 0.
\end{eqnarray*}
Finally, we estimate
\[
\bigg|\mathbb{E}\bigg[F(X_{[0,s]}^\theta) 
\int_s^t c_i^\theta(X^\theta(r))\,\partial_{x_i}f(X^\theta(r))\,\mathrm{d}r 
- F(X_{[0,s]})\int_s^t c_i^{\theta_0}(X(r))\,\partial_{x_i}f(X(r))\,\mathrm{d}r\bigg]\bigg|
\]
\begin{equation*}
\begin{split}
&\le \mathbb{E}\bigg[\big|F(X_{[0,s]}^\theta)\big| \int_s^t \big|\partial_{x_i}f(X^\theta(r))\big|\, \big|c_i^\theta(X^\theta(r))-c_i^{\theta_0}(X^\theta(r))\big|\,\mathrm{d}r\bigg] 
\qquad\qquad\qquad\qquad\quad \\
&\quad
+\mathbb{E}\bigg[\big|F(X_{[0,s]}^\theta)\big| \int_s^t \big|c_i^{\theta_0}(X^\theta(r))\big|\,
\big|\partial_{x_i}f(X^\theta(r))-\partial_{x_i}f(X(r))\big|\,\mathrm{d}r\bigg] \\
&\quad
+\mathbb{E}\bigg[\big|F(X_{[0,s]}^\theta)\big| \int_s^t \big|\partial_{x_i}f(X(r))\big|\,
\big|c_i^{\theta_0}(X^\theta(r)) - c_i^{\theta_0}(X(r))\big|\,\mathrm{d}r \bigg] \\
&\quad
+ \mathbb{E}\bigg[\big|F(X_{[0,s]}^\theta) - F(X_{[0,s]})\big| 
\int_s^t \big|\partial_{x_i}f(X(r))\big|\,\big|c_i^{\theta_0}(X(r))\big|\,\mathrm{d}r\bigg].
\end{split}
\end{equation*}

\smallskip

In the first summand, we bound $|F|$ and $|\partial_{x_i}f|$ by a constant, and use Lemma \ref{DBMint}(a) to find that
\[
\mathbb{E}\Big[\big|c_i^\theta(X^\theta(r))- c_i^{\theta_0}(X^\theta(r))\big|\Big] 
= |\theta - \theta_0|\,\mathbb{E}\bigg[\bigg|\sum_{j \neq i} \frac{1}{X_i^\theta(r)-X_j^\theta(r)}\bigg|\bigg]
\le C_{\theta_0}\,|\theta-\theta_0|\,r^{-\frac{1}{2}}.
\]
In the second summand, $|F|$ can be bounded by a constant, and we can apply Hölder's inequality with $p\in(1,2)$ and $q^{-1}:=1-p^{-1}$ and then Lemma \ref{DBMint}(a) to obtain
\[
C_{\theta_0,p} \int_{s}^t r^{-\frac{1}{2}}\,
\mathbb{E}\Big[\big|\partial_{x_i}f(X^\theta(r)) - \partial_{x_i}f(X(r))\big|^q\Big]^{\frac{1}{q}}\,\mathrm{d}r,
\]
which goes to $0$ as $\theta\to\theta_0$ by the Bounded Convergence Theorem. In the third summand, we bound $|F|$ and $|\partial_{x_i}f|$ by a constant, and recall that $c_i^{\theta_0}(X^\theta(r))\to c_i^{\theta_0}(X(r))$ in $L^1$ for any fixed $r\in(0,T]$. Since $\mathbb{E}[|c_i^{\theta_0}(X^\theta(r))-c_i^{\theta_0}(X(r))|]\le C_{\theta_0}\,r^{-\frac{1}{2}}$ by Lemma \ref{DBMint}(a), the third summand goes to $0$ as $\theta\to\theta_0$ by the Dominated Convergence Theorem. In the fourth summand, we bound $|\partial_{x_i}f|$ by a constant, and use Hölder's inequality with $p\in(1,2)$ and $q^{-1}:=1-p^{-1}$ followed by Jensen's inequality and Lemma \ref{DBMint}(a) to get 
\begin{equation*}
\begin{split}
&\,C\,\mathbb{E}\bigg[\bigg(\int_s^t \big|c_i^{\theta_0}(X(r))\big|\,\mathrm{d}r\bigg)^p\bigg]^{\frac{1}{p}}\, \mathbb{E}\Big[\big|F(X_{[0,s]}^\theta) - F(X_{[0,s]})\big|^q\Big]^{\frac{1}{q}} \\
&\le C(t-s)^{1-\frac{1}{p}}\,\mathbb{E}\bigg[\int_s^t \big|c_i^{\theta_0}(X(r))\big|^p \,\mathrm{d}r\bigg]^{\frac{1}{p}}\,
\mathbb{E}\Big[\big|F(X_{[0,s]}^\theta) - F(X_{[0,s]})\big|^q\Big]^{\frac{1}{q}} \\
&\leq C(t-s)^{1-\frac{1}{p}}\,\bigg(\int_s^t r^{-\frac{p}{2}}\,\mathrm{d}r\bigg)^{\frac{1}{p}}\,
\mathbb{E}\Big[\big|F(X_{[0,s]}^\theta) - F(X_{[0,s]})\big|^q\Big]^{\frac{1}{q}},
\end{split}
\end{equation*}
which goes to $0$ as $\theta\to\theta_0$. 

\medskip

Altogether, we have shown \eqref{DBMMART}. Thus, for every $f \in C_c^{\infty}(\mathbb{R}^N)$, the process 
\[
f(X(t)) - f(X(0)) - \int_0^t \frac{1}{2}\Delta f(X(r)) + 
\sum_{i=1}^Nc_i^{\theta_0}(X(r))\,\partial_{x_i}f(X(r))\,\mathrm{d}r 
\]
is a martingale.~A standard localization argument renders $X$ a solution of the local martingale problem, and hence a weak solution of the SDE \eqref{DBM} (see, e.g., \cite[Chapter 5, Proposition 4.6]{karatzas1991brownian}). Such a weak solution is unique, so the law of $X$ is $\mu^{\theta_0}$. \qed

\medskip

We also need some regularity of the densities of $\theta$-DBMs at fixed times, which is the subject of the next lemma. 

\begin{lemma}\label{reg} 
Fix a $\theta>\frac{1}{2}$, and let $X$ be a $\theta$-DBM with a compactly supported initial distribution~$\mu_0$ such that $\mathrm{dist}(\mathrm{supp}(\mu_0),\partial \mathcal{W})\ge\varepsilon$ for some $\varepsilon>0$ (with respect to $|\cdot|_\infty$, say). If $\mu_0(\mathrm{d}x)=p_0\,\mathrm{d}x$ and $\int_{\mathcal{W}} | \log p_0|\, p_0\,\mathrm{d}x< \infty$, then necessarily  $\mathbb{P}(X(t) \in \mathrm{d}x) = p_t \,\mathrm{d}x$ with $p_t \in W_1^{1}(\mathbb{R}^N)$ for all $t>0$. Moreover, for any $T\in(0,\infty)$, $p \in L^{1,\frac{N}{N-2}}([0,T]\times\mathbb{R}^N)$ if $N > 2$, and $p \in L^{1,q}([0,T]\times\mathbb{R}^N)$ for all $q\in[1,\infty)$ if $N=2$. In both cases,
\begin{equation}
\int_{0}^T\int_{\mathcal{W}}\frac{|\nabla p_t|^2}{p_t}\,\mathrm{d}x \,\mathrm{d}t < \infty.
\end{equation}
\end{lemma}

\noindent\textbf{Proof.} We aim to apply \cite[Theorem 7.4.1]{bogachev2015fokker}. Write $c = (c_i)_{i=1}^N$ and $\{\mu_t\}_{t\ge0}$ for the drift coefficient vector and the fixed time distributions of $X$, respectively. For $f\in C_c^\infty((0,T)\times\stackrel{\circ}{\mathcal{W}})$, 
\[ 
\int_{0}^T \int_{\mathcal{W}} \partial_t f + \frac{1}{2}\Delta f + c \cdot \nabla f \,\mathrm{d}\mu_t \,\mathrm{d}t = 0
\]
by Itô's formula. The other assumption in \cite[Theorem 7.4.1]{bogachev2015fokker} is that, with the measure $\mathrm{d}\mu:=\mathrm{d}\mu_t \,\mathrm{d}t$ on $(0,T)\times\stackrel{\circ}{\mathcal{W}}$, it holds $|c|,\log \max(|x|,1) \in L^2(\mu)$. In our case,
\[
\int_{(0,T)\times\stackrel{\circ}{\mathcal{W}}} |c|^2 \,\mathrm{d}\mu 
\le C\int_{(0,T)\times\stackrel{\circ}{\mathcal{W}}} \sum_{i =1}^N \sum_{j \neq i} \frac{1}{(x_i-x_j)^2} \,\mathrm{d}\mu 
= C\sum_{i =1}^N \sum_{j \neq i}\int_0^T \mathbb{E}\big[(X_i(t)-X_j(t))^{-2}\big]\,\mathrm{d}t < \infty
\]
by Lemma \ref{DBMint}(d). We turn to $\log \max(|x|,1)$ and notice that $(\log \max(|x|,1))^2 \le |x|\leq \sum_{i=1}^N |x_i|$. Furthermore,
\[
|X_i(t)| \leq |X_i(0)| + \theta\int_0^t\sum_{j \neq i} \frac{1}{|X_i(s)-X_j(s)|} \,\mathrm{d}s + |B_i(t)|.
\]
Therefore, the compactness of $\mathrm{supp}(\mu_0)$ and Lemma \ref{DBMint}(b) imply
\[
\mathbb{E}[|X_i(t)|]\leq C_{\mu_0} + Ct^{\frac{1}{2}} + \sqrt{\frac{2t}{\pi}}.
\]
Integrating in $t$ we obtain $\log \max(|x|,1) \in L^2(\mu)$, and \cite[Theorem 7.4.1]{bogachev2015fokker} yields the lemma. \qed

\medskip

We conclude with an important consequence of Lemma \ref{reg}. 

\begin{lemma}\label{TIBP}
In the setting of Lemma \ref{reg}, 
\begin{equation} \label{FPmodified}
\int_{0}^T \int_{\mathcal{W}} (\partial_t f + c \cdot \nabla_x f)\, p_t - \frac{1}{2}\nabla_x f \cdot \nabla p_t \,\mathrm{d}x \,\mathrm{d}t = 0,\quad f \in C_c^{1,1}((0,T)\,\times\stackrel{\circ}{\mathcal W}). 
\end{equation}
\end{lemma}

\noindent\textbf{Proof.} Equation \eqref{FPmodified} is clear when $f \in C_c^{1,2}((0,T)\,\times\stackrel{\circ}{\mathcal W})$, for which one can apply Itô's formula and then integrate by parts in $x$. Hereby, we note that the integrability of $\nabla_x f \cdot \nabla p_t$ follows from
\[
\int_0^T\int_{\mathcal{W}} |\nabla_x f \cdot \nabla p_t| \,\mathrm{d}x \,\mathrm{d}t 
\le \bigg(\int_0^T \int_{\mathcal{W}} |\nabla_x f|^2\, p_t \, \mathrm{d}x \,\mathrm{d}t \bigg)^{1/2}\, 
\bigg(\int_0^T \int_{\mathcal{W}}\frac{|\nabla p_t|^2}{p_t} \,\mathrm{d}x \,\mathrm{d}t \bigg)^{1/2}.
\]
For $f \in C_c^{1,1}((0,T)\,\times\!\stackrel{\circ}{\mathcal W})$, let $\{\zeta_\varepsilon\}_{\varepsilon>0}$ be scaled standard bump functions in $\mathbb{R}^N$, set $f_\varepsilon=f\ast\zeta_\varepsilon$, and observe that $\bigcup_{t \in [0,T]} \mathrm{supp}(f_\varepsilon(t,\cdot))$ is contained in a compact $K\subset\,\stackrel{\circ}{\mathcal W}$ for all $\varepsilon>0$ small enough. From here on, we only consider $N>2$, with the case $N=2$ being similar. Take $0<a<b<T$ satisfying $f(t,\cdot)\equiv0$ unless $t\in[a,b]$. Since \eqref{FPmodified} holds for $f_\varepsilon$, it suffices to check that the following two terms go to $0$ with $\varepsilon$:
\begin{eqnarray}
&& \int_{0}^T\int_{\mathcal{W}} |\partial_t(f - f_{\varepsilon}) + c\cdot \nabla_x(f - f_{\varepsilon})|\, p_t \,\mathrm{d}x \,\mathrm{d}t,  \label{term1} \\
&& \int_{0}^T\int_{\mathcal{W}}|\nabla_x (f-f_{\varepsilon}) \cdot \nabla p_t| \,\mathrm{d}x \,\mathrm{d}t. \label{term2}
\end{eqnarray}

\smallskip

Using Hölder's inequality twice we bound the term in \eqref{term1} by
\begin{equation*}
\begin{split}
&\,\sup_{t\in[a,b]} \bigg(\int_K |\partial_t(f-f_\varepsilon)|^{\frac{N}{2}}\,\mathrm{d}x\bigg)^{2/N}\,
\|p\|_{L^{1,\frac{N}{N-2}}([a,b]\times K)} \\
&+C\,\sup_{t\in[a,b]} \bigg(\int_K |\nabla_x(f-f_\varepsilon)|^{\frac{N}{2}}\,\mathrm{d}x\bigg)^{2/N}\,
\|p\|_{L^{1,\frac{N}{N-2}}([a,b]\times K)}.
\end{split}
\end{equation*}
For the term in \eqref{term2}, we first observe that 
\[
\int_a^b \int_K |\nabla_x (f-f_{\varepsilon}) \cdot \nabla p_t| \,\mathrm{d}x \,\mathrm{d}t 
\le \bigg(\int_a^b \int_K |\nabla_x(f-f_\varepsilon)|^2\, p_t \,\mathrm{d}x \,\mathrm{d}t\bigg)^{1/2}\, 
\bigg(\int_{a}^b \int_K \frac{|\nabla p_t|^2}{p_t}\,\mathrm{d}x \,\mathrm{d}t\bigg)^{1/2},
\]
and then apply Hölder's inequality again to estimate 
\[
\int_a^b \int_K |\nabla_x(f-f_\varepsilon)|^2\, p_t \,\mathrm{d}x \,\mathrm{d}t
\le \sup_{t\in[a,b]} \bigg(\int_K |\nabla_x(f-f_\varepsilon)|^N \,\mathrm{d}x\bigg)^{2/N}\,
\|p\|_{L^{1,\frac{N}{N-2}}([a,b]\times K)}.
\]
Since $(\partial_t f_\varepsilon,\nabla_x f_\varepsilon)\to(\partial_t f,\nabla_x f)$ uniformly on $[a,b]\times K$ as $\varepsilon\downarrow0$, the obtained bounds on the terms in \eqref{term1} and \eqref{term2} go to $0$ with $\varepsilon$. The extension of \eqref{FPmodified} to all $f \in C_c^{1,1}((0,T)\,\times\stackrel{\circ}{\mathcal W})$ is complete. \qed

\bigskip\bigskip

\end{document}